# THE SCALING LIMITS OF PLANAR LERW IN FINITELY CONNECTED DOMAINS

By Dapeng Zhan

*University of California, Berkeley*

We define a family of stochastic Loewner evolution-type processes in finitely connected domains, which are called continuous LERW (loop-erased random walk). A continuous LERW describes a random curve in a finitely connected domain that starts from a prime end and ends at a certain target set, which could be an interior point, or a prime end, or a side arc. It is defined using the usual chordal Loewner equation with the driving function being $\sqrt{2}B(t)$ plus a drift term. The distributions of continuous LERW are conformally invariant. A continuous LERW preserves a family of local martingales, which are composed of generalized Poisson kernels, normalized by their behaviors near the target set. These local martingales resemble the discrete martingales preserved by the corresponding LERW on the discrete approximation of the domain. For all kinds of targets, if the domain satisfies certain boundary conditions, we use these martingales to prove that when the mesh of the discrete approximation is small enough, the continuous LERW and the corresponding discrete LERW can be coupled together, such that after suitable reparametrization, with probability close to 1, the two curves are uniformly close to each other.

**1. Introduction.** LERW (loop-erased random walk) (cf. [4]) is obtained by removing loops, in the order they are created, from a simple random walk on a graph that is stopped at some hitting time. Since the loops are erased, so an LERW is a simple lattice path. In this paper, we will consider the loop-erasures of conditional random walks. They have properties that are very similar to loop-erased random walks, so we still call them LERW.

In [16], Schramm introduced stochastic Loewner evolution (SLE), a family of random growth processes of closed fractal subsets in simply connected plane domains. The evolution is described by the classical Loewner equation with the driving term being $\sqrt{\kappa}$ times a standard linear Brownian motion for some $\kappa \geq 0$. SLE behaves differently for different values of $\kappa$. Schramm conjectured that $\mathrm{SLE}_2$ is the scaling limit of a kind of LERW on the grid approximation of the domain. And he proved the conjecture in that paper under the assumption that the scaling limits of LERW are conformally invariant.

Schramm's processes turned out to be very useful. On the one hand, they are amenable to computations; on the other hand, they are related with some statistical physics models. In a series of papers [6, 7, 8], Lawler, Schramm and Werner used SLE to determine the Brownian motion intersection exponents in the plane. In [10], the conjecture in [16] is completely solved, where no additional assumption is added. In the same paper, $\mathrm{SLE}_8$ is proved to be the scaling limits of UST (uniform spanning tree) Peano curve. Smirnov proved in [18] that chordal $\mathrm{SLE}_6$ is the scaling limit of critical site percolation on the triangular lattice. And Schramm and Sheffield proved in [17] that the harmonic explorer converges to chordal $\mathrm{SLE}_4$. In [9], $\mathrm{SLE}_{8/3}$ is proved to have the restriction property, and so is conjectured to be the scaling limits of self-avoiding walk. For the properties of SLE, see [5, 15] and [19].

At the beginning, the SLE is only defined in simply connected domains, because the definition uses the Riemann mapping theorem. In [20], a kind of SLE-type process, which is called annulus SLE, is defined in doubly connected domains. The definition uses the rotation symmetry and reflection symmetry of an annulus. It is proved there that annulus $\mathrm{SLE}_2$ is the scaling limit of the LERW in the grid approximation of a doubly connected domain that starts from a vertex that is close to a boundary point and stops when it hits the other boundary component.

The definitions of LERW on grid approximations of simply or doubly connected domains could be easily extended to multiply connected domains. It is interesting to study the scaling limits of the LERW in multiply connected domains. This may help us to extend the SLE to multiply connected domains.

In this paper, we will define a family of SLE-type processes, which are called continuous LERW, in finitely connected domains. They are defined using the usual chordal Loewner equation with the driving function being $\sqrt{2}B(t) + S(t)$, where $B(t)$ is a standard linear Brownian motion, and the drift term $S(t)$ is continuously differentiable in $t$. The drift term is carefully chosen, so that the continuous LERW satisfy the conformal invariance, and preserve a family of local martingales generated by generalized Poisson kernels. The local martingales resemble the discrete martingales preserved by the corresponding discrete LERW on the discrete approximation of that domain. And this resemblance is used to prove the convergence of discrete LERW to continuous LERW.

This paper is organized as follows. In Section 2, we define some notation that will be used in this paper. In Section 3, three kinds of continuous LERW are defined, which are continuous LERW "aimed" at interior points, prime ends and side arcs. And we prove that they all satisfy the conformal invariance. In Section 4, we present the continuous and discrete martingales preserved by continuous and discrete LERW, respectively, and explain the similarity between these martingales.

In Section 5, we give a rigorous proof of the existence and uniqueness of the solution to the equation that is used to define a continuous LERW. The lemmas that are used for the proof are interesting. We first use the idea of Carathéodory topology to define the convergence of plane domains. Then we define a metric on the space of hulls in the upper half plane, so that the set of hulls that are contained in a fixed hull is compact. This compactness property is frequently used in the remaining part of this paper. In this section, we use it to derive many uniform constants without working on concrete functions.

In Section 6, we first consider one kind of LERW, whose targets are interior points. The method given in [10] is used to get a coupling of the driving process for the discrete LERW and that for the continuous LERW such that the two driving processes are uniformly close to each other in probability. In Section 7, we first use some regular properties of the discrete LERW curve to get a local coupling of the LERW curve and the continuous LERW trace so that the two curves are close to each other, before either of them leaves a hull bounded by a crosscut. Finally, we glue all local couplings to get a global coupling of the curves. In the last section, we study the convergence of the other two kinds of LERW. And we get the similar results of the convergence.

## 2. Some notation.

2.1. *Loop-erased random walk.* In general, an LERW is defined on a connected locally finite graph $G=(V,E)$. We will usually consider the graphs that are discrete approximations of some plane domains. A loop-erasure of a finite lattice path $v=(v(0),\dots,v(n))$ on $G$ is defined as follows. Let $n_0=\max\{m: v(m)=v(0)\}$. Define the sequence $(n_j)$ inductively by $n_{j+1}=\max\{m: v(m)=v(n_j+1)\}$ if $n_j$ is defined and $n_j<n$. Let $\chi$ be the first $j$ such that $n_j=n$. Let $w(j)=v(n_j)$ for $0\le j\le \chi$. Then $w=(w(0),\dots,w(\chi))$ is called the loop-erasure of $(v(0),\dots,v(n))$ (see [4]), and is denoted by $\mathrm{LE}(v)$. It is a simple lattice path with $w(0)=v(0)$ and $w(\chi)=v(n)$.

A subset $S$ of $V$ is called reachable in $G$ if for any $v\in V\setminus S$, a (simple) random walk on $G$ started from $v$ will hit $S$ in finitely many steps almost surely. Suppose $A$ and $B$ are disjoint subsets of $V$ such that $A\cup B$ is reachable in $G$. Suppose $v_0\in V\setminus(A\cup B)$ and there is a lattice path on $G$ connecting $v_0$ and $A$ without passing through $B$. Then the probability that a random walk started from $v_0$ hits $A$ before $B$ is positive. We now consider this random walk stopped on hitting $A \cup B$ and conditioned to hit $A$. It is a random finite lattice path. The loop-erasure of this path is called the LERW on $G$ started from $x$ conditioned to hit $A$ before $B$.

For a function $f$ defined on $V$, and $v \in V$, let $\Delta_G f(v) = \sum_{w\sim v}(f(w) - f(v))$, where $w \sim v$ means that $w$ and $v$ are adjacent. If $\Delta_G f(v) = 0$, then we say $f$ is discrete harmonic at $v$. The proof of the following lemma is easy, and can be found in [20].

Lemma 2.1. *Suppose $A$ and $B$ are disjoint subsets of $V$ and $A \cup B$ is reachable in $G$. Let $x \in V \setminus (A \cup B)$ be such that there is a lattice path connecting $x$ and $A$ without passing through any vertex on $B$. Then there is a unique nonnegative bounded function $h$ on $V$ such that $h \equiv 0$ on $A \cup B$; $\Delta_G h \equiv 0$ on $V \setminus (A \cup B \cup \{x\})$; and $\sum_{v\in A} \Delta_G h(v) = 1$. Moreover, if either $A$ or $B$ is a finite set, then there is a unique nonnegative bounded function $g$ on $V$ such that $g \equiv 0$ on $B$; $g \equiv 1$ on $A$; $\Delta_G g \equiv 0$ on $V \setminus (A \cup B \cup \{x\})$; and $\sum_{v\in A} \Delta_G g(v) = 0$.*

Suppose $E_{-1}$ and $F$ are disjoint subsets of $V$ and $E_{-1} \cup F$ is reachable in $G$. Let $x_0 \in N$ be such that there is a lattice path connecting $x_0$ and $F$ without passing through any vertex on $E_{-1}$. Let $(q(0), \dots, q(\chi))$ be the LERW on $G$ started from $x_0$ conditioned to hit $F$ before $E_{-1}$. So $q(0) = x_0$ and $q(\chi) \in F$. For $0 \le j < \chi$, let $E_j = E_{-1} \cup \{q(0), \dots, q(j)\}$. Then $E_j$ and $F$ are disjoint. Since $E_j \cup F$ is bigger than $E_{-1} \cup F$, so it is also reachable. Note that for any $0 \le j < \chi$, $(q(j), \dots, q(\chi))$ is a lattice path connecting $q(j)$ with $F$ without passing through $E_{j-1}$. Let $h_j$ be as in Lemma 2.1 with $A = F$, $B = E_{j-1}$ and $x = q(j)$. If either $E_{-1}$ or $F$ is finite, then either $E_j$ or $F$ is finite. Let $g_j$ be the $g$ in Lemma 2.1 with $A = F$, $B = E_{j-1}$ and $x = q(j)$. Let $\overline{F}$ be the union of $F$ with the set of vertices of $V$ that are adjacent to $F$. Then we have:

Proposition 2.1. *Fix any $v_0 \in V$. Then $(g_k(v_0))$ (if $E_{-1}$ or $F$ is finite) and $(h_k(v_0))$ are discrete martingales up to the first time $x_k$ hits $\overline{F}$, or $E_k$ disconnects $v_0$ from $F$ in $G$.*

Proof. The result for $(g_k)$ in a special case is Proposition 3.2 in [20]. The proof of that proposition applies to general cases. The proof for $(h_k)$ is similar. □

2.2. *Finitely connected domains.* In this paper, a domain is a nonempty connected open subset of the Riemann sphere $\widehat{\mathbb{C}} = \mathbb{C} \cup \{\infty\}$. Here we allow that the domain contains $\infty$. For $n \in \mathbb{Z}_{\ge 0}$, an $n$-connected domain is a domain $D$ such that $\widehat{\mathbb{C}} \setminus D$ is the disjoint union of $n$ connected compact sets, each of which contains more than one point. A finitely connected domain is an $n$-connected domain for some $n \in \mathbb{Z}_{\ge 0}$. A 0-connected domain is just $\widehat{\mathbb{C}}$. A 1-connected domain is conformally equivalent to the unit disc.

We will use dist (resp. $\text{dist}^\#$) to denote the Euclidean (resp. spherical) distance; use diam (resp. $\text{diam}^\#$) to denote the Euclidean (resp. spherical) diameter; and use $\mathbf{B}(z_0; r)$ [resp. $\mathbf{B}^\#(z_0; r)$] to denote the ball centered at $z_0$ with radius $r$, in the Euclidean (resp. spherical) metric. Let $\partial^\# D$ denote the boundary of $D$ in $\widehat{\mathbb{C}}$; and let $\partial D = \partial^\# D \cap \mathbb{C}$.

Suppose $D$ is an $n$-connected domain. Then $\partial^\# D$ has $n$ connected components, each of which is the boundary of a connected component of $\widehat{\mathbb{C}} \setminus D$. If $f$ maps $D$ conformally into $\widehat{\mathbb{C}}$, then $D' := f(D)$ is also an $n$-connected domain. And $f$ induces a one-to-one correspondence $\breve{f}$ from the set of components of $\partial^\# D$ to the set of components of $\partial^\# D'$ such that for any component $A$ of $\partial^\# D$ and $z \in D$, $z \to A$ iff $f(z) \to \breve{f}(A)$. There exists some $f$ that maps $D$ conformally onto a plane domain that is bounded by $n$ mutually disjoint analytic Jordan curves. We call such $f$ a boundary smoothing map of $D$. Suppose $f_1$ and $f_2$ are two boundary smoothing maps of $D$, and $E_j = f_j(D)$, $j = 1, 2$. Then $f_2 \circ f_1^{-1}$ maps $E_1$ conformally onto $E_2$, and $f_2 \circ f_1^{-1}$ induces a one-to-one correspondence $J$ from the set of Jordan curves that bound $E_1$ to the set of Jordan curves that bound $E_2$ such that for any Jordan curve $\sigma$ that bounds $E_1$ and $z \in E_1$, $z \to \sigma$ iff $f_2 \circ f_1^{-1}(z) \to J(\sigma)$. Since $\sigma$ and $J(\sigma)$ are both analytic, from the Schwarz reflection principle, $f_2 \circ f_1^{-1}$ can be extended conformally across $\sigma$, and maps $\sigma$ onto $J(\sigma)$.

Now consider the set of all pairs $(f, z)$ such that $f$ is a boundary smoothing map of $D$, and $z \in \overline{f(D)}$. Two pairs $(f_1, z_1)$ and $(f_2, z_2)$ are equivalent if the extension of $f_2 \circ f_1^{-1}$ maps $z_1$ to $z_2$. Let $\widehat{D}$ be the set of all equivalent classes. There is a unique conformal structure on $\widehat{D}$ such that $z \mapsto [(f, z)]$ maps $\overline{f(D)}$ conformally onto $\widehat{D}$ for any boundary smoothing map $f$. Then $z \mapsto [(f, f(z))]$ is a conformal map from $D$ into $\widehat{D}$ independent of the choice of $f$. So we may view $D$ as a subset of $\widehat{D}$, and call $\widehat{D}$ the conformal closure of $D$. It is clear that a conformal map between two finitely connected domains extends uniquely to a conformal map between their conformal closures.

We call $\widehat{\partial} D := \widehat{D} \setminus D$ the conformal boundary of $D$. Then $\widehat{\partial} D$ is a union of $n$ disjoint analytic Jordan curves, each of which is called a side of $D$. Each side $\sigma$ corresponds to a component $A$ of $\partial^\# D$ such that for $z \in D$, $z \to \sigma$ in $\widehat{D}$ iff $z \to A$. Each point on $\sigma$ is called a prime end of $D$ on $A$. This is equivalent to the prime ends defined in [1] and [13]. In fact, the definition in [1] describes the property of a sequence of points in $D$ that converges to a point on $\widehat{\partial} D$, and the definition in [12] describes a neighborhood basis bounded by crosscuts of a point on $\partial \widehat{D}$. A connected subset of a side that contains more than one point is called a side arc.

If $z_0 \in \widehat{\mathbb{C}}$ and a prime end $w_0$ of $D$ satisfies that for $z \in D$, $z \to z_0$ iff $z \to w_0$ in $\widehat{D}$, then we say the point $z_0$ and the prime end $w_0$ correspond to each other, and we do not distinguish the point $z_0$ from the prime end $w_0$. For example, if a boundary component of $D$ is a Jordan curve, then each point on this curve corresponds to a prime end. If $z_0 \in \partial D$ and for some $\varepsilon > 0$, $\mathbf{B}(z_0;\varepsilon) \setminus D$ is a simple curve $\gamma$ connecting $z_0$ with $\{|z - z_0| = \varepsilon\}$, then $z_0$ corresponds to a prime end of $D$. But every other point on $\gamma$ corresponds to two prime ends of $D$.

If $\alpha: (a,b) \to D$ is a curve in $D$, and for some $z_0 \in \partial^{\#} D$, $\alpha(t) \to z_0$ as $t \to a$, then there is some prime end $w_0$ of $D$ such that $\alpha(t) \to w_0$ in $\widehat{D}$ as $t \to a$. Such $w_0$ is called the prime end determined by $\alpha$ at one end. In general, not every prime end of $D$ can be determined by a curve in $D$ in this way.

2.3. *Positive harmonic functions.* Suppose $D$ is a finitely connected domain, and $z_0 \in D$. The Green function $G(D, z_0; \cdot)$ in $D$ with the pole at $z_0$ is the continuous function defined on $\widehat{D} \setminus \{z_0\}$ which vanishes on $\widehat{\partial} D$, is positive and harmonic in $D \setminus \{z_0\}$, and $G(D, z_0; z)$ behaves like $-\ln|z - z_0|/(2\pi)$ near $z_0$ if $z_0 \neq \infty$; behaves like $\ln|z|/(2\pi)$ near $\infty$ if $z_0 = \infty$.

Suppose $w_0$ is a prime end of $D$. There is a continuous function $P$ defined on $\widehat{D} \setminus \{w_0\}$ which vanishes on $\widehat{D} \setminus \{w_0\}$, and is harmonic and positive in $D$. It is called a generalized Poisson kernel in $D$ with the pole at $w_0$. Such $P$ is not unique. But any two generalized Poisson kernels in $D$ with the pole at $w_0$ differ by a positive multiple constant. Suppose $z_0 \in \partial D$, and $\partial D$ is analytic near $z_0$; then $z_0$ corresponds to a prime end of $D$, and the Poisson kernel in $D$ with the pole at $z_0$ in the usual sense is well defined, and is an example of a generalized Poisson kernel in $D$ with the pole at $z_0$.

Suppose $I$ is a side arc of $D$. The harmonic measure function $H(D, I; \cdot)$ is a bounded continuous function defined on $\widehat{D}$ taking away the end points of $I$, which is harmonic in $D$, vanishes on $\widehat{\partial} D \setminus \overline{I}$, and takes constant value 1 on $I$ except the end points. For any $z \in D$, $H(D, I; z)$ is equal to the probability that the plane Brownian motion started from $z$ first hits $\partial D$ at $I$.

2.4. *Hulls and Loewner chains.* Suppose $D$ is an $n$-connected domain, and $\sigma$ is a side of $D$. Let $A(\sigma)$ be the connected component of $\widehat{\mathbb{C}} \setminus D$ that corresponds to $\sigma$. A closed subset $H$ is called a hull of $D$ on $\sigma$ if $D \setminus H$ is also an $n$-connected domain, and $A(\sigma) \cup H$ is a component of $\widehat{\mathbb{C}} \setminus (D \setminus H)$. Then other components of $\widehat{\mathbb{C}} \setminus (D \setminus H)$ are the components of $\widehat{\mathbb{C}} \setminus D$ other than $A(\sigma)$.

In this paper, we define a crosscut to be an open simple curve $\alpha$ in $D$, whose two ends approach to two points on $\partial D$, in the Lebesgue metric, such that $D \setminus \alpha$ has two components, one of which is simply connected. If $U$ is a simply connected component of $D \setminus \alpha$, then $U \cup \alpha$ is a hull in $D$. If $n > 1$, that is, $D$ is not simply connected, then $U$ is determined by $\alpha$, and let $H(\alpha) := U \cup \alpha$ be the hull bounded by $\alpha$. If $n = 1$, then the two components of $D \setminus \alpha$ are both simply connected, so we need some other restrictions to determine $H(\alpha)$. For example, if we say that $H(\alpha)$ is a neighborhood of some prime end $w_0$ in $D$, then there is no ambiguity.

Suppose $\sigma$ is a side of $D$. A Loewner chain in $D$ on $\sigma$ is a function $L$ from $[0, T)$ for some $T \in (0, +\infty]$ into the set of hulls in $D$ on $\sigma$ such that $L(0) = \varnothing$, $L(t_1) \subsetneqq L(t_2)$ if $0 \le t_1 < t_2 < T$, and for any fixed $b \in [0, T)$ and any compact subset $F$ of $D \setminus L(b)$, the extremal length (see [1]) of the family of curves in $D \setminus L(t+\varepsilon)$ that separates $F$ from $L(t+\varepsilon) \setminus L(t)$ tends to 0 as $\varepsilon \to 0^+$, uniformly w.r.t. $t \in [0, b]$. Suppose $L(t)$, $0 \le t < T$, is a Loewner chain in $D$ on $\sigma$. For each $t \in [0, T)$, let $d_t$ be any metric on $\widehat{D \setminus L(t)}$. From the definition, the $d_t$-diameter of $L(t+\varepsilon) \setminus L(t)$ tends to 0 as $\varepsilon \to 0^+$. Thus there is a unique prime end $w(t)$ of $D \setminus L(t)$ that lies on the closure of $L(t+\varepsilon) \setminus L(t)$ in $\widehat{D \setminus L(t)}$ for all $\varepsilon > 0$. We call $w(t)$ the prime end determined by $L$ at time $t$. Especially, $w(0)$ is a prime end on $\sigma$. We say $L$ is a Loewner chain started from $w(0)$. It is clear that for any $b \in [0, T)$, $t \mapsto L(b+t)$, $0 \le t < T - b$, is a Loewner chain in $D \setminus L(b)$ started from $w(b)$.

Suppose $L(t)$, $0 \le t < T$, is a Loewner chain in $D$. Suppose $u$ is a continuous (strictly) increasing function defined on $[0, T)$ with $u(0) = 0$. Let $u(T) := \sup u([0, T))$. Then $L'(t) := L(u^{-1}(t))$, $0 \le t < u(T)$, is also a Loewner chain in $D$. Such $L'$ is called a time-change of $L$ through $u$. Moreover, the prime end determined by $L'$ at time $u(t)$ is the same as the prime end determined by $L$ at time $t$.

One example of a Loewner chain is constructed by a simple curve. Suppose $\gamma : [0, T) \to \widehat{D}$ is a simple curve that satisfies $\gamma(0) \in \widehat{\partial} D$ and $\gamma(t) \in D$ for $0 \le t < T$. Let $L(t) = \gamma((0, t])$, $0 \le t < T$. Then $L$ is a Loewner chain in $D$ started from $\gamma(0)$, and $\gamma(t)$ corresponds to the prime end determined by $L$ at time $t$. We say that $L$ is the Loewner chain generated by $\gamma$.

## 3. Continuous LERW.

3.1. *Chordal Loewner equation.* Let $\mathbb{H} = \{z \in \mathbb{C} : \operatorname{Im} z > 0\}$. Then $\mathbb{H}$ is a 1-connected domain whose side is $\widehat{\mathbb{R}} := \mathbb{R} \cup \{\infty\}$. We say $H$ is a hull in $\mathbb{H}$ w.r.t. $\infty$ if $H$ is a hull in $\mathbb{H}$ and $H$ is bounded (i.e., bounded away from $\infty$). A Loewner chain $L$ in $\mathbb{H}$ w.r.t. $\infty$ is a Loewner chain in $\mathbb{H}$ such that each $L(t)$ is a hull in $\mathbb{H}$ w.r.t. $\infty$. For each hull $H$ in $\mathbb{H}$ w.r.t. $\infty$, there is a unique function $\varphi_H$ that maps $\mathbb{H} \setminus H$ conformally onto $\mathbb{H}$ such that for some $c \ge 0$,

$$\varphi_H(z) = z + \frac{c}{z} + O\left(\frac{1}{z^2}\right),$$

as $z \to \infty$. Such $c$ is called the capacity of $H$ in $\mathbb{H}$ w.r.t. $\infty$, denoted by $\operatorname{hcap}(H)$. The empty set is a hull in $\mathbb{H}$ w.r.t. $\infty$, and $\varphi_\varnothing = \mathrm{id}$, so $\operatorname{hcap}(\varnothing) = 0$.

PROPOSITION 3.1. *Suppose $\Omega$ is an open neighborhood of $x_0 \in \mathbb{R}$ in $\mathbb{H}$. Suppose $W$ maps $\Omega$ conformally into $\mathbb{H}$ such that for some $r > 0$, if $z \to (x_0 - r, x_0 + r)$ in $\Omega$, then $W(z) \to \mathbb{R}$. So $W$ extends conformally across $(x_0 - r, x_0 + r)$ by the Schwarz reflection principle. Then for any $\varepsilon > 0$, there is some $\delta > 0$ such that if a hull $H$ in $\mathbb{H}$ w.r.t. $\infty$ is contained in $\{z \in \mathbb{H} : |z - x_0| < \delta\}$, then $W(H)$ is also a hull in $\mathbb{H}$ w.r.t. $\infty$, and*

$$|\operatorname{hcap}(W(H)) - W'(x_0)^2 \operatorname{hcap}(H)| \le \varepsilon |\operatorname{hcap}(H)|.$$

PROOF. This is Lemma 2.8 in [6]. □

For $T \in (0, +\infty]$, let $C([0,T))$ denote the space of real-valued continuous functions on $[0,T)$. Suppose $\xi \in C([0,T))$. We solve the chordal Loewner equation:

$$\partial_t \varphi_t(z) = \frac{2}{\varphi_t(z) - \xi(t)}, \qquad \varphi_0(z) = z,$$

for $0 \le t < T$. For each $t \in [0,T)$, let $K_t$ be the set of $z \in \mathbb{H}$ such that the solution $\varphi_s(z)$ blows up before or at time $t$. We say that $\varphi_t$ and $K_t$, $0 \le t < T$, are chordal Loewner maps and hulls, respectively, driven by $\xi$.

For $0 \le t < T$, $K_t$ is a bounded closed subset of $\mathbb{H}$, $\varphi_t$ maps $\mathbb{H} \setminus K_t$ conformally onto $\mathbb{H}$, and satisfies

$$\varphi_t(z) = z + \frac{2t}{z} + O\left(\frac{1}{z^2}\right)$$

as $z \to \infty$. So $K_t$ is a hull in $\mathbb{H}$ w.r.t. $\infty$, $\operatorname{hcap}(K_t) = 2t$ and $\varphi_{K_t} = \varphi_t$.

PROPOSITION 3.2. (i) *Suppose $\varphi_t$ and $K_t$, $0 \le t < T$, are chordal Loewner maps and hulls, respectively, driven by $\xi$. Then $t \mapsto K_t$, $0 \le t < T$, is a Loewner chain in $\mathbb{H}$ w.r.t. $\infty$ started from $\xi(0)$. And for each $t \in [0,T)$, $\operatorname{hcap}(K_t) = 2t$, $\varphi_t = \varphi_{K_t}$, and*

$$\{\xi(t)\} = \bigcap_{\varepsilon > 0} \overline{\varphi_t(K_{t+\varepsilon} \setminus K_t)}.$$

(ii) *Suppose $L(t)$, $0 \le t < T$, is a Loewner chain in $\mathbb{H}$ w.r.t. $\infty$. Let $v(t) = \operatorname{hcap}(L(t))/2$, $0 \le t < T$. Then $v$ is a continuous increasing function with $v(0) = 0$. And $K_t := L(v^{-1}(t))$, $0 \le t < v(T)$, are chordal Loewner hulls driven by some $\xi \in C([0, v(T)))$.*

PROOF. This is almost the same as Theorem 2.6 in [6]. $\square$

Fix $b \in [0,T)$. Let $\varphi_{b,t} = \varphi_t \circ \varphi_b^{-1}$ and $K_{b,t} = \varphi_b(K_t \setminus K_b)$ for $b \le t < T$. Then it is easy to check that $K_{b,b+t}$ and $\varphi_{b,b+t}$, $0 \le t < T-b$, are chordal Loewner hulls and maps driven by $t \mapsto \xi(b+t)$, $0 \le t < T-b$. Thus for any $s < t \in [0,T)$, $\varphi_s(K_t \setminus K_s)$ is a hull in $\mathbb{H}$ w.r.t. $\infty$, and its capacity in $\mathbb{H}$ w.r.t. $\infty$ is $2(t-s)$.

3.2. *Continuous LERW aiming at an interior point.* We define an almost $\mathbb{H}$ domain to be a finitely connected domain in $\mathbb{H}$ that is bounded by $\widehat{\mathbb{R}}$ and mutually disjoint analytic Jordan curves in $\mathbb{H}$. Let $\Omega$ be an almost $\mathbb{H}$ domain, and $p \in \Omega$. If $K$ is a hull in $\mathbb{H}$ w.r.t. $\infty$ such that $K \subset \Omega \setminus \{p\}$, let $\Omega_K = \varphi_K(\Omega \setminus K)$. Then $\Omega_K$ is also an almost $\mathbb{H}$ domain, and $\varphi_K(p) \in \Omega_K$.

For $a \ge 0$, let $C([0,a])$ be the space of all real-valued continuous functions defined on $[0,a]$ with norm $\|\xi\|_a := \sup\{|\xi(t)| : 0 \le t \le a\}$. For $\xi \in C([0,a])$, let $K_t^\xi$ and $\varphi_t^\xi$, $0 \le t \le a$, be chordal Loewner hulls and maps, respectively, driven by $\xi$. If $K_t^\xi \subset \Omega \setminus \{p\}$, we write $\Omega_t^\xi$ for $\Omega_{K_t^\xi}$. Define

$$J_t^\xi(z) = G(\Omega \setminus K_t^\xi, p; \cdot) \circ (\varphi_t^\xi)^{-1}. \tag{3.1}$$

Since $J_t^\xi = G(\Omega_t^\xi, \varphi_t^\xi(p); \cdot)$ is positive and harmonic in $\Omega_t^\xi \setminus \{\varphi_t^\xi(p)\}$, and vanishes on $\mathbb{R}$, so it extends harmonically across $\mathbb{R}$. Let

$$X_t^\xi = (\partial_x \partial_y / \partial_y) J_t^\xi(\xi(t)) = \partial_x \partial_y J_t^\xi(\xi(t)) / \partial_y J_t^\xi(\xi(t)).$$

We begin with a theorem. The proof is postponed to Section 5 in this paper.

THEOREM 3.1. *For any $A \in C([0,\infty))$ and $\lambda \in \mathbb{R}$, the equation*

$$\xi(t) = A(t) + \lambda \int_0^t X_s^\xi \, ds \tag{3.2}$$

*has a unique maximal solution $\xi(t) = \xi_A(t)$, $0 \le t < T_A$, where $T_A \in (0,\infty]$. Here "maximal" means that the solution cannot be extended. Moreover, we have:*

(i) *For any $a \in (0,\infty)$, the set $\{A \in C([0,\infty)) : T_A > a\}$ is open w.r.t. the metric $\|\cdot\|_a$, and $A \mapsto \xi_A$ is $(\|\cdot\|_a, \|\cdot\|_a)$ continuous on $\{A \in C([0,\infty)) : T_A > a\}$.*

(ii) *There is no crosscut $\alpha$ in $\mathbb{H}$ such that $\bigcup_{0 \le t < T_A} K_t^\xi \subset H(\alpha) \subset \Omega \setminus \{p\}$.*

Suppose $D$ is a finitely connected domain, $w_0$ is a prime end of $D$, and $z_e \in D$. There is $f$ that maps $D$ conformally onto an almost $\mathbb{H}$ domain $\Omega$, such that $f(w_0) = 0$. Let $p = f(z_e)$, $B(t)$ be a Brownian motion, and $\xi(t)$, $0 \le t < T$, be the maximal solution to (3.2) with $A(t) = \sqrt{2}B(t)$ and $\lambda = 2$.

Let $\{\mathcal{F}_t\}$ be the filtration generated by $B(t)$. From Theorem 3.1(i), $T$ is an $\{\mathcal{F}_t\}$-stopping time, and $(\xi(t))$ is $\{\mathcal{F}_t\}$-adapted. For $0 \le t < T$, let

$$u(t) = \int_0^t (\partial_y J_s^\xi(\xi(s)))^2 \, ds. \tag{3.3}$$

Let $S = u(T)$, and $F(t) = f^{-1}(K^\xi_{u^{-1}(t)})$, $0 \le t < S$. In the next subsection, we will prove the following theorem.

Theorem 3.2. *For $j = 1, 2$, suppose $f_j$ maps $D$ conformally onto some almost $\mathbb{H}$ domain $\Omega_j$ such that $f_j(w_0) = 0$. For $j = 1, 2$, let $p_j = f_j(z_e)$, $B_j(t)$ be a Brownian motion, and $\xi_j(t)$, $0 \le t < T_j$, be the maximal solution to*

$$\begin{aligned} \xi_j(t) &= \sqrt{2} B_j(t) \\ &\quad + 2\int_0^t (\partial_x \partial_y / \partial_y)(G(\Omega_j \setminus K_s^{\xi_j}, p_j; \cdot) \circ (\varphi_s^{\xi_j})^{-1})(\xi_j(s)) \, ds; \end{aligned} \tag{3.4}$$

*and let $u_j(t)$, $0 \le t < T_j$, be defined by*

$$u_j(t) = \int_0^t \partial_y (G(\Omega_j \setminus K_s^{\xi_j}, p_j; \cdot) \circ (\varphi_s^{\xi_j})^{-1})(\xi_j(s))^2 \, ds.$$

*Let $S_j = u_j(T)$ and $F_j(t) = f_j^{-1}(K^{\xi_j}_{u_j^{-1}(t)})$, $0 \le t < S_j$, $j = 1, 2$. Then $(F_1(t), 0 \le t < S_1)$ and $(F_2(t), 0 \le t < S_2)$ have the same distribution.*

Thus the distribution of $(F(t), 0 \le t < S)$ does not depend on the choice of $f$, and is conformally invariant. We call $(F(t), 0 \le t < S)$ a continuous LERW in $D$ from $w_0$ to $z_e$, and let it be denoted by $\mathrm{LERW}(D; w_0 \to z_e)$. From the property of chordal $\mathrm{SLE}_2$ (cf. [15]) and Girsanov's theorem [11, 14], almost surely there is a simple curve $\gamma(t) : [0, S) \to \widehat{D}$ such that $\gamma(0) = w_0$, $\gamma(t) \in D$ for $0 < t < S$, and $F(t) = \gamma((0, t])$ for $0 \le t < S$, that is, $F$ is the Loewner chain generated by $\gamma$. We call such $\gamma$ an $\mathrm{LERW}(D; w_0 \to z_e)$ trace.

Remark. If $D$ is a 1-connected domain, $w_0$ is a prime end of $D$ and $z_e \in D$, then an $\mathrm{LERW}(D; w_0 \to z_e)$ has the same distribution as a radial $\mathrm{SLE}_2(D; w_0 \to z_e)$ up to a linear time-change.

### 3.3. *Conformal invariance.*

Proof of Theorem 3.2. For $j = 1, 2$, let $v_j = u_j^{-1}$ and $L_j(t) = K^{\xi_j}_{u_j^{-1}(t)}$. Then $F_j(t) = f_j^{-1}(L_j(t))$, $0 \le t < S_j$. Let $W = f_2 \circ f_1^{-1}$. Then $W$ maps $\Omega_1$ conformally onto $\Omega_2$, $W(0) = 0$ and $W(p_1) = p_2$. Let $L_{2'}(t) = W(L_1(t))$, $0 \le t < S_1$. It suffices to show that $(L_{2'}(t), 0 \le t < S_1)$ has the same distribution as $(L_2(t), 0 \le t < S_2)$. Let $\beta_1(t)$ be the random simple curve that generates $L_1(t)$, that is, $\beta_1(0)=0$, $\beta_1(t)\in\Omega_1$, $0<t<S_1$, and $L_1(t)=\beta_1((0,t])$, $0\le t<S_1$. Let $\beta_{2'}(t)=W(\beta_1(t))$, $0\le t<S_1$. Then $\beta_{2'}$ is a simple curve, $\beta_{2'}(0)=0$, $\beta_{2'}(t)\in\Omega_2\subset\mathbb{H}$, $0<t<S_1$, and $L_{2'}(t)=\beta_{2'}((0,t])$, $0\le t<S_1$. Thus $L_{2'}$ is a Loewner chain in $\mathbb{H}$ w.r.t. $\infty$. Let $v_{2'}(t)=\mathrm{hcap}(L_{2'}(t))/2$, $0\le t<S_1$. Let $T_{2'}=v_{2'}(S_1)$ and $u_{2'}=v_{2'}^{-1}$. Then from Proposition 3.2, $L_{2'}(u_{2'}(t))=K_t^{\xi_{2'}}$, $0\le t<T_{2'}$, for some $\xi_{2'}\in C([0,T_{2'}))$.

Let $\{\mathcal{F}_t^1\}$ be the filtration generated by $B_1(t)$. Let

$$R_1(t,x)=\partial_y(G(\Omega_1\setminus K_t^{\xi_1},p_1;\cdot)\circ(\varphi_t^{\xi_1})^{-1})(x).$$

From Theorem 3.1(i), $(\xi_1(t))$ and $R_1(t,x)$ are $\mathcal{F}_t^1$-adapted, and $T_1$ is an $\mathcal{F}_t^1$-stopping time. Thus for $0\le t<T_1$, we have

$$d\xi_1(t)=\sqrt{2}\,dB_1(t)+2\,\frac{\partial_x R_1(t,\xi_1(t))}{R_1(t,\xi_1(t))}\,dt$$

and

$$u_1'(t)=R_1(t,\xi_1(t))^2.$$

So there is another Brownian motion $\breve{B}_1(t)$ such that for $0\le t<S_1$,

$$d\xi_1(v_1(t))=\frac{\sqrt{2}}{R_1(v_1(t),\xi_1(v_1(t)))}\,d\breve{B}_1(t)+2\,\frac{\partial_x R_1(v_1(t),\xi_1(v_1(t)))}{R_1(v_1(t),\xi_1(v_1(t)))^3}\,dt. \tag{3.5}$$

Note that $W$ maps $\Omega_1\setminus L_1(t)$ conformally onto $\Omega_2\setminus L_{2'}(t)$. Let $\Omega_1(t)=\varphi^{\xi_1}_{v_1(t)}(\Omega_1\setminus L_1(t))$, $\Omega_{2'}(t)=\varphi^{\xi_{2'}}_{v_{2'}(t)}(\Omega_2\setminus L_{2'}(t))$ and $W_t=\varphi^{\xi_{2'}}_{v_{2'}(t)}\circ W\circ(\varphi^{\xi_1}_{v_1(t)})^{-1}$. Then both $\Omega_1(t)$ and $\Omega_{2'}(t)$ are almost $\mathbb{H}$ domains, and $W_t$ maps $\Omega_1(t)$ conformally onto $\Omega_{2'}(t)$, and maps $\widehat{\mathbb{R}}$ onto itself. For $t\in[0,S_1)$ and $\varepsilon\in[0,S_1-t)$, define $L_1(t,\varepsilon)=\varphi^{\xi_1}_{v_1(t)}(K^{\xi_1}_{v_1(t+\varepsilon)}\setminus K^{\xi_1}_{v_1(t)})$ and $L_{2'}(t,\varepsilon)=\varphi^{\xi_{2'}}_{v_{2'}(t)}(K^{\xi_{2'}}_{v_{2'}(t+\varepsilon)}\setminus K^{\xi_{2'}}_{v_{2'}(t)})$. Then $\mathrm{hcap}(L_1(t,\varepsilon))=2(v_1(t+\varepsilon)-v_1(t))$, $\mathrm{hcap}(L_{2'}(t,\varepsilon)=2(v_{2'}(t+\varepsilon)-v_{2'}(t)))$, and $\overline{W_t(L_1}(t,\varepsilon))=L_{2'}(t,\varepsilon)$. From Proposition 3.2, we have $\{\xi_1(v_1(t))\}=\bigcap_{\varepsilon>0}\overline{L_1(t,\varepsilon)}$ and $\{\xi_{2'}(v_{2'}(t))\}=\bigcap_{\varepsilon>0}\overline{L_{2'}(t,\varepsilon)}$. Thus $\xi_{2'}(v_{2'}(t))=W_t(\xi_1(v_1(t)))$. From Proposition 3.1, we have $v_{2'}'(t)=W_t'(\xi_1(v_1(t)))^2v_1'(t)$.

Differentiate the equality $W_t\circ\varphi^{\xi_1}_{v_1(t)}=\varphi^{\xi_{2'}}_{v_{2'}(t)}\circ W$ w.r.t. $t$. We get

$$\partial_t W_t(\varphi^{\xi_1}_{v_1(t)}(z))+\frac{2W_t'(\varphi^{\xi_1}_{v_1(t)}(z))v_1'(t)}{\varphi^{\xi_1}_{v_1(t)}(z)-\xi_1(v_1(t))}=\frac{2v_{2'}'(t)}{\varphi^{\xi_{2'}}_{v_{2'}(t)}\circ W(z)-\xi_{2'}(v_{2'}(t))}$$

for any $z\in\Omega_1\setminus L_1(t)$. Since $\varphi^{\xi_1}_{v_1(t)}$ maps $\Omega_1\setminus L_1(t)$ conformally onto $\Omega_1(t)$, so for any $w\in\Omega_1(t)$, we have

$$\partial_t W_t(w)=\frac{2W_t'(\xi_1(v_1(t)))^2v_1'(t)}{W_t(w)-W_t(\xi_1(v_1(t)))}-\frac{2W_t'(w)v_1'(t)}{w-\xi_1(v_1(t))}.$$

Let $w \to \xi_1(v_1(t))$ in $\Omega_1(t)$; from Taylor expansion of $W_t$ at $\xi_1(v_1(t))$, we get

$$\begin{aligned}\partial_t W_t(\xi_1(v_1(t))) &= -3W_t''(\xi_1(v_1(t)))v_1'(t)\\ &= -3W_t''(\xi_1(v_1(t)))/R_1(v_1(t),\xi_1(v_1(t)))^2.\end{aligned}$$

Since $\xi_{2'}(v_{2'}(t)) = W_t(\xi_1(v_1(t)))$, so from (3.5) and Itô's formula [11, 14], we have

$$\begin{aligned}d\xi_{2'}(v_{2'}(t)) &= \partial_t W_t(\xi_1(v_1(t)))\,dt + W_t'(\xi_1(v_1(t)))\,d\xi_1(v_1(t))\\
&\quad + W_t''(\xi_1(v_1(t)))\,d\langle \xi_1(v_1(t))\rangle/2\\
&= \frac{\sqrt{2}W_t'(\xi_1(v_1(t)))}{R_1(v_1(t),\xi_1(v_1(t)))}\,d\breve{B}_1(t)\\
&\quad + 2\,\frac{W_t'(\xi_1(v_1(t)))\,\partial_x R_1(v_1(t),\xi_1(v_1(t)))}{R_1(v_1(t),\xi_1(v_1(t)))^3}\,dt\\
&\quad + \partial_t W_t(\xi_1(v_1(t)))\,dt + \frac{W_t''(\xi_1(v_1(t)))}{R_1(v_1(t),\xi_1(v_1(t)))^2}\,dt\\
&= \frac{\sqrt{2}W_t'(\xi_1(v_1(t)))}{R_1(v_1(t),\xi_1(v_1(t)))}\,d\breve{B}_1(t)\\
&\quad + 2\Big(\frac{W_t'(\xi_1(v_1(t)))\,\partial_x R_1(v_1(t),\xi_1(v_1(t)))}{R_1(v_1(t),\xi_1(v_1(t)))^3}\\
&\qquad\qquad - \frac{W_t''(\xi_1(v_1(t)))}{R_1(v_1(t),\xi_1(v_1(t)))^2}\Big)\,dt.\end{aligned} \tag{3.6}$$

Since $\varphi^{\xi_1}_{v_1(t)}$ maps $\Omega_1 \setminus L_1(t)$ conformally onto $\Omega_1(t)$, so

$$\begin{aligned}R_1(v_1(t),x) &= \partial_y(G(\Omega_1\setminus K^{\xi_1}_{v_1(t)},p_1;\cdot)\circ(\varphi^{\xi_1}_{v_1(t)})^{-1})(x)\\ &= \partial_y G(\Omega_1(t),\varphi^{\xi_1}_{v_1(t)}(p_1);\cdot)(x).\end{aligned}$$

Since $W_t$ maps $\Omega_1(t)$ conformally onto $\Omega_{2'}(t)$, and $W_t(\varphi^{\xi_1}_{v_1(t)}(p_1)) = \varphi^{\xi_{2'}}_{v_{2'}(t)}(p_2)$, so

$$G(\Omega_1(t),\varphi^{\xi_1}_{v_1(t)}(p_1);\cdot) = G(\Omega_{2'}(t),\varphi^{\xi_{2'}}_{v_{2'}(t)}(p_2);\cdot)\circ W_t.$$

Thus

$$\begin{aligned}R_1(v_1(t),x) &= \partial_y G(\Omega_{2'}(t),\varphi^{\xi_{2'}}_{v_{2'}(t)}(p_2);W_t(x))W_t'(x);\\
\partial_x R_1(v_1(t),x) &= \partial_x\partial_y G(\Omega_{2'}(t),\varphi^{\xi_{2'}}_{v_{2'}(t)}(p_2);W_t(x))(W_t'(x))^2\\
&\quad + \partial_y G(\Omega_{2'}(t),\varphi^{\xi_{2'}}_{v_{2'}(t)}(p_2);W_t(x))W_t''(x).\end{aligned}$$

Plugging these equalities into (3.6) and letting $x = \xi_1(v_1(t))$, we get

$$d\xi_{2'}(v_{2'}(t)) = \frac{\sqrt{2}W_t'(\xi_1(v_1(t)))}{R_1(v_1(t),\xi_1(v_1(t)))}\, d\breve{B}_1(t) + 2\frac{\partial_x\partial_y G(\Omega_{2'}(t),\varphi^{\xi_{2'}}_{v_{2'}(t)}(p_2);\xi_{2'}(v_{2'}(t)))}{\partial_y G(\Omega_{2'}(t),\varphi^{\xi_{2'}}_{v_{2'}(t)}(p_2);\xi_{2'}(v_{2'}(t)))^3}\, dt.$$

Since

$$v_{2'}'(t) = W_t'(\xi_1(v_1(t)))^2 v_1'(t) = \frac{W_t'(\xi_1(v_1(t)))^2}{R_1(v_1(t),\xi_1(v_1(t)))^2} = \partial_y G(\Omega_{2'}(t),\varphi^{\xi_{2'}}_{v_{2'}(t)}(p_2);\xi_{2'}(v_{2'}(t)))^{-2}, \tag{3.7}$$

and $G(\Omega_{2'}(t),\varphi^{\xi_{2'}}_{v_{2'}(t)}(p_2);\cdot) = G(\Omega_2\setminus K^{\xi_{2'}}_{v_{2'}(t)},p_2;\cdot)\circ\varphi^{\xi_{2'}}_{v_{2'}(t)}$, so for $0\le t< T_{2'}$,

$$d\xi_{2'}(t) = \sqrt{2}\, dB_{2'}(t) + 2(\partial_x\partial_y/\partial_y)(G(\Omega_2\setminus K^{\xi_{2'}}_t,p_2;\cdot)\circ\varphi^{\xi_{2'}}_t)(\xi_{2'}(t))\, dt$$

for another Brownian motion $B_{2'}(t)$. Since $\xi_{2'}(0) = W_0(\xi_1(0)) = W(0) = 0$, so for $0\le t<T_{2'}$,

$$\xi_{2'}(t) = \sqrt{2}B_{2'}(t) + 2\int_0^t (\partial_x\partial_y/\partial_y)(G(\Omega_2\setminus K^{\xi_{2'}}_s,p_2;\cdot)\circ\varphi^{\xi_{2'}}_s)(\xi_{2'}(s))\, ds. \tag{3.8}$$

We claim that $\xi_{2'}(t)$, $0\le t<T_{2'}$, is the maximal solution to (3.8). Suppose the claim is not true. Then it may happen that the solution $\xi_{2'}$ extends to $[0,T_{2'}]$. Note that $W(\infty)$ is a prime end on $\widehat{\mathbb{R}}$ other than $W(0)=0$. We may find a crosscut $\alpha$ in $\mathbb{H}$ such that $K^{\xi_{2'}}_{T_{2'}}\subset H(\alpha)\subset\Omega_2\setminus\{p_2\}$, and $W(\infty)\notin\overline{H(\alpha)}$. Then $W^{-1}(\alpha)$ is also a crosscut in $\mathbb{H}$, and $H(W^{-1}(\alpha)) = W^{-1}(H(\alpha))\subset\Omega_1\setminus\{p_1\}$. So $W^{-1}(K^{\xi_{2'}}_t)\subset H(W^{-1}(\alpha))$ for $0\le t<T_{2'}$, which implies that $K^{\xi_1}_t\subset H(W^{-1}(\alpha))$ for $0\le t<T_1$. This contradicts Theorem 3.1(ii). So the claim is justified.

Since $\xi_{2'}(t)$, $0\le t<T_{2'}$ [resp. $\xi_2(t)$, $0\le t<T_2$], is the maximal solution to (3.8) [resp. equation (3.2) when $j=2$], and $(B_{2'}(t))$ has the same distribution as $(B_2(t))$, so $(\xi_{2'}(t), 0\le t<T_{2'})$ has the same distribution as $(\xi_2(t), 0\le t< T_2)$. From (3.7), $u_{2'} = v_{2'}^{-1}$, and that $u_{2'}(0)=0$, we see that for $0\le t<T_{2'}$,

$$u_{2'}(t) = \int_0^t \partial_y(G(\Omega_2\setminus K^{\xi_{2'}}_s,p_2;\cdot)\circ\varphi^{\xi_{2'}}_s)(\xi_{2'}(s))^2\, ds.$$

Thus $((\xi_{2'}(t),u_{2'}(t)), 0\le t<T_{2'})$ has the same distribution as $((\xi_2(t),u_2(t)),0\le t<T_2)$. Since $L_{2'}(t) = K^{\xi_{2'}}_{u_{2'}^{-1}(t)}$ for $0\le t<S_1 = u_{2'}(T_{2'})$, and $L_2(t) = K^{\xi_2}_{u_2^{-1}(t)}$ for $0 \le t < S_2 = u_2(T_2)$, so $(L_{2'}(t), 0 \le t < S_1)$ has the same distribution as $(L_2(t), 0 \le t < S_2)$. $\square$

3.4. *Continuous LERW with other kinds of targets.* Suppose $D$ is a finitely connected domain, $w_0$ is a prime end of $D$, and $I_e$ is a side arc of $D$ that is bounded away from $w_0$. Then there is $f$ that maps $D$ conformally onto an almost $\mathbb{H}$ domain $\Omega$ such that $f(w_0) = 0$. If a hull $K$ in $\mathbb{H}$ w.r.t. $\infty$ is bounded away from $f(I_e)$, and $K \subset \Omega$, then $f(I_e)$ is a side arc of $\Omega \setminus K$. We have the harmonic measure function $H(\Omega \setminus K, f(I_e); \cdot)$.

Now we change the definition of $J_t^\xi$ by replacing $G(\Omega \setminus K_t^\xi, p; \cdot)$ by $H(\Omega \setminus K_t^\xi, f(I_e); \cdot)$ in (3.1), and still let $X_t^\xi = (\partial_x \partial_y / \partial_y) J_t^\xi(\xi(t))$. Let everything else in Section 3.2 be unchanged. Then Theorem 3.1 still holds if the condition on $\alpha$ is replaced by that $\alpha$ is a crosscut in $\mathbb{H}$ such that $H(\alpha) \subset \Omega$ and $H(\alpha)$ is bounded away from $f(I_e)$. Let $u(t)$ be defined by (3.3). Then $(F(t) = f^{-1}(K^\xi_{u^{-1}(t)}), 0 \le t < S = u(T))$ is called a continuous LERW in $D$ from $w_0$ to $I_e$, and is denoted by LERW$(D; w_0 \to I_e)$. It is almost surely generated by a random simple curve, which is called an LERW$(D; w_0 \to I_e)$ trace. The variation of Theorem 3.2 for LERW$(D; w_0 \to I_e)$ still holds. Thus the distribution of LERW$(D; w_0 \to I_e)$ does not depend on the choice of $f$, and is conformally invariant.

Suppose $D$ is a finitely connected domain, $w_0$ and $w_e$ are two different prime ends of $D$. There is $f$ that maps $D$ conformally onto an almost $\mathbb{H}$ domain $\Omega$ such that $f(w_0) = 0$. Then $p := f(w_e)$ is a prime end of $\Omega$ other than 0. If a hull $K$ in $\mathbb{H}$ w.r.t. $\infty$ is bounded away from $p$, and $K \subset \Omega$, then $p$ is a prime end of $\Omega \setminus K$.

A normalization function is a function $h$ that maps a neighborhood $U$ of $p$ in $\widehat{\Omega}$ conformally onto a neighborhood $V$ of 0 in $\overline{\mathbb{H}}$ such that $h(p) = 0$ and $h(U \cap \widehat{\partial} D) \subset \mathbb{R}$. There is a unique generalized Poisson kernel $P(z)$ in $\Omega \setminus K$ with the pole at $p$ such that the principal part of $P \circ h^{-1}(z)$ at 0 is $\operatorname{Im} \frac{-1}{z}$. Let $P(\Omega \setminus K, p, h; z)$ denote this function.

Now fix a normalizing function $h$. Change the definition of $J_t^\xi$ by replacing $G(\Omega \setminus K_t^\xi, p; \cdot)$ by $P(\Omega \setminus K_t^\xi, p, h; \cdot)$ in (3.1), and still let $X_t^\xi = (\partial_x \partial_y / \partial_y) J_t^\xi(\xi(t))$. Let everything else in Section 3.2 be unchanged. Then Theorem 3.1 still holds if the condition on $\alpha$ is replaced by that $\alpha$ is a crosscut in $\mathbb{H}$ such that $H(\alpha) \subset \Omega$, and $H(\alpha)$ is bounded away from $p = f(w_e)$. Let $u(t)$ be defined by (3.3). Then $(F(t) = f^{-1}(K^\xi_{u^{-1}(t)}), 0 \le t < S = u(T))$ is called a continuous LERW in $D$ from $w_0$ to $w_e$, normalized by $h$, and is denoted by LERW$(D; w_0 \to w_e)$. It is almost surely generated by a random simple curve, which is called an LERW$(D; w_0 \to w_e)$ trace normalized by $h$. The variation of Theorem 3.2 for LERW$(D; w_0 \to w_e)$ holds with simple modification: $(F_1(t), 0 \le t < S_1)$ and $(F_2(t/a^2), 0 \le t < a^2 S_2)$ have the same distribution, where $a=(h_2\circ h_1^{-1})'(0)$ and $h_j$, $j=1,2$, are normalization functions. Thus the distribution of $\mathrm{LERW}(D;w_0\to w_e)$ up to a linear time-change does not depend on the choices of $f$ and $h$, and is conformally invariant.

REMARK. (i) If $D$ is a 1-connected domain, and $w_0\neq w_e$ are two prime ends of $D$, then an $\mathrm{LERW}(D;w_0\to w_e)$ has the same distribution as a chordal $\mathrm{SLE}_2(D;w_0\to w_e)$ up to a linear time-change.

(ii) If $D$ is a 1-connected domain, $w_0$ is a prime end of $D$, and $I_e$ is a side arc of $D$ that is bounded away from $w_0$, then an $\mathrm{LERW}(D;w_0\to I_e)$ has the same distribution as a strip or dipolar $\mathrm{SLE}_2(D;w_0\to I_e)$ (cf. [2, 21]) up to a linear time-change.

(iii) If $D$ is a 2-connected domain, $w_0$ is a prime end of $D$, and $I_e$ is a side of $D$ that does not contain $w_0$, then an $\mathrm{LERW}(D;w_0\to I_e)$ has the same distribution as an annulus $\mathrm{SLE}_2(D;w_0\to I_e)$ (cf. [20]) up to a deterministic time-change.

## 4. Observables generated by martingales.

4.1. *Local martingales for continuous LERW.* Suppose $D$ is a finitely connected domain, $z_e\in D$, and $w_0$ is a prime end of $D$. Let $\gamma(t)$, $0\le t<S$, be an $\mathrm{LERW}(D;w_0\to z_e)$ trace. So $\gamma$ is a simple curve in $\widehat{D}$ with $\gamma(0)=w_0$ and $\gamma(t)\in D$ for $0<t<S$. For $0\le t<S$, let $P_t$ be the generalized Poisson kernel in $D\setminus\gamma((0,t])$ with the pole at $\gamma(t)$, normalized by $P_t(z_e)=1$.

THEOREM 4.1. *For any fixed $z\in D$, $(P_t(z))$ is a local martingale.*

Let $\Omega$ be an almost $\mathbb{H}$ domain, and $p\in\Omega$. If $K$ is a hull in $\mathbb{H}$ w.r.t. $\infty$ such that $K\subset\Omega\setminus\{p\}$, let $P(K,x,\cdot)$ be the generalized Poisson kernel in $\Omega_K$ with the pole at $x$, normalized by $P(K,x,\varphi_K(p))=1$. Suppose $\xi\in C([0,T))$ satisfies $\bigcup_{0\le t<T}K_t^\xi\subset\Omega\setminus\{p\}$. We write $P^\xi(t,\cdot,\cdot)$ for $P(K_t^\xi,\cdot,\cdot)$, $t\in[0,T)$. It is standard to check that $P^\xi$ is $C^{1,2,h}$ differentiable, where "$h$" means harmonic.

LEMMA 4.1. *For any $t\in[0,T)$ and $z\in\Omega\setminus K_t^\xi$, we have $\mathcal{V}_t(z)=0$, where*

$$\mathcal{V}_t(z)=\partial_1P^\xi(t,\xi(t),\varphi_t^\xi(z))+2\partial_2P^\xi(t,\xi(t),\varphi_t^\xi(z))X_t^\xi$$
$$+\partial_2^2P^\xi(t,\xi(t),\varphi_t^\xi(z))+2\,\mathrm{Re}\Big(\partial_{3,z}P^\xi(t,\xi(t),\varphi_t^\xi(z))\cdot\frac{2}{\varphi_t^\xi(z)-\xi(t)}\Big).$$

*Here $\partial_1$ and $\partial_2$ are partial derivatives w.r.t. the first two (real) variables, and $\partial_{3,z}=(\partial_{3,x}-i\partial_{3,y})/2$ is the partial derivative w.r.t. the third (complex) variable.*

Proof. For $t\in[0,T)$ and $z\in\partial\Omega\setminus\mathbb{R}$, since $\varphi_t^\xi(z)\in\partial\Omega_t^\xi\setminus\mathbb{R}$, so $P^\xi(t,x,\varphi_t^\xi(z))=0$ for any $x\in\mathbb{R}$, which implies that $\partial_2 P^\xi=\partial_2^2P^\xi=0$ at $(t,x,\varphi_t^\xi(z))$, and

$$\partial_1 P^\xi(t,x,\varphi_t^\xi(z))+2\operatorname{Re}\left(\partial_{3,z}P^\xi(t,x,\varphi_t^\xi(z))\cdot\frac{2}{\varphi_t^\xi(z)-\xi(t)}\right)=0.$$

Thus $\mathcal{V}_t$ vanishes on $\partial\Omega\setminus\mathbb{R}$ for $t\in[0,T)$. Let $\mathcal{W}_t=\mathcal{V}_t\circ(\varphi_t^\xi)^{-1}$. Then $\mathcal{W}_t$ vanishes on $\partial\Omega_t^\xi\setminus\mathbb{R}$ for $t\in[0,T)$. Note that for $t\in[0,T)$ and $w\in\Omega_t^\xi$,

$$\begin{aligned}\mathcal{W}_t(w)&=\partial_1P^\xi(t,\xi(t),w)+2\partial_2P^\xi(t,\xi(t),w)X_t^\xi\\&\quad+\partial_2^2P^\xi(t,\xi(t),w)+2\operatorname{Re}\left(\partial_{3,z}P^\xi(t,\xi(t),w)\cdot\frac{2}{w-\xi(t)}\right).\end{aligned}$$

Since $P^\xi(t,\xi(t),\cdot)$ vanishes on $\mathbb{R}\setminus\{\xi(t)\}$ and $\frac{2}{w-\xi(t)}$ is real on $\mathbb{R}\setminus\{\xi(t)\}$, so $\mathcal{W}_t$ vanishes on $\mathbb{R}\setminus\{\xi(t)\}$. As $w\to\infty$ in $\mathbb{H}$, $\partial_1$, $\partial_2$, $\partial_2^2$ and $\partial_{3,z}$ of $P^\xi$ at $(t,\xi(t),w)$ all tend to 0, and $\frac{2}{w-\xi(t)}$ tends to 0 as well. Thus $\mathcal{W}_t$ vanishes on $\widehat{\mathbb{R}}\setminus\{\xi(t)\}$.

Suppose for some $c(t,x)\in\mathbb{R}$, $\operatorname{Im}\frac{c(t,x)}{w-x}$ is the principal part of $P^\xi(t,x,w)$ at $x$. So there is some analytic function $F(t,x,\cdot)$ defined in some neighborhood of $x$ such that in that neighborhood, $P^\xi(t,x,w)=\operatorname{Im}(F(t,x,w)+\frac{c(t,x)}{w-x})$. Then we have

$$\begin{aligned}\partial_1P^\xi(t,\xi(t),w)&=\operatorname{Im}\left(\partial_1F(t,\xi(t),w)+\frac{\partial_1c(t,\xi(t))}{w-\xi(t)}\right),\\ \partial_2P^\xi(t,\xi(t),w)&=\operatorname{Im}\left(\partial_2F(t,\xi(t),w)+\frac{\partial_2c(t,\xi(t))}{w-\xi(t)}\right.\\&\qquad\left.+\frac{c(t,\xi(t))}{(w-\xi(t))^2}\right),\\ \partial_2^2P^\xi(t,\xi(t),w)&=\operatorname{Im}\left(\partial_2^2F(t,\xi(t),w)+\frac{\partial_2^2c(t,\xi(t))}{w-\xi(t)}\right.\\&\qquad\left.+\frac{2\partial_2c(t,\xi(t))}{(w-\xi(t))^2}+\frac{2c(t,\xi(t))}{(w-\xi(t))^3}\right)\end{aligned}$$

and

$$2\operatorname{Re}\left(\partial_{3,z}P^\xi(t,\xi(t),w)\cdot\frac{2}{w-\xi(t)}\right)=\operatorname{Im}\left(\frac{2F'(t,\xi(t),w)}{w-\xi(t)}-\frac{2c(t,\xi(t))}{(w-\xi(t))^3}\right).$$

Thus $\mathcal{W}_t$ equals the imaginary part of

$$\partial_1F(t,\xi(t),w)+\frac{\partial_1c(t,\xi(t))}{w-\xi(t)}$$

$$+2\Big(\partial_2 F(t,\xi(t),w)+\frac{\partial_2 c(t,\xi(t))}{w-\xi(t)}+\frac{c(t,\xi(t))}{(w-\xi(t))^2}\Big)X^\xi_t$$
$$+\partial_2^2 F(t,\xi(t),w)+\frac{\partial_2^2 c(t,\xi(t))}{w-\xi(t)}+\frac{2\partial_2 c(t,\xi(t))}{(w-\xi(t))^2}+\frac{2c(t,\xi(t))}{(w-\xi(t))^3}$$
$$+\frac{2F'(t,\xi(t),w)}{w-\xi(t)}-\frac{2c(t,\xi(t))}{(w-\xi(t))^3}$$
$$=\partial_1 F(t,\xi(t),w)+2\partial_2 F(t,\xi(t),w)X^\xi_t$$
$$+\partial_2^2 F(t,\xi(t),w)+\frac{A_1(t)}{w-\xi(t)}+\frac{A_2(t)}{(w-\xi(t))^2}$$

for some functions $A_1(t)$ and $A_2(t)$, where $A_2(t)=2c(t,\xi(t))X^\xi_t+2\partial_2 c(t,\xi(t))$.

Since $J^\xi_t=G(\Omega^\xi_t,\varphi^\xi_t(p);\cdot)$, so for $x\in\mathbb{R}$, $\partial_y J^\xi_t(x)$ equals the value at $\varphi^\xi_t(p)$ of the (usual) Poisson kernel in $\Omega^\xi_t$ with the pole at $x$. Note that $P^\xi(t,x,\cdot)$ equals some constant times the Poisson kernel in $\Omega^\xi_t$ with the pole at $x$, of which the principal part at $x$ is $\operatorname{Im}\frac{-1/\pi}{w-x}$. So we have

$$\partial_y J^\xi_t(x)/(-1/\pi)=P^\xi(t,x,\varphi^\xi_t(p))/c(t,x)=1/c(t,x).$$

Thus $c(t,x)\,\partial_y J^\xi_t(x)=-1/\pi$ for any $x\in\mathbb{R}$, which implies that

$$0=c(t,\xi(t))\,\partial_x\,\partial_y J^\xi_t(\xi(t))+\partial_2 c(t,\xi(t))\partial_y J^\xi_t(\xi(t))=A_2(t)\,\partial_y J^\xi_t(\xi(t))/2.$$

So $A_2(t)=0$, and $\mathcal{W}_t$ equals the imaginary part of some analytic function plus $\frac{A_1(t)}{w-\xi(t)}$ near $\xi(t)$. Since $\mathcal{W}_t$ is harmonic in $\Omega^\xi_t$, and vanishes at every prime end of $\Omega^\xi_t$ other than $\xi(t)$, so $\mathcal{W}_t=C(t)P^\xi(t,\xi(t),\cdot)$ for some $C(t)\in\mathbb{R}$. From $P^\xi(t,x,\varphi^\xi_t(p))=1$ for any $t\in[0,T)$ and $x\in\mathbb{R}$, we get $\mathcal{W}_t(\varphi^\xi_t(p))=0$. So for $t\in[0,T)$, we have $C(t)=0$, which implies that $\mathcal{W}_t$ vanishes on $\Omega^\xi_t$, and so $\mathcal{V}_t$ vanishes on $\Omega\setminus K^\xi_t$. □

Suppose $f$ maps $D$ conformally onto an almost $\mathbb{H}$ domain $\Omega$ such that $f(w_0)=0$. Let $p=f(z_e)$. Let $v(t)=\operatorname{hcap}(f(\gamma((0,t])))/2$, $0\le t<S$. Let $T=v(S)$, and $u$ be the reversal of $v$. Then $f(\gamma((0,u(t)]))=K^\xi_t$, $0\le t<T$, where $\xi\in C([0,T))$ solves equation (3.2) with $\lambda=2$ and $A(t)=\sqrt{2}B(t)$ for some Brownian motion $B(t)$. Since $\varphi^\xi_t\circ f$ maps $D\setminus f(\gamma((0,u(t)]))$ conformally onto $\Omega^\xi_t$, $\varphi^\xi_t\circ f(\gamma(u(t)))=\xi(t)$ and $\varphi^\xi_t\circ f(z_e)=\varphi^\xi_t(p)$, so from the conformal invariance, $P_{u(t)}\circ f^{-1}\circ(\varphi^\xi_t)^{-1}$ is the generalized Poisson kernel in $\Omega^\xi_t$ with the pole at $\xi(t)$, whose value at $\varphi^\xi_t(p)$ is 1, that is,

$$P_{u(t)}\circ f^{-1}\circ(\varphi^\xi_t)^{-1}=P^\xi(t,\xi(t),\cdot). \tag{4.1}$$

Proof of Theorem 4.1. Let $Q_t(z) = P^\xi(t, \xi(t), \varphi_t^\xi(z))$ for $z \in \Omega \setminus K_t^\xi$. From Itô's formula, $(Q_t(z))$ is a semimartingale, and the drift term equals $\mathcal{V}_t(z)$, which vanishes on $\Omega \setminus K_t^\xi$ by Lemma 4.1. Thus $(Q_t(z))$ is a local martingale for any fixed $z \in \Omega$. From (4.1), $P_t(z) = Q_{v(t)}(f(z))$ for $z \in D$. Since $f(D) = \Omega$, and a time-change preserves a local martingale, so $(P_t(z))$ is a local martingale for any fixed $z \in D$.

Second, we consider an LERW$(D; w_0 \to I_e)$ trace: $\gamma(t)$, $0 \le t < S$, where $w_0$ is a prime end of $D$, and $I_e$ is a side arc of $D$. Let $P_t$ be the generalized Poisson kernel in $D \setminus \gamma((0,t])$ with the pole at $\gamma(t)$, normalized by $\int_{I_e} \partial_{\mathbf{n}} P_t(z)\, ds(z) = 1$. Here the equality means that if $g$ maps a neighborhood $U$ of $I_e$ in $\widehat{D}$ conformally into $\mathbb{C}$ such that $g(I_e)$ is an analytic arc, then $\int_{g(I_e)} \partial_{\mathbf{n}} (P_t \circ g^{-1})(z)\, ds(z) = 1$, where $\mathbf{n}$ is the unit normal vector pointing inward, and $ds$ is the length of the curve. In fact, the value of the integral does not depend on the choice of $g$.

Suppose $f$ maps $D$ conformally onto an almost $\mathbb{H}$ domain $\Omega$ such that $f(w_0) = 0$. Let $J = f(I_e)$. If $K_t^\xi \subset \Omega$, and is bounded away from $J$, let $P^\xi(t, x, \cdot)$ be the generalized Poisson kernel in $\Omega_t^\xi$, with the pole at $x$, normalized by $\int_{\varphi_t^\xi(J)} \partial_{\mathbf{n}} P^\xi(t, x, z)\, ds(z) = 1$. Then Lemma 4.1 holds in this setting, and the proof is similar. Formula (4.1) still holds, so we have Theorem 4.1.

Third, we consider an LERW$(D; w_0 \to w_e)$ trace: $\gamma(t)$, $0 \le t < S$, where $w_0 \ne w_e$ are prime ends of $D$. Fix $g$ that maps a neighborhood $U$ of $w_e$ in $\widehat{D}$ conformally onto a neighborhood $V$ of $0$ in $\overline{\mathbb{H}}$ such that $g(w_e) = 0$ and $g(U \cap \widehat{\partial} D) \subset \mathbb{R}$. Let $P_t$ be the generalized Poisson kernel in $D \setminus \gamma((0,t])$ with the pole at $\gamma(t)$, normalized by $\partial_y (P_t \circ g^{-1})(0) = 1$.

Suppose $f$ maps $D$ conformally onto an almost $\mathbb{H}$ domain $\Omega$ such that $f(w_0) = 0$. Let $p = f(w_e)$. If $K_t^\xi \subset \Omega$, and is bounded away from $p$, let $P^\xi(t, x, \cdot)$ be the generalized Poisson kernel in $\Omega_t^\xi$, with the pole at $x$, normalized by $\partial_y (P^\xi(t, x, \cdot) \circ f \circ g^{-1})(0) = 1$. Then Lemma 4.1 holds in this setting, and the proof is similar. Formula (4.1) still holds, so we have Theorem 4.1. □

4.2. *Discrete approximations.* Let $D$ be a finitely connected domain. Suppose $0 \in \partial D$, and there is some $\delta_D > 0$ such that the half open line segment $[\delta_D, 0)$ is contained in $D$. As $z \to 0$ along $[\delta_D, 0)$, $z$ tends to a prime end of $D$. We use $0_+$ to denote this prime end.

For $\delta > 0$, let $\delta\mathbb{Z}^2 = \{(j + ik)\delta : j, k \in \mathbb{Z}\} \subset \mathbb{C}$. We also view $\delta\mathbb{Z}^2$ as a graph whose vertices are $(j + ik)\delta$, $j, k \in \mathbb{Z}$, and two vertices are adjacent iff the distance between them is $\delta$. We define a graph $\breve{D}^\delta$ that approximates $D$ in $\delta\mathbb{Z}^2$ as follows. The vertex set $V(\breve{D}^\delta)$ is the union of interior vertex set $V_I(\breve{D}^\delta)$ and boundary vertex set $V_\partial(\breve{D}^\delta)$, where $V_I(\breve{D}^\delta) := \delta\mathbb{Z}^2 \cap D$, and $V_\partial(\breve{D}^\delta)$ is the set of ordered pairs $\langle z_1, z_2 \rangle$ such that $z_1 \in V_I(\breve{D}^\delta)$, $z_2 \in \partial D$, and there is $z_3 \in \delta\mathbb{Z}^2$ that is adjacent to $z_1$ in $\delta\mathbb{Z}^2$, such that $[z_1, z_2) \subset [z_1, z_3) \cap D$. Two vertices $w_1$ and $w_2$ in $V(\breve{D}^\delta)$ are adjacent iff either $w_1, w_2 \in V_I(\breve{D}^\delta)$, $w_1$ and $w_2$ are adjacent in $\delta\mathbb{Z}^2$, and $[w_1, w_2] \subset D$; or for $j = 1$ or 2, $w_j \in V_I(\breve{D}^\delta)$ and $w_{3-j} = \langle w_j, z_3 \rangle \in V_\partial(\breve{D}^\delta)$ for some $z_3 \in \partial D$.

Every interior vertex of $\breve{D}^\delta$ has exactly four adjacent vertices, and every boundary vertex $w = \langle z_1, z_2 \rangle$ has exactly one adjacent vertex, which is the interior vertex $z_1$. So $\breve{D}^\delta$ is locally finite. If $\langle z_1, z_2 \rangle$ is a boundary vertex, then it determines a boundary point, which is $z_2$, and a prime end of $D$, which is the limit in $\widehat{D}$ as $z \to z_2$ along $[z_1, z_2)$. If there is no ambiguity, we do not distinguish a boundary vertex from the boundary point or prime end it determines. Suppose $\delta \in (0, \delta_D]$. Then $\delta$ is an interior vertex of $\breve{D}^\delta$, and $\langle \delta, 0 \rangle$ is a boundary vertex of $\breve{D}^\delta$. A random walk on $\breve{D}^\delta$ started from an interior vertex $w_0$ up to the first time it leaves $D$ agrees with a random walk on $\delta\mathbb{Z}^2$ started from $w_0$ up to the first time it uses an edge that intersects $\partial D$. Let $D^\delta$ be the connected component of $\breve{D}^\delta$ that contains $\delta$. Let $V_I(D^\delta) := V(D^\delta) \cap V_I(\breve{D}^\delta)$ and $V_\partial(D^\delta) := V(D^\delta) \cap V_\partial(\breve{D}^\delta)$ be the set of interior and boundary vertices, respectively, of $D^\delta$.

Fix $z_e \in D \setminus \{\infty\}$. Let $w_e^\delta$ be the vertex in $\delta\mathbb{Z}^2$ that is closest to $z_e$. If such vertex is not unique, we choose the one that maximizes $\operatorname{Re} z + \pi \operatorname{Im} z$ to break the tie. Suppose $\delta \in (0, \delta_D]$ is small enough. Then there is a lattice path on $\breve{D}^\delta$ that connects $\delta$ with $w_e^\delta$, which does not pass through any boundary vertex. So $w_e^\delta$ is an interior vertex of $D^\delta$. Let $F = \{w_e^\delta\}$ and $E_{-1} = V_\partial(D^\delta)$. From the recurrence of the random walks on $\mathbb{Z}^2$, we know that $E \cup F$ is reachable in $D^\delta$. Let $(q_\delta(0), \ldots, q_\delta(\chi_\delta))$ be the LERW on $D^\delta$ started from $\delta$ conditioned to hit $F$ before $E_{-1}$. So $q_\delta(0) = \delta$ and $q_\delta(\chi_\delta) = w_e^\delta$. Let $q_\delta(-1) = 0$. Extend $q_\delta$ to $[-1, \chi_\delta]$ such that $q_\delta$ is linear on $[k-1, k]$ for each $k \in \mathbb{Z}_{[0,\chi_\delta]}$. Then $q_\delta$ is a simple curve in $D \cup \{0\}$ that connects 0 and $w_e^\delta$.

Since $F$ contains only one point, we may define $g_k$ as in Proposition 2.1. Then for any fixed vertex $v_0$ on $D^\delta$, $(g_k(v_0))$ is a martingale up to the time $q_\delta(k)$ is next to $w_e^\delta$ or $E_k := E_{-1} \cup \{q_\delta(0), \ldots, q_\delta(k)\}$ disconnects $v_0$ from $z_e$. Note that $g_k$ vanishes on $E_k \setminus \{q_\delta(k)\}$, is discrete harmonic at every interior vertex of $D^\delta$ except $q_\delta(0), \ldots, q_\delta(k)$, and $g_k(w_e^\delta) = 1$. For $0 \le k \le \chi_\delta - 1$, let $D_k = D \setminus q_\delta([-1, k])$. Then $q_\delta(k)$ corresponds to a prime end of $D_k$. When $\delta$ is small, the function $g_k$ approximates the generalized Poisson kernel $P_k$ in $D_k$ with the pole at $q_\delta(k)$, normalized by $P_k(z_e) = 1$. Note the resemblance of the discrete martingales preserved by (discrete) LERW and the local martingales preserved by continuous LERW. Suppose $\gamma_0(t)$, $0 \le t < S_0$, is an LERW$(D; 0_+ \to z_e)$ trace. In the last several sections, we will prove the following theorem. Note that we do not require that the boundary of $D$ is good.

THEOREM 4.2. (i) *Suppose $U$ is a neighborhood of $0_+$ in $D$. Then for any $\varepsilon > 0$, there is $\delta_0 > 0$ such that if $\delta < \delta_0$, then there are a coupling of*

*$q_\delta$ and $\gamma_0$, and a continuous increasing function $\breve{u}$ that maps $(-1,\chi_\delta)$ onto $(0,S_0)$ such that*

$$\mathbf{P}[\sup\{|q_\delta(\breve{u}^{-1}(t))-\gamma_0(t)|: T_U(\gamma_0)\le t<S_0\}<\varepsilon]>1-\varepsilon,$$

*where $T_U(\gamma_0)$ is the first time that $\gamma_0$ leaves $U$.*

(ii) *If the prime end $0_+$ is degenerate (see [13]), then* (i) *holds with "$T_U(\gamma_0)\le t$" replaced by "$0<t$."*

Now suppose $w_e\in\partial D\setminus\{0\}$ satisfies $w_e\in\delta_e\mathbb{Z}^2$ for some $\delta_e>0$, and $\partial D$ is flat near $w_e$, which means that there is $r>0$ such that $D\cap\{z\in\mathbb{C}:|z-w_e|<r\}=(w_e+a\mathbb{H})\cap\{z\in\mathbb{C}:|z-w_e|<r\}$ for some $a\in\{\pm1,\pm i\}$. For $\delta>0$, let $w_e^\delta=w_e+ia\delta$.

Let $\mathcal{M}$ be the set of $\delta>0$ such that $w_e\in\delta\mathbb{Z}^2$. If $\delta\in\mathcal{M}$ is small enough, then $\langle w_e^\delta,w_e\rangle$ is a boundary vertex of $\breve{D}^\delta$, which determines the boundary point and prime end $w_e$, and there is a lattice path on $D^\delta$ that connects $\delta$ with $w_e$ without passing through any other boundary vertex. Here we do not distinguish $w_e$ from the boundary vertex $\langle w_e^\delta,w_e\rangle$. Let $F=\{w_e\}$ and $E_{-1}=V_\partial(D^\delta)\setminus F$. Then $E\cup F=V_\partial(D^\delta)$ is reachable in $D^\delta$. Let $(q_\delta(0),\dots,q_\delta(\chi_\delta))$ be the LERW on $D^\delta$ started from $\delta$ conditioned to hit $F$ before $E_{-1}$. So $q_\delta(0)=\delta$ and $q_\delta(\chi_\delta)=w_e$. Let $q_\delta(-1)=0$. Extend $q_\delta$ to be defined on $[-1,\chi_\delta]$ such that $q_\delta$ is linear on $[k-1,k]$ for each $k\in\mathbb{Z}_{[0,\chi_\delta]}$. Then $q_\delta$ is a simple curve in $D\cup\{0,w_e\}$ that connects 0 and $w_e$.

Let $h_k$ be as in Proposition 2.1. Then for any fixed vertex $v_0$ on $D^\delta$, $(h_k(v_0))$ is a martingale up to the time when $q_\delta(k)=w_e^\delta$ or $E_k=E_{-1}\cup\{q_\delta(0),\dots,q_\delta(k)\}$ disconnects $v_0$ from $w_e$. Let $D_k=D\setminus q_\delta([-1,k])$. Then $q_\delta(k)$ is a prime end of $D_k$. Note that $h_k$ vanishes on $q_\delta(-1),\dots,q_\delta(k-1)$ and all boundary vertices of $D^\delta$, is discrete harmonic at all interior vertices of $D^\delta$ except $q_\delta(0),\dots,q_\delta(k)$, and $h_k(w_e^\delta)=1$. So when $\delta$ is small, $\delta\cdot h_k$ is close to the generalized Poisson kernel $P_k$ in $D_k$ with the pole at $q_\delta(k)$ normalized by $\partial_{\mathbf{n}}P_k(w_e)=1$. Suppose $\gamma_0(t)$, $0\le t<S$, is an LERW$(D;0_+\to w_e)$ trace. Then Theorem 4.2 still holds for $q_\delta$ and $\gamma_0$ defined here if we replace "$\delta<\delta_0$" by "$\delta\in\mathcal{M}$ and $\delta<\delta_0$."

Now suppose $I_e$ is a side arc of $D$ that is bounded away from $0_+$. Let $I_e^\delta$ be the set of boundary vertices of $D^\delta$ which determine prime ends that lie on $I_e$. If $\delta$ is small enough, $I_e^\delta$ is nonempty, and there is a lattice path on $D^\delta$ that connects $\delta$ with $I_e^\delta$ without passing through any boundary vertex not in $I_e^\delta$. Then we let $F=I_e^\delta$ and $E_{-1}=V_\partial(D^\delta)\setminus F$. Let $(q_\delta(0),\dots,q_\delta(\chi_\delta))$ be the LERW on $D^\delta$ started from $\delta$ conditioned to hit $F$ before $E_{-1}$. So $q_\delta(0)=\delta$ and $q_\delta(\chi_\delta)\in I_e$.

Let $h_k$ be as in Proposition 2.1. Then for any fixed vertex $v_0$ on $D^\delta$, $(h_k(v_0))$ is a martingale up to the time $q_\delta(k)$ is close to $I_e$ or $E_k:=E_{-1}\cup\{q_\delta(0),\dots,q_\delta(k)\}$ disconnects $v_0$ from $I_e$. Note that $h_k$ vanishes on $q_\delta(-1),\dots,q_\delta(k-1)$ and all boundary vertices of $D^\delta$, $h_k$ is discrete harmonic at every interior vertex of $D^\delta$ except $q_\delta(0),\dots,q_\delta(k)$, and $\sum_{v\in I_e^\delta}\Delta h_k(v)=1$. So when $\delta$ is small, the function $h_k$ seems to be close to the generalized Poisson kernel $P_k$ in $D_k$ with the pole at $q_\delta(k)$ normalized by $\int_{I_e}\partial_{\mathbf{n}}P_k(z)\,ds(z)=1$.

If $I_e$ is a whole side of $D$, then Theorem 4.2 still holds for $q_\delta$ and $\gamma_0$ defined here. If $I_e$ is not a whole side, for the purpose of convergence, we may need some additional boundary conditions. Suppose the two ends of $I_e$ correspond to $w_e^1,w_e^2\in\partial D$, near which $\partial D$ is flat, and $w_e^1,w_e^2\in\delta_e\mathbb{Z}^2$ for some $\delta_e>0$. Let $\mathcal{M}$ be the set of $\delta>0$ such that $w_e^1,w_e^2\in\delta\mathbb{Z}^2$. Then Theorem 4.2 still holds for $q_\delta$ and $\gamma_0$ defined here if we replace "$\delta<\delta_0$" by "$\delta\in\mathcal{M}$ and $\delta<\delta_0$."

**5. Existence and uniqueness.** In this section we will prove Theorem 3.1. The proof is somehow similar to that of the existence and uniqueness of the solution of an ordinary differential equation.

5.1. *Convergence of domains.*

DEFINITION 5.1. Suppose $D_n$ is a sequence of domains and $D$ is a domain. We say that $(D_n)$ converges to $D$, denoted by $D_n\xrightarrow{\text{Cara}}D$, if for every $z\in D$, $\text{dist}^\#(z,\partial^\# D_n)\to\text{dist}^\#(z,\partial^\# D)$. This is equivalent to the followings:

(i) every compact subset of $D$ is contained in all but finitely many $D_n$'s; and
(ii) for every point $z_0\in\partial^\# D$, $\text{dist}^\#(z_0,\partial^\# D_n)\to 0$ as $n\to\infty$.

A sequence of domains may converge to two different domains. For example, let $D_n=\mathbb{C}\setminus((-\infty,n])$. Then $D_n\xrightarrow{\text{Cara}}\mathbb{H}$, and $D_n\xrightarrow{\text{Cara}}-\mathbb{H}$ as well. But two different limit domains of the same domain sequence must be disjoint from each other, because if they have nonempty intersection, then one contains some boundary point of the other, which implies a contradiction.

If only condition (i) in the definition is satisfied, then for any $z\in D$, $\text{dist}^\#(z,\partial^\# D)\le\liminf\text{dist}^\#(z,\partial^\# D_n)$. Thus $D_n\cap D\xrightarrow{\text{Cara}}D$. If $D_n\xrightarrow{\text{Cara}}D$, $E_n\xrightarrow{\text{Cara}}E$, and $z_0\in D\cap E$. Let $F_n$ (resp. $F$) be the connected component of $D_n\cap E_n$ (resp. $D\cap E$) that contains $z_0$. Then for any $z\in F$, $\text{dist}^\#(z,\partial^\# F_n)=\text{dist}^\#(z,\partial^\# D_n)\wedge\text{dist}^\#(z,\partial^\# E_n)$ for each $n$, and $\text{dist}^\#(z,\partial^\# F)=\text{dist}^\#(z,\partial^\# D)\wedge\text{dist}^\#(z,\partial^\# E)$, which implies $F_n\xrightarrow{\text{Cara}}F$. Thus if $D_n\xrightarrow{\text{Cara}}D$, $E_n\xrightarrow{\text{Cara}}E$, $D_n\subset E_n$ for each $n$, and $D\cap E\ne\varnothing$, then we have $D\subset E$.

Suppose $D_n\xrightarrow{\text{Cara}}D$, and for each $n$, $f_n$ is a $\widehat{\mathbb{C}}$-valued function on $D_n$, and $f$ is a $\widehat{\mathbb{C}}$-valued function on $D$. We say that $f_n$ converges to $f$ locally uniformly in $D$, or $f_n\xrightarrow{\text{l.u.}}f$ in $D$, if for each compact subset $F$ of $D$, $f_n$ converges to $f$ in the spherical metric uniformly on $F$. If every $f_n$ is analytic (resp. harmonic), then $f$ is also analytic (resp. harmonic).

LEMMA 5.1. *Suppose* $D_n \xrightarrow{\text{Cara}} D$, $f_n$ *maps* $D_n$ *conformally onto some domain* $E_n$ *for each* $n$, *and* $f_n \xrightarrow{\text{l.u.}} f$ *in* $D$. *Then either* $f$ *is constant on* $D$, *or* $f$ *maps* $D$ *conformally onto some domain* $E$. *And in the latter case,* $E_n \xrightarrow{\text{Cara}} E$ *and* $f_n^{-1} \xrightarrow{\text{l.u.}} f^{-1}$ *in* $E$.

This lemma is similar to Theorem 1.8, the Carathéodory kernel theorem, in [13], and the proof is also similar. When applying this lemma, we will usually first exclude the possibility that $f$ is constant, and then obtain the convergence of the image domains and the inverse functions.

5.2. *Topology on the space of hulls.* If $H$ is a nonempty hull in $\mathbb{H}$ w.r.t. $\infty$, then $\overline{H}\cap\mathbb{R}$ is nonempty. Let $a_H=\inf(\overline{H}\cap\mathbb{R})$ and $b_H=\sup(\overline{H}\cap\mathbb{R})$. Let

$$\Sigma_H=\mathbb{C}\setminus(H\cup\{\overline{z}:z\in H\}\cup[a_H,b_H]).$$

By the reflection principle, $\varphi_H$ extends to $\Sigma_H$, and maps $\Sigma_H$ conformally onto $\mathbb{C}\setminus[c_H,d_H]$ for some $c_H<d_H\in\mathbb{R}$. Moreover, $\varphi_H$ is increasing on $(-\infty,a_H)$ and $(b_H,+\infty)$, and maps them onto $(-\infty,c_H)$ and $(d_H,+\infty)$, respectively. So $\varphi_H^{-1}$ extends conformally to $\mathbb{C}\setminus[c_H,d_H]$. And $[c_H,d_H]$ is the smallest in the sense that if $\varphi_H^{-1}$ extends conformally to $\mathbb{C}\setminus I$ for some closed interval $I$, then $[c_H,d_H]\subset I$. If $H=\varnothing$, we do not define $a_H,b_H,c_H,d_H$, but still use the notation $[a_\varnothing,b_\varnothing]$ and $[c_\varnothing,d_\varnothing]$ to denote empty sets. Then $\Sigma_\varnothing=\mathbb{C}$, so it is true that $\varphi_\varnothing$ maps $\Sigma_\varnothing$ conformally onto $\mathbb{C}\setminus[c_\varnothing,d_\varnothing]$.

If $\gamma$ is a crosscut in $\mathbb{H}$, we define $H(\gamma)$ to be $\gamma$ unions the bounded component of $\mathbb{H}\setminus\gamma$. Then $H(\gamma)$ is a hull in $\mathbb{H}$ w.r.t. $\infty$. We call it the hull bounded by $\gamma$. If $A\subset\overline{H(\gamma)}$, then we say $\gamma$ encloses $A$. If $A\subset\overline{H(\gamma)}$ and $\overline{A}\cap\overline{\gamma}=\varnothing$, then we say $\gamma$ strictly encloses $A$. For simplicity, we write $x_\gamma$ instead of $x_{H(\gamma)}$ when $x$ is one of the following symbols: $a,b,c,d,\Sigma,\varphi$.

Since $\widehat{\partial}(\mathbb{H}\setminus H(\gamma))=(\widehat{\mathbb{R}}\setminus(a_\gamma,b_\gamma))\cup\gamma$ is a simple curve, so $\varphi_\gamma$ extends to a homeomorphism of $\overline{\mathbb{H}\setminus H(\gamma)}$, and maps $\gamma$ onto $[c_\gamma,d_\gamma]$. So $\varphi_{H(\gamma)}^{-1}$ has a continuous extension to $\mathbb{H}\cup\mathbb{R}$, and maps $(c_\gamma,d_\gamma)$ onto $\gamma$. From the results about Poisson kernel, we have

$$\varphi_\gamma^{-1}(z)-z=\int_{c_\gamma}^{d_\gamma}\frac{-1}{z-x}\,\frac{\operatorname{Im}\varphi_\gamma^{-1}(x)}{\pi}\,dx,$$

for any $z\in\Sigma_\gamma$. From the behavior of $\varphi_\gamma$ near $\infty$, we have $\int_{c_\gamma}^{d_\gamma}\operatorname{Im}\varphi_\gamma^{-1}(x)/\pi\,dx=\operatorname{hcap}(H(\gamma))$. If $H$ is a general nonempty hull in $\mathbb{H}$ w.r.t. $\infty$, then $\varphi_H^{-1}$ may not have continuous extension to $[c_H,d_H]$. We may use a sequence of hulls bounded by crosscuts to approximate $H$. Then we conclude that there is a positive measure $\mu_H$ supported by $[c_H,d_H]$ with total mass $|\mu_H|=\operatorname{hcap}(H)$ such that for any $z\in\Sigma_H$,

$$\varphi_H^{-1}(z)-z=\int_{c_H}^{d_H}\frac{-1}{z-x}\,d\mu_H(x). \tag{5.1}$$

EXAMPLE. Suppose $x_0 \in \mathbb{R}$ and $r_0 > 0$. Let $\alpha = \{z \in \mathbb{H} : |z - x_0| = r_0\}$. Then $\alpha$ is a crosscut in $\mathbb{H}$, $H(\alpha) = \{z \in \mathbb{H} : |z - x_0| \le r_0\}$ and $[a_\alpha, b_\alpha] = [x_0 - r_0, x_0 + r_0]$. It is clear that $\varphi_\alpha(z) = z + \frac{r_0^2}{z - x_0}$. Thus $\operatorname{hcap}(H(\alpha)) = r_0^2$ and $[c_\alpha, d_\alpha] = [x_0 - 2r_0, x_0 + 2r_0]$.

LEMMA 5.2. *If $H$ is a nonempty hull in $\mathbb{H}$ w.r.t. $\infty$, then $\varphi_H^{-1}(x) > x$ for any $x \in (-\infty, c_H)$; $\varphi_H^{-1}(x) < x$ for any $x \in (d_H, +\infty)$; $\varphi_H(x) < x$ for any $x \in (-\infty, a_H)$; $\varphi_H(x) > x$ for any $x \in (b_H, +\infty)$. So if $H$ is any hull in $\mathbb{H}$ w.r.t. $\infty$, then $[a_H, b_H] \subset [c_H, d_H]$.*

PROOF. This follows from (5.1) and that $\varphi_H$ maps $(-\infty, a_H)$ and $(b_H, +\infty)$ onto $(-\infty, c_H)$ and $(d_H, +\infty)$, respectively. $\square$

If $H_1 \subset H_2$ are two hulls in $\mathbb{H}$ w.r.t. $\infty$, we call $H_1$ a sub-hull of $H_2$. Then $H_2/H_1 := \varphi_{H_1}(H_2 \setminus H_1)$ is also a hull in $\mathbb{H}$ w.r.t. $\infty$. We call $H_2/H_1$ a quotient-hull of $H_2$. It is clear that $\varphi_{H_2} = \varphi_{H_2/H_1} \circ \varphi_{H_1}$. Thus $\operatorname{hcap}(H_2) = \operatorname{hcap}(H_2/H_1) + \operatorname{hcap}(H_1)$, and so $\operatorname{hcap}(H_1), \operatorname{hcap}(H_2/H_1) \le \operatorname{hcap}(H_2)$.

LEMMA 5.3. *If $H_1 \subset H_2$ are two hulls in $\mathbb{H}$ w.r.t. $\infty$, then $[c_{H_1}, d_{H_1}] \subset [c_{H_2}, d_{H_2}]$ and $[c_{H_2/H_1}, d_{H_2/H_1}] \subset [c_{H_2}, d_{H_2}]$.*

PROOF. If $H_1 = \varnothing$ or $H_1 = H_2$, then $H_2/H_1 = H_2$ or $H_2/H_1 = \varnothing$, so it is trivial. Now suppose $\varnothing \subsetneqq H_1 \subsetneqq H_2$. Then $H_2/H_1 \ne \varnothing$. Since $\varphi_{H_2/H_1}^{-1}(z) = \varphi_{H_1} \circ \varphi_{H_2}^{-1}(z)$ for $z \in \mathbb{H}$, $\varphi_{H_2}^{-1}$ maps $\mathbb{C} \setminus [c_{H_2}, d_{H_2}]$ onto $\Sigma_{H_2}$, and $\varphi_{H_1}$ extends conformally to $\Sigma_{H_1} \supset \Sigma_{H_2}$, so $\varphi_{H_2/H_1}^{-1}$ extends conformally to $\mathbb{C} \setminus [c_{H_2}, d_{H_2}]$. From the minimum property of $[c_{H_2/H_1}, d_{H_2/H_1}]$, we have $[c_{H_2/H_1}, d_{H_2/H_1}] \subset [c_{H_2}, d_{H_2}]$.

If $x \in (-\infty, a_{H_2})$, then $\varphi_{H_2}(x) \in (-\infty, c_{H_2}) \subset (-\infty, c_{H_2/H_1})$. Since $\varphi_{H_2/H_1}^{-1}(x) > x$ on $(-\infty, c_{H_2/H_1})$, so $\varphi_{H_1}(x) = \varphi_{H_2/H_1}^{-1} \circ \varphi_{H_2}(x) > \varphi_{H_2}(x)$ on $(-\infty, a_{H_2})$. Thus

$$\begin{aligned} c_{H_1} &= \sup \varphi_{H_1}((-\infty, a_{H_1})) \ge \sup \varphi_{H_1}((-\infty, a_{H_2})) \\ &\ge \sup \varphi_{H_2}((-\infty, a_{H_2})) = c_{H_2}. \end{aligned}$$

Similarly, we have $d_{H_1} \le d_{H_2}$. Thus $[c_{H_1}, d_{H_1}] \subset [c_{H_2}, d_{H_2}]$. $\square$

COROLLARY 5.1. *If $H_1 \subset H_2 \subset H_3$ are hulls in $\mathbb{H}$ w.r.t. $\infty$, then $\operatorname{hcap}(H_2/H_1) \le \operatorname{hcap}(H_3)$ and $[c_{H_2/H_1}, d_{H_2/H_1}] \subset [c_{H_3}, d_{H_3}]$. We call $H_2/H_1$ a sub-quotient-hull of $H_3$.*

Let $H$ be a nonempty hull in $\mathbb{H}$ w.r.t. $\infty$. Let $\mathcal{H}(H)$ denote the set of all sub-hulls of $H$. Let $\mathcal{H}_{\mathrm{sq}}(H)$ denote the set of all sub-quotient-hulls of $H$. If $\alpha$ is a crosscut in $\mathbb{H}$, we write $\mathcal{H}(\alpha)$ for $\mathcal{H}(H(\alpha))$, and $\mathcal{H}_{\rm sq}(\alpha)$ for $\mathcal{H}_{\rm sq}(H(\alpha))$. Choose $d>0$. Let $\alpha=\{z\in\mathbb{C}: |z-(c_H+d_H)/2|=|d_H-c_H|/2+d\}$. Then $\alpha$ is a Jordan curve that encloses $[c_H,d_H]$, and $d$ is the distance between $\alpha$ and $[c_H,d_H]$. Suppose $K\in\mathcal{H}_{\rm sq}(H)$. Then $[c_K,d_K]\subset[c_H,d_H]$. If $z$ lies on or outside $\alpha$, from (5.1),

$$|\varphi_K^{-1}(z)-z|\le|\mu_K|/d=\operatorname{hcap}(K)/d\le\operatorname{hcap}(H)/d.$$

If $z\in\mathbb{C}\setminus[c_K,d_K]$ lies inside $\alpha$, then $\varphi_K^{-1}(z)$ lies inside $\varphi_K^{-1}(\alpha)$. Choose $w\in\alpha$; then

$$\begin{aligned}|\varphi_K^{-1}(z)-z|&\le|z-w|+|w-\varphi_K^{-1}(w)|+|\varphi_K^{-1}(w)-\varphi_K^{-1}(z)|\\&\le\operatorname{diam}(\alpha)+\operatorname{hcap}(H)/d+\operatorname{diam}(\varphi_K^{-1}(\alpha))\\&\le 2|d_H-c_H|+4d+3\operatorname{hcap}(H)/d.\end{aligned}$$

Let $d=\sqrt{\operatorname{hcap}(H)}$ and $M_H=2|d_H-c_H|+7\sqrt{\operatorname{hcap}(H)}$. Then for any $z\in\mathbb{C}\setminus[c_K,d_K]$, $|\varphi_K^{-1}(z)-z|\le M_H$. Since $\varphi_K^{-1}$ maps $\mathbb{C}\setminus[c_K,d_K]$ onto $\Sigma_K$, so for any $z\in\Sigma_K$, $|\varphi_K(z)-z|\le M_H$. Since $\mathbb{C}\setminus[c_K,d_K]\supset\mathbb{C}\setminus[c_H,d_H]$, so $\{\varphi_K^{-1}(z)-z: K\in\mathcal{H}_{\rm sq}(H)\}$ is uniformly bounded in $\mathbb{C}\setminus[c_H,d_H]$ by $M_H$, and so is a normal family.

Let $\mathcal{H}$ denote the set of all hulls in $\mathbb{H}$ w.r.t. $\infty$. Choose a sequence of compact subsets $(F_n)$ of $\mathbb{H}$ such that $F_n\subset\operatorname{int}F_{n+1}$ for each $n\in\mathbb{N}$, and $\bigcup_n F_n=\mathbb{H}$. We may define a distant function $d_{\mathcal{H}}$ on $\mathcal{H}$ such that

$$d_{\mathcal{H}}(H_1,H_2)=\sum_{n=1}^\infty\frac{1}{2^n}\Big(1\wedge\sup_{z\in F_n}\{|\varphi_{H_1}^{-1}(z)-\varphi_{H_2}^{-1}(z)|\}\Big).$$

We use $\xrightarrow{\mathcal{H}}$ to denote the convergence w.r.t. $d_{\mathcal{H}}$. It is clear that $H_n\xrightarrow{\mathcal{H}}H$ iff $\varphi_{H_n}^{-1}\xrightarrow{\rm l.u.}\varphi_H^{-1}$ in $\mathbb{H}$. So the topology does not depend on the choice of $(F_n)$.

From Lemma 5.1, if $H_n\xrightarrow{\mathcal{H}}H$, then $\mathbb{H}\setminus H_n\xrightarrow{\rm Cara}\mathbb{H}\setminus H$ and $\varphi_{H_n}\xrightarrow{\rm l.u.}\varphi_H$ in $\mathbb{H}\setminus H$. However, $\mathbb{H}\setminus H_n\xrightarrow{\rm Cara}\mathbb{H}\setminus H$ does not imply $H_n\xrightarrow{\mathcal{H}}H$. For example, let $H_n=\{z\in\mathbb{H}: |z-2n|\le n\}$ for $n\in\mathbb{N}$. Then $\mathbb{H}\setminus H_n\xrightarrow{\rm Cara}\mathbb{H}=\mathbb{H}\setminus\varnothing$, but $\varphi_{H_n}(z)=z+n^2/(z-2n)\not\to z=\varphi_\varnothing(z)$. And $H_n\xrightarrow{\mathcal{H}}H$ does not imply $\Sigma_{H_n}\xrightarrow{\rm Cara}\Sigma_H$. For example, let $H_n=\{z\in\mathbb{H}: |\operatorname{Re}z|\le1,\operatorname{Im}z\le1/n\}$ for $n\in\mathbb{N}$. Then $H_n\xrightarrow{\mathcal{H}}\varnothing$, but $\Sigma_{H_n}\xrightarrow{\rm Cara}\mathbb{C}\setminus[-1,1]\ne\mathbb{C}=\Sigma_\varnothing$.

Suppose $H_n\xrightarrow{\mathcal{H}}H$, $K_n\xrightarrow{\mathcal{H}}K$ and $K_n\subset H_n$ for each $n$. Then $\mathbb{H}\setminus H_n\xrightarrow{\rm Cara}\mathbb{H}\setminus H$, $\mathbb{H}\setminus K_n\xrightarrow{\rm Cara}\mathbb{H}\setminus K$ and $\mathbb{H}\setminus H_n\subset\mathbb{H}\setminus K_n$ for each $n$. Since $(\mathbb{H}\setminus H)\cap(\mathbb{H}\setminus K)=\mathbb{H}\setminus(H\cup K)\ne\varnothing$, so $\mathbb{H}\setminus H\subset\mathbb{H}\setminus K$. Thus $K\subset H$. Let $L_n=H_n/K_n$ for each $n$ and $L=K/H$. Then $\varphi_{L_n}^{-1}=\varphi_{K_n}\circ\varphi_{H_n}^{-1}$ and $\varphi_L^{-1}=\varphi_K\circ\varphi_H^{-1}$. Since $\varphi_{H_n}^{-1}\xrightarrow{\rm l.u.}\varphi_H^{-1}$ in $\mathbb{H}$, and $\varphi_{K_n}\xrightarrow{\rm l.u.}\varphi_K$ in $\mathbb{H}\setminus K\supset\mathbb{H}\setminus H=\varphi_H(\mathbb{H})$, so $\varphi_{L_n}^{-1}\xrightarrow{\rm l.u.}\varphi_L^{-1}$ in $\mathbb{H}$. Thus $L_n\xrightarrow{\mathcal{H}}L$, that is, $H_n/K_n\xrightarrow{\mathcal{H}}H/K$.

LEMMA 5.4 (Compactness). *$\mathcal{H}(H)$ and $\mathcal{H}_{\rm sq}(H)$ are compact. Moreover, we have:*

(i) *Suppose $(K_n)$ is a sequence in $\mathcal{H}(H)$, then it has a subsequence $(L_n)$ that converges to some $K\in\mathcal{H}(H)$ w.r.t. $d_{\mathcal{H}}$, and $\varphi_{L_n}^{-1}\xrightarrow{\rm l.u.}\varphi_K^{-1}$ in $\mathbb{C}\setminus[c_H,d_H]$, $\Sigma_{L_n}\setminus[a_H,b_H]\xrightarrow{\rm Cara}\Sigma_K\setminus[a_H,b_H]$, and $\varphi_{L_n}\xrightarrow{\rm l.u.}\varphi_K$ in $\Sigma_K\setminus[a_H,b_H]$.*

(ii) *Suppose $(K_n)$ is a sequence in $\mathcal{H}_{sq}(H)$; then it has a subsequence $(L_n)$ that converges to some $K\in\mathcal{H}_{sq}(H)$ w.r.t. $d_{\mathcal{H}}$, and $\varphi_{L_n}^{-1}\xrightarrow{\rm l.u.}\varphi_K^{-1}$ in $\mathbb{C}\setminus[c_H,d_H]$, $\Sigma_{L_n}\setminus[c_H,d_H]\xrightarrow{\rm Cara}\Sigma_K\setminus[c_H,d_H]$, and $\varphi_{L_n}\xrightarrow{\rm l.u.}\varphi_K$ in $\Sigma_K\setminus[c_H,d_H]$.*

PROOF. (i) Since $\{\varphi_{K_n}^{-1}(z)-z:n\in\mathbb{N}\}$ is uniformly bounded in $\mathbb{C}\setminus[c_H,d_H]$, so $(K_n)$ has a subsequence $(L_n)$ such that $\varphi_{L_n}^{-1}(z)-z$ converges to some function $f$ locally uniformly in $\mathbb{C}\setminus[c_H,d_H]$. Then $|f(z)|\le M$ for any $z\in\mathbb{C}\setminus[c_H,d_H]$. Let $g(z)=f(z)+z$ for $z\in\mathbb{C}\setminus[c_H,d_H]$. Then $\varphi_{L_n}^{-1}\xrightarrow{\rm l.u.}g$ in $\mathbb{C}\setminus[c_H,d_H]$. There are $z_1,z_2\in\mathbb{C}\setminus[c_H,d_H]$ with $|z_1-z_2|>2M$. Then $|g(z_1)-g(z_2)|\ge|z_1-z_2|-|g(z_1)-z_1|-|g(z_2)-z_2|>2M-M-M=0$. So $g$ is not constant. From Lemma 5.1, $g$ is a conformal map. Since for each $n$, $\mathbb{H}\supset\varphi_{L_n}^{-1}(\mathbb{H})=\mathbb{H}\setminus L_n\supset\mathbb{H}\setminus H$, so $\mathbb{H}\supset g(\mathbb{H})\supset\mathbb{H}\setminus H$. Let $K=\mathbb{H}\setminus g(\mathbb{H})$. Then $K\in\mathcal{H}(H)$, and $g$ maps $\mathbb{H}$ conformally onto $\mathbb{H}\setminus K$. Since $\varphi_{L_n}^{-1}(z)-z=O(1/z)$ as $z\to\infty$, so $g(z)-z=O(1/z)$ as $z\to\infty$. Thus $g(z)=\varphi_K^{-1}(z)$ for $z\in\mathbb{C}\setminus[c_H,d_H]$. So $\varphi_{L_n}^{-1}\xrightarrow{\rm l.u.}\varphi_K^{-1}$ in $\mathbb{C}\setminus[c_H,d_H]$. Especially, $\varphi_{L_n}^{-1}\xrightarrow{\rm l.u.}\varphi_K^{-1}$ in $\mathbb{H}$. So $K$ is a subsequential limit of $(K_n)$. Thus $\mathcal{H}(H)$ is compact.

For $L\in\mathcal{H}(H)$, let $\Sigma_L^1:=\Sigma_L\setminus[a_H,b_H]$, $\Sigma_L^2:=\Sigma_L\setminus[c_H,d_H]$. Then $\Sigma_L^2\subset\Sigma_L^1$, and

$$\Sigma_L^1=(\mathbb{H}\setminus L)\cup\{z\in\mathbb{C}:\overline{z}\in\mathbb{H}\setminus L\}\cup(-\infty,a_H)\cup(b_H,+\infty),\tag{5.2}$$

$$\Sigma_L^2=(\mathbb{H}\setminus L)\cup\{z\in\mathbb{C}:\overline{z}\in\mathbb{H}\setminus L\}\cup(-\infty,c_H)\cup(d_H,+\infty),\tag{5.3}$$

because $(\mathbb{C}\setminus\Sigma_L)\cap\mathbb{R}\subset[a_L,b_L]\subset[a_H,b_H]\subset[c_H,d_H]$. So from $\mathbb{H}\setminus L_n\xrightarrow{\rm Cara}\mathbb{H}\setminus K$, we have $\Sigma_{L_n}^j\xrightarrow{\rm Cara}\Sigma_K^j$ for $j=1,2$. From Lemma 5.1, $\varphi_{L_n}^{-1}(\mathbb{C}\setminus[c_H,d_H])\xrightarrow{\rm Cara}\varphi_K^{-1}(\mathbb{C}\setminus[c_H,d_H])$ and $\varphi_{L_n}\xrightarrow{\rm l.u.}\varphi_K$ in $\varphi_K^{-1}(\mathbb{C}\setminus[c_H,d_H])$. Note that $\varphi_K^{-1}(\mathbb{C}\setminus[c_H,d_H])\supset\Sigma_K^2$, where the inclusion follows from Lemma 5.2. Thus $\varphi_{L_n}\xrightarrow{\rm l.u.}\varphi_K$ in $\Sigma_K^2$.

Since $|\varphi_{L_n}(z)-z|\le M$ for all $n\in\mathbb{N}$ and $z\in\Sigma_{L_n}$, and $\Sigma_{L_n}^1\subset\Sigma_{L_n}$, so every subsequence of $(\varphi_{L_n})$ has a subsequence that converges to some analytic function $h$ locally uniformly in $\Sigma_K^1$. Since $\varphi_{L_n}\xrightarrow{\rm l.u.}\varphi_K$ in $\Sigma_K^2\subset\Sigma_K^1$, so $h$ agrees with $\varphi_K$ on $\Sigma_K^2$. Since they are both analytic, so $h$ agrees with $\varphi_K$ on $\Sigma_K^1$. Since all subsequential limits of $\varphi_{L_n}$ in $\Sigma_K^1$ are the same function $\varphi_K$, so $\varphi_{L_n}\xrightarrow{\rm l.u.}\varphi_K$ in $\Sigma_K^1=\Sigma_K\setminus[a_H,b_H]$.

(ii) Suppose $K_n = K_n^2/K_n^1$ with $K_n^1 \subset K_n^2 \subset H$. From (i), $(K_n)$ has a subsequence $(L_n = L_n^2/L_n^1)$ such that $L_n^j \xrightarrow{\mathcal{H}} K^j$ for some $K^j \in \mathcal{H}(H)$, $j = 1, 2$. Since $L_n^1 \subset L_n^2$ for each $n$, so $K^1 \subset K^2$. Let $K = K^2/K^1$. Then $K \in \mathcal{H}_{sq}(H)$, and $L_n = L_n^2/L_n^1 \xrightarrow{\mathcal{H}} K^2/K^1 = K$. So $K$ is a subsequential limit of $(K_n)$. Thus $\mathcal{H}_{sq}(H)$ is compact.

Since $\{\varphi_{L_n}^{-1}(z) - z : n \in \mathbb{N}\}$ is uniformly bounded in $\mathbb{C} \setminus [c_H, d_H]$, so every subsequence of $(\varphi_{L_n}^{-1})$ has a subsequence which converges to some $h$ locally uniformly in $\mathbb{C} \setminus [c_H, d_H]$. Then $h$ agrees with $\varphi_K^{-1}$ on $\mathbb{H}$. Since they are both analytic in $\mathbb{C} \setminus [c_H, d_H]$, so $h$ agrees with $\varphi_K^{-1}$ on $\mathbb{C} \setminus [c_H, d_H]$. Thus $(\varphi_{L_n}^{-1}) \xrightarrow{\text{l.u.}} \varphi_K^{-1}$ in $\mathbb{C} \setminus [c_H, d_H]$.

For $L \in \mathcal{H}_{\text{sq}}(H)$, we define $\Sigma_L^j$, $j = 1, 2$, as in (i). Then (5.3) still holds because $[a_L, b_L] \subset [c_L, d_L] \subset [c_H, d_H]$, but (5.2) does not because $[a_L, b_L] \subset [a_H, b_H]$ may not be true. A similar argument gives that $\varphi_{L_n} \xrightarrow{\text{l.u.}} \varphi_K$ in $\Sigma_K^2 = \Sigma_K \setminus [c_H, d_H]$. □

5.3. *Lipschitz conditions.* Suppose $\xi \in C([0, a])$ for some $a > 0$, and $K_a^\xi \in \mathcal{H}(\alpha)$. Then for each $t \in [0, a]$, $\varphi_t^\xi = \varphi_{K_t^\xi}$, and $K_t^\xi \in \mathcal{H}(\alpha)$. For $0 \le t_1 < t_2 \le a$, let $K_{t_1,t_2}^\xi = K_{t_2}^\xi / K_{t_1}^\xi$. Then $K_{t_1,t_2}^\xi \in \mathcal{H}_{\text{sq}}(\alpha)$ and $\varphi_{K_{t_1,t_2}^\xi} = \varphi_{t_2}^\xi \circ (\varphi_{t_1}^\xi)^{-1}$, $\varphi_{K_{t_1,t_2}^\xi}^{-1} = \varphi_{t_1}^\xi \circ (\varphi_{t_2}^\xi)^{-1}$. Since $\xi(t_1) \in \overline{K_{t_1,t_2}^\xi}$, so

$$\xi(t_1) \in [a_{K_{t_1,t_2}^\xi}, b_{K_{t_1,t_2}^\xi}] \subset [c_{K_{t_1,t_2}^\xi}, d_{K_{t_1,t_2}^\xi}] \subset [c_\alpha, d_\alpha].$$

This holds for any $t_1 \in [0, a)$. Since $\xi$ is continuous, so we also have $\xi(a) \in [c_\alpha, d_\alpha]$.

Lemma 5.5. *Suppose $\alpha_0$ and $\alpha_1$ are crosscuts in $\mathbb{H}$, and $\alpha_0$ is strictly enclosed by $\alpha_1$. Then there are $\delta, C > 0$ such that if $\zeta, \eta \in C([0, a])$, $\|\zeta - \eta\|_a < \delta$, and $K_a^\zeta \subset H(\alpha_0)$, then $K_a^\eta \subset H(\alpha_1)$, and for any $z \in \overline{\mathbb{H} \setminus H(\alpha_1)}$,*

$$|\varphi_a^\zeta(z) - \varphi_a^\eta(z)| \le Ca\|\zeta - \eta\|_a.$$

Proof. Suppose $\zeta, \eta \in C([0, a])$ and $K_a^\zeta \subset H(\alpha_0)$. Choose a crosscut $\alpha_{0.5}$ in $\mathbb{H}$ that strictly encloses $\alpha_0$, and is strictly enclosed by $\alpha_1$. Then $\overline{\alpha_{0.5}}$ and $\overline{\alpha_1}$ are disjoint compact subsets of $\Sigma_{\alpha_0}$, which contains $\Sigma_K \setminus [a_{\alpha_0}, b_{\alpha_0}]$ for any $K \in \mathcal{H}(\alpha_0)$. From the compactness of $\mathcal{H}(\alpha_0)$, there is $d > 0$, such that the distance between $\varphi_K(\alpha_{0.5})$ and $\varphi_K(\alpha_1)$ is at least $d$ for any $K \in \mathcal{H}(\alpha_0)$. For $t \in [0, a]$, since $K_t^\zeta \in \mathcal{H}(\alpha_0)$, so the distance between $\varphi_t^\zeta(\alpha_{0.5})$ and $\varphi_t^\zeta(\alpha_1)$ is at least $d$. Since $K_a^\zeta$ is enclosed by $\alpha_{0.5}$, so $K_{t,a}^\zeta = \varphi_t^\zeta(K_a^\zeta \setminus K_t^\zeta)$ is enclosed by $\varphi_t^\zeta(\alpha_{0.5})$, which implies that $\zeta(t) \in \overline{K_{t,a}^\zeta}$ is enclosed by $\varphi_t^\zeta(\alpha_{0.5})$. Thus the distance between $\zeta(t)$ and $\varphi_t^\zeta(z)$ is at least $d$ for any $z \in \overline{\mathbb{H}\setminus H(\alpha_1)}$ and $t \in [0,a]$. Fix $z \in \overline{\mathbb{H}\setminus H(\alpha_1)}$ and $\delta \in (0, d/3]$. Then $|\varphi_t^\zeta(z) - \zeta(t)| \ge d$ for any $t \in [0,a]$. Suppose $\|\zeta - \eta\|_a < \delta$. Note that $\varphi_0^\zeta(z) = z = \varphi_0^\eta(z)$. Let $[0,b)$ be the maximal subinterval of $[0,a)$ on which $\varphi_t^\eta(z)$ is defined and $|\varphi_t^\zeta(z) - \varphi_t^\eta(z)| \le d/3$. Then for any $t \in [0,b)$,

$$|\varphi_t^\eta(z) - \eta(t)| \ge |\varphi_t^\zeta(z) - \zeta(t)| - |\varphi_t^\zeta(z) - \varphi_t^\eta(z)| - |\zeta(t) - \eta(t)| \ge d/3.$$

Thus $\varphi_b^\eta(z)$ is also defined. From the chordal Loewner equation, for $t \in [0,b]$,

$$\begin{aligned}
|\varphi_t^\zeta(z) - \varphi_t^\eta(z)| &\le \int_0^t \left| \frac{2}{\varphi_s^\zeta(z) - \zeta(s)} - \frac{2}{\varphi_s^\eta(z) - \eta(s)} \right| ds \\
&\le \int_0^t \left| \frac{2(\zeta(s) - \eta(s))}{(\varphi_s^\zeta(z) - \zeta(s))(\varphi_s^\eta(z) - \eta(s))} \right| ds \\
&\quad + \int_0^t \left| \frac{2(\varphi_s^\eta(z) - \varphi_s^\zeta(z))}{(\varphi_s^\zeta(z) - \zeta(s))(\varphi_s^\eta(z) - \eta(s))} \right| ds \\
&\le \frac{6t}{d^2}\|\zeta - \eta\|_a + \frac{6}{d^2}\int_0^t |\varphi_s^\zeta(z) - \varphi_s^\eta(z)|\, ds \qquad (5.4)\\
&\le \frac{6\delta t}{d^2} + \frac{6}{d^2}\int_0^t |\varphi_s^\zeta(z) - \varphi_s^\eta(z)|\, ds. \qquad (5.5)
\end{aligned}$$

Solving inequality (5.5), we get

$$|\varphi_b^\zeta(z) - \varphi_b^\eta(z)| \le \delta(e^{3b/d^2} - 1) \le \delta(e^{3a/d^2} - 1).$$

Let $h = \mathrm{hcap}(H(\alpha_0))$. Then $a = \mathrm{hcap}(K_a^\zeta)/2 \le h/2$. Choose $\delta = \min\{d/3, \frac{d/6}{e^{3h/d^2}-1}\}$. Then $|\varphi_b^\zeta(z) - \varphi_b^\eta(z)| \le d/6$. So we have $b = a$, which implies that $\varphi_t^\eta(z)$ is defined on $[0,a]$, that is, $z \notin K_a^\eta$. Since this is true for any $z \in \mathbb{H}\setminus H(\alpha_1)$, so $K_a^\eta \subset H(\alpha_1)$. Finally, let $C = (\exp(\frac{3h}{d^2}) - 1)/(h/2)$. Solving inequality (5.4) for $t \in [0,a]$, we get

$$|\varphi_a^\zeta(z) - \varphi_a^\eta(z)| \le (e^{6a/d^2} - 1)\|\zeta - \eta\|_a \le Ca\|\zeta - \eta\|_a$$

for any $z \in \overline{\mathbb{H}\setminus H(\alpha_1)}$, where the second "$\le$" holds because $a \le h/2$. □

Lemma 5.6. *Suppose $\alpha$ and $\rho$ are crosscuts in $\mathbb{H}$, and $[c_\alpha, d_\alpha]$ is strictly enclosed by $\rho$. Then there are $\delta, C > 0$ such that if $\zeta, \eta \in C([0,a])$, $\|\zeta - \eta\|_a < \delta$, and $K_a^\zeta \subset H(\alpha)$, then $K_a^\eta$ is enclosed by $(\varphi_a^\zeta)^{-1}(\rho)$, and for any $w \in \overline{\mathbb{H}\setminus H(\rho)}$,*

$$|w - \varphi_a^\eta \circ (\varphi_a^\zeta)^{-1}(w)| \le Ca\|\zeta - \eta\|_a.$$

Proof. Suppose $\zeta,\eta\in C([0,a])$ and $K_a^\zeta\subset H(\alpha)$. Choose $\rho_0$ that strictly encloses $[c_\alpha,d_\alpha]$, and is strictly enclosed by $\rho$. Then for any $t\in[0,a)$, $\zeta(t)\in\overline{K_{t,a}^\zeta}$ is enclosed by $\varphi^{-1}_{K_{t,a}^\zeta}(\rho_0)$. Note that $K_{t,a}^\zeta\in\mathcal{H}_{sq}(\alpha)$ and $\varphi^{-1}_{K_{t,a}^\zeta}=\varphi_t^\zeta\circ(\varphi_a^\zeta)^{-1}$. From the compactness of $\mathcal{H}_{\rm sq}(\alpha)$ and an argument that is similar to the first paragraph of the last proof, we see that there is $d>0$ depending only on $\alpha$ and $\rho$ such that $|\varphi_t^\zeta\circ(\varphi_a^\zeta)^{-1}(w)-\zeta(t)|\ge d$ for any $t\in[0,a]$ and $w\in\overline{\mathbb{H}\setminus H(\rho)}$. Fix $w\in\overline{\mathbb{H}\setminus H(\rho)}$. Applying the argument of the proof of the last lemma to $z=(\varphi_a^\zeta)^{-1}(w)$, we have $\delta,C>0$ depending only on $\alpha$ and $\rho$ such that if $\|\zeta-\eta\|_a<\delta$, then $\varphi_a^\eta(z)$ is well defined, and

$$|w-\varphi_a^\eta\circ(\varphi_a^\zeta)^{-1}(w)|=|\varphi_a^\zeta(z)-\varphi_a^\eta(z)|\le Ca\|\zeta-\eta\|_a.$$

That $\varphi_a^\eta(z)$ is well defined implies that $(\varphi_a^\zeta)^{-1}(w)=z\notin K_a^\eta$. Since this holds for any $w\in\mathbb{H}\setminus H(\rho)$, so $K_a^\eta$ is enclosed by $(\varphi_a^\zeta)^{-1}(\rho)$. □

Now suppose $\Omega$ is an almost $\mathbb{H}$ domain, and $p\in\Omega$. Suppose $\alpha$ is a crosscut in $\mathbb{H}$ such that $H(\alpha)\subset\Omega\setminus\{p\}$. From the compactness of $\mathcal{H}(\alpha)$, there is $\mathbf{h}>0$ depending only on $\Omega,p,\alpha$, such that if $K\in\mathcal{H}(\alpha)$, then $\operatorname{dist}(\varphi_K(\{p\}\cup\partial\Omega\setminus\mathbb{R}),\mathbb{R})\ge\mathbf{h}$. Let $\rho_1$ and $\rho_2$ be crosscuts in $\mathbb{H}$ with height smaller than $\mathbf{h}/2$ such that $\rho_1$ strictly encloses $\rho_2$, and $\rho_2$ strictly encloses $[c_\alpha,d_\alpha]$. Then for any $K\in\mathcal{H}(\alpha)$, $H(\varphi_K^{-1}(\rho))\subset\Omega_K\setminus\{\varphi_K(p)\}$.

Lemma 5.7. *There are $\delta,C>0$ such that if $\zeta,\eta\in C([0,a])$, $\|\zeta-\eta\|_a<\delta$, and $K_a^\zeta\subset H(\alpha)$, then for any $z\in\rho_1$,*

$$|J_a^\zeta(z)-J_a^\eta(z)|\le Ca\|\zeta-\eta\|_a. \tag{5.6}$$

Proof. Choose a crosscut $\alpha_1$ in $\mathbb{H}$ that strictly encloses $\alpha$ such that $H(\alpha_1)\subset\Omega\setminus\{p\}$. Suppose $\zeta,\eta\in C([0,a])$ and $K_a^\zeta\subset H(\alpha)$. From Lemma 5.5, there is $\delta_0>0$ depending only on $\alpha$ and $\alpha_1$ such that if $\|\zeta-\eta\|_a<\delta_0$, then $K_a^\eta\subset H(\alpha_1)$.

From Lemma 5.6, there are $\delta_1,C_1>0$ depending only on $\alpha,\rho_1,\rho_2$, such that if $\|\zeta-\eta\|_a<\delta_1$, then $K_a^\eta$ is enclosed by $(\varphi_a^\zeta)^{-1}(\rho_2)$, and for any $z\in\rho_1\cup\rho_2$,

$$|z-\varphi_a^\eta\circ(\varphi_a^\zeta)^{-1}(z)|\le C_1a\|\zeta-\eta\|_a. \tag{5.7}$$

Let $F=\{z\in\overline{\mathbb{H}}:\operatorname{dist}(z,H(\rho_1))\le\mathbf{h}/4\}$. From the compactness of $\mathcal{H}(\alpha_1)$, there is $D>0$ depending only on $\Omega,p,\alpha_1,F$, such that if $K_t^\xi\in\mathcal{H}(\alpha_1)$, then for any $z\in F$,

$$|\nabla J_t^\xi(z)|\le D. \tag{5.8}$$

Let $h_0=\operatorname{hcap}(H(\alpha))$. Then $a=\operatorname{hcap}(K_a^\zeta)/2\le h_0/2$. Let $\delta=\min\{\delta_0,\delta_1,\mathbf{h}/(2C_1h_0)\}$. Suppose $\|\zeta-\eta\|_a<\delta$. Then for any $z\in\rho_1\cup\rho_2$,

$$|z-\varphi_a^\eta\circ(\varphi_a^\zeta)^{-1}(z)|\le C_1a\delta\le C_1h_0\delta/2\le\mathbf{h}/4,$$

which implies that $[z,\varphi^\eta_a\circ(\varphi^\zeta_a)^{-1}(z)]\subset F$.

Define $G^\xi_t=G(\Omega\setminus K^\xi_t,p;\cdot)$ if $K^\xi_t\subset\Omega\setminus\{p\}$. For $j=1,2$, let

$$N_j=\sup_{z\in\rho_j}\{|J^\zeta_a(z)-J^\eta_a(z)|\}=\sup_{z\in\rho_j}\{|G^\zeta_a\circ(\varphi^\zeta_a)^{-1}(z)-G^\eta_a\circ(\varphi^\eta_a)^{-1}(z)|\};$$

$$N'_j=\sup_{z\in(\varphi^\zeta_a)^{-1}(\rho_j)}\{|G^\zeta_a(z)-G^\eta_a(z)|\}=\sup_{z\in\rho_j}\{|G^\zeta_a\circ(\varphi^\zeta_a)^{-1}(z)-G^\eta_a\circ(\varphi^\zeta_a)^{-1}(z)|\}.$$

There is $q\in(0,1)$ depending only on $\rho_1$ and $\rho_2$ such that for any $z\in\rho_2$, the probability that a plane Brownian motion started from $z$ hits $\rho_1$ before $\mathbb{R}$ is less than $q$. Since both $J^\zeta_a$ and $J^\eta_a$ are harmonic in $H(\rho_1)$, have continuations to $\overline{H(\rho_1)}$, and vanish on $\mathbb{R}$, so $J^\zeta_a-J^\eta_a$ also has these properties. Since $\rho_2\subset H(\rho_1)$, so

$$N_2\le qN_1. \tag{5.9}$$

Since $K^\zeta_a$ and $K^\eta_a$ are enclosed by $(\varphi^\zeta_a)^{-1}(\rho_2)$, so $G^\zeta_a$ and $G^\eta_a$ are harmonic in $\Omega\setminus\{p\}\setminus H((\varphi^\zeta_a)^{-1}(\rho_2))$. Since they both behave like $-\ln(z-p)/(2\pi)+O(1)$ near $p$, so $G^\zeta_a-G^\eta_a$ has a harmonic extension in $\Omega\setminus H((\varphi^\zeta_a)^{-1}(\rho_2))$. Since $G^\zeta_a-G^\eta_a$ vanishes at every boundary point of $\Omega\setminus H((\varphi^\zeta_a)^{-1}(\rho_2))$ including $\infty$, except on $(\varphi^\zeta_a)^{-1}(\rho_2)$, and $(\varphi^\zeta_a)^{-1}(\rho_1)\subset\Omega\setminus H((\varphi^\zeta_a)^{-1}(\rho_2))$, so from the maximum principle for harmonic functions,

$$N'_1\le N'_2. \tag{5.10}$$

Fix $j\in\{1,2\}$. From $[z,\varphi^\eta_a\circ(\varphi^\zeta_a)^{-1}(z)]\subset F$ for $z\in\rho_j$, $K^\eta_a\in\mathcal{H}(\alpha_1)$, and (5.7) and (5.8), we have

$$\begin{aligned}|N_j-N'_j|&\le\sup_{z\in\rho_j}\{|G^\eta_a\circ(\varphi^\eta_a)^{-1}(z)-G^\eta_a\circ(\varphi^\zeta_a)^{-1}(z)|\}\\&=\sup_{z\in\rho_j}\{|J^\eta_a(z)-J^\eta_a(\varphi^\eta_a\circ(\varphi^\zeta_a)^{-1}(z))|\}\\&\le\sup_{w\in F}|\nabla J^\eta_a(w)|\sup_{z\in\rho_j}\{|z-\varphi^\eta_a\circ(\varphi^\zeta_a)^{-1}(z)|\}\\&\le DC_1a\|\zeta-\eta\|_a.\end{aligned}$$

From (5.9), (5.10) and the above inequality, we have

$$\begin{aligned}N_1&\le N'_1+DC_1a\|\zeta-\eta\|_a\le N'_2+DC_1a\|\zeta-\eta\|_a\\&\le N_2+2DC_1a\|\zeta-\eta\|_a\le qN_1+2DC_1a\|\zeta-\eta\|_a,\end{aligned}$$

which implies that $N_1\le Ca\|\zeta-\eta\|_a$, where $C=2DC_1/(1-q)$. So we get (5.6). $\square$

Lemma 5.8. *There are $\delta,C>0$ such that if $\zeta,\eta\in C([0,a])$, $\|\zeta-\eta\|_a<\delta$, and $K^\zeta_a\subset H(\alpha)$, then $|X^\zeta_a-X^\eta_a|\le C\|\zeta-\eta\|_a$.*

PROOF. Suppose $\zeta,\eta\in C([0,a])$ and $K_a^\zeta\subset H(\alpha)$. Choose a crosscut $\alpha_1$ in $\mathbb{H}$ that strictly encloses $\alpha$ such that $H(\alpha_1)\subset\Omega\setminus\{p\}$. Let $\rho$ be a crosscut in $\mathbb{H}$ with height smaller than $\mathbf{h}/2$ that strictly encloses $[c_\alpha,d_\alpha]$. From Lemmas 5.5 and 5.7, there are $\delta_0,C_0>0$ depending only on $\Omega,p,\alpha,\alpha_1,\rho$, such that if $\|\zeta-\eta\|_a<\delta_0$, then $K_a^\eta\subset H(\alpha_1)$ and for any $z\in\rho$, $|J_a^\zeta(z)-J_a^\eta(z)|\le C_0a\|\zeta-\eta\|_a$. Let $d_0=\text{dist}([c_\alpha,d_\alpha],\rho)/2>0$, and $\delta=\delta_0\wedge d_0$.

Suppose $\|\zeta-\eta\|_a<\delta$. Then $K_a^\zeta,K_a^\eta\subset H(\alpha_1)$. From the compactness of $\mathcal{H}(\alpha_1)$, there are $m,M_1,M_2,M_3>0$ depending only on $\Omega,p,\alpha,\alpha_1,\rho$, such that for any $x\in[c_\alpha-d_0,d_\alpha+d_0]$,

$$m\le\partial_yJ_a^\zeta(x),\partial_yJ_a^\eta(x)\le M_1\quad\text{and}\quad|\partial_x^{j-1}\partial_yJ_a^\zeta(x)|,|\partial_x^{j-1}\partial_yJ_a^\eta(x)|\le M_j,$$

for $j=2,3$. Let $C_1=M_3/m+M_2^2/m^2$. So for any $x\in[c_\alpha-d_0,d_\alpha+d_0]$,

$$\begin{aligned}&|\partial_x(\partial_x\partial_y/\partial_y)J_a^\zeta(x)|\\&\quad=|(\partial_x^2\partial_y/\partial_y-((\partial_x\partial_y\cdot\partial_x\partial_y)/(\partial_y\cdot\partial_y)))J_a^\zeta(x)|\le C_1.\end{aligned}\tag{5.11}$$

Since $\text{dist}([c_\alpha-d_0,d_\alpha+d_0],\rho)\ge d_0$, so for any $x\in[c_\alpha-d_0,d_\alpha+d_0]$,

$$|\partial_x^{j-1}\partial_y(J_a^\zeta-J_a^\eta)(x)|\le\frac{2j\,!}{d_0^j}\sup_{z\in\rho}|J_a^\zeta(z)-J_a^\eta(z)|\le\frac{2j\,!}{d_0^j}C_0a\|\zeta-\eta\|_a,$$

for $j=1,2$, from which follows that

$$\begin{aligned}&|(\partial_x\partial_y/\partial_y)J_a^\zeta(x)-(\partial_x\partial_y/\partial_y)J_a^\eta(x)|\\&\quad=|\partial_x\partial_yJ_a^\zeta(x)\,\partial_yJ_a^\eta(x)-\partial_x\partial_yJ_a^\eta(x)\,\partial_yJ_a^\zeta(x)|/|\partial_yJ_a^\zeta(x)\partial_yJ_a^\eta(x)|\\&\quad\le|\partial_x\partial_yJ_a^\zeta(x)\,\partial_yJ_a^\eta(x)-\partial_x\partial_yJ_a^\zeta(x)\,\partial_yJ_a^\zeta(x)|/m^2\\&\qquad+|\partial_x\partial_yJ_a^\zeta(x)\,\partial_yJ_a^\zeta(x)-\partial_x\partial_yJ_a^\eta(x)\,\partial_yJ_a^\zeta(x)|/m^2\\&\quad\le M_2|\partial_y(J_a^\zeta-J_a^\eta)(x)|/m^2+M_1|\partial_x\partial_y(J_a^\zeta-J_a^\eta)(x)|/m^2\\&\quad\le(2M_2/d_0+4M_1/d_0^2)C_0a\|\zeta-\eta\|_a/m^2\le C_2\|\zeta-\eta\|_a,\end{aligned}\tag{5.12}$$

if we let $C_2=(M_2/d_0+2M_1/d_0^2)C_0\,\text{hcap}(H(\alpha))/m^2$.

Since $K_a^\zeta\in\mathcal{H}(\alpha)$, so $\zeta(a)\in[c_\alpha,d_\alpha]$. From $|\eta(a)-\zeta(a)|\le\delta\le d_0$, we have $\eta(a)\in[c_\alpha-d_0,d_\alpha+d_0]$. Let $C=C_1+C_2$. From (5.11) and (5.12), we have

$$\begin{aligned}|X_a^\zeta-X_a^\eta|&=|(\partial_x\partial_y/\partial_y)J_a^\zeta(\zeta(a))-(\partial_x\partial_y/\partial_y)J_a^\eta(\eta(a))|\\&\le|(\partial_x\partial_y/\partial_y)J_a^\zeta(\zeta(a))-(\partial_x\partial_y/\partial_y)J_a^\zeta(\eta(a))|\\&\quad+|(\partial_x\partial_y/\partial_y)J_a^\zeta(\eta(a))-(\partial_x\partial_y/\partial_y)J_a^\eta(\eta(a))|\\&\le C_1|\zeta(a)-\eta(a)|+C_2\|\zeta-\eta\|_a\le C\|\zeta-\eta\|_a.\end{aligned}$$

□

PROOF OF THEOREM 3.1. Let $\xi_0(t)=A(0)$, $t\in[0,\infty)$. We may have $a_0>0$ such that $K_{a_0}^{\xi_0}\subset\Omega\setminus\{p\}$. Choose crosscuts $\alpha_0$ and $\alpha_1$ in $\mathbb{H}$ such that $K^{\xi_0}_{a_0}$ is enclosed by $\alpha_0$, $\alpha_0$ is strictly enclosed by $\alpha_1$, and $H(\alpha_1)\subset\Omega\setminus\{p\}$. From Lemma 5.5, there is $\delta_1>0$ such that for any $t\in[0,a_0]$, if $\eta\in C([0,t])$ satisfies $\|\eta-\xi_0\|_t<\delta_1$, then $K^\eta_t\subset H(\alpha_1)$. Let $\delta_2, C>0$ be the constants given by Lemma 5.8 with $\alpha=\alpha_1$. Let $\delta=\delta_1\wedge(\delta_2/2)$. Then for any $t\in(0,a_0]$, if $\eta_j\in C([0,t])$ and $\|\eta_j-\xi_0\|_t<\delta$, $j=1,2$, then

$$|X^{\eta_1}_t-X^{\eta_2}_t|\le C\|\eta_1-\eta_2\|_t. \tag{5.13}$$

Define a sequence of functions $(\xi_n(t))$ by induction:

$$\xi_{n+1}(t)=A(t)+\lambda\int_0^t X^{\xi_n}_s\,ds, \tag{5.14}$$

as long as $X^{\xi_n}_s$, $0\le s\le t$, are defined. From Lemma 6.3, we see that $X^\xi_t$ is continuous, and so the integral makes sense. We may choose $a\in(0,a_0]$ such that $|\lambda|Ca<1/2$ and $\|\xi_1-\xi_0\|_a<\delta/2$. For $n=1$, we have $\|\xi_n-\xi_0\|_a<(1-1/2^n)\delta$ and $\|\xi_n-\xi_{n-1}\|_a<\delta/2^n$. Suppose this is true for some $n\in\mathbb{N}$. Then from (5.13) and (5.14), we have

$$\begin{aligned}|\xi_{n+1}(t)-\xi_n(t)|&\le|\lambda|\int_0^t|X^{\xi_n}_s-X^{\xi_{n-1}}_s|\,ds\le|\lambda|\int_0^t C\|\xi_n-\xi_{n-1}\|_a\,ds\\&\le|\lambda|Ca\|\xi_n-\xi_{n-1}\|_a<\|\xi_n-\xi_{n-1}\|_a/2<\delta/2^{n+1},\end{aligned}$$

for $t\in[0,a]$. Thus $\|\xi_{n+1}-\xi_n\|_a<\delta/2^{n+1}$, and $\|\xi_{n+1}-\xi_0\|_b<\delta/2^{n+1}+\|\xi_n-\xi_0\|_b<(1-1/2^{n+1})\delta$. From induction, we have $\|\xi_{n+1}-\xi_n\|_a<\delta/2^{n+1}$ for any $n\in\mathbb{N}$. Thus $(\xi_n)$ restricted to $[0,a]$ is a Cauchy sequence in $C([0,a])$. Let $\xi_\infty=\lim_{n\to\infty}\xi_n|_{[0,a]}\in C([0,a])$. Let $n\to\infty$ in (5.14); we see that $\xi_\infty$ solves (3.2) for $t\in[0,a]$.

Let $\mathcal{S}$ be the set of all couples $(\xi,T)$ such that $T>0$ and $\xi$ solves (3.2) for $t\in[0,T]$. We have proved that $\mathcal{S}$ is nonempty. We claim that if $(\xi,T)\in\mathcal{S}$, then there is $(\xi_e,T_e)\in\mathcal{S}$ such that $T_e>T$ and $\xi_e(t)=\xi(t)$ for $t\in[0,T]$. To prove this claim, let $\breve\Omega=\varphi^\xi_T(\Omega\setminus K^\xi_T)$ and $\breve p=\varphi^\xi_T(p)$. If $K^{\breve\xi}_t\subset\breve\Omega\setminus\{\breve p\}$, let $\breve J^{\breve\xi}_t=G(\breve\Omega\setminus K^{\breve\xi}_t,\breve p;\cdot)\circ(\varphi^{\breve\xi}_t)^{-1}$, and $\breve X^{\breve\xi}_t=(\partial_x\partial_y/\partial_y)\breve J^{\breve\xi}_t(\breve\xi(t))$. From the first part of the proof, the solution to

$$\breve\xi(t)=\xi(T)+A(T+t)-A(T)+\lambda\int_0^t\breve X^{\breve\xi}_s\,ds \tag{5.15}$$

exists on $[0,\breve a]$ for some $\breve a>0$. Let $T_e=T+\breve a>T$. Define $\xi_e(t)=\xi(t)$ for $t\in[0,T]$ and $\xi_e(t)=\breve\xi(t-T)$ for $t\in[T,T_e]$. It is clear that $\xi_e\in C([0,T_e])$. Since $\xi_e$ agrees with $\xi$ on $[0,T]$, so $\xi_e$ solves (3.2) for $t\in[0,T]$. For $t\in[0,T_e-T]$, we have $\varphi^{\xi_e}_{T+t}=\varphi^{\breve\xi}_t\circ\varphi^\xi_T$ and $K^{\xi_e}_{T+t}=K^\xi_T\cup(\varphi^\xi_T)^{-1}(K^{\breve\xi}_t)$. Since $\varphi^\xi_T$ maps $p$ to $\breve p$, and $\Omega\setminus K^{\xi_e}_{T+t}$ onto $\breve\Omega\setminus K^{\breve\xi}_t$, so

$$\begin{aligned}\breve J^{\breve\xi}_t&=G(\breve\Omega\setminus K^{\breve\xi}_t,\breve p;\cdot)\circ(\varphi^{\breve\xi}_t)^{-1}=G(\Omega\setminus K^{\xi_e}_{T+t},p;\cdot)\circ(\varphi^\xi_T)^{-1}\circ(\varphi^{\breve\xi}_t)^{-1}\\&=G(\Omega\setminus K^{\xi_e}_{T+t},p;\cdot)\circ(\varphi^{\xi_e}_{T+t})^{-1}=J^{\xi_e}_{T+t}.\end{aligned}$$

So $\breve{X}_t^{\breve{\xi}} = X_{T+t}^{\xi_e}$. Thus for $t \in [0, T_e - T]$,

$$\xi_e(T+t) = \breve{\xi}(t) = \xi(T) + A(T+t) - A(T) + \lambda \int_0^t \breve{X}_s^{\breve{\xi}}\, ds$$
$$= A(T+t) + \lambda \int_0^T X_s^{\xi}\, ds + \lambda \int_0^t X_{T+s}^{\xi_e}\, ds$$
$$= A(T+t) + \lambda \int_0^{T+t} X_s^{\xi_e}\, ds.$$

So $\xi_e$ solves (3.2) for $t \in [T, T_e]$. Thus $(\xi_e, T_e) \in \mathcal{S}$. So the claim is justified.

Suppose $(\xi_1, T_1), (\xi_2, T_2) \in \mathcal{S}$. For $j = 1, 2$, since $\xi_j(0) = A(0) = \xi_0(0)$, so there is $S_j \in (0, T_j \wedge a_0]$ such that $\|\xi_j - \xi_0\|_{S_j} < \delta$. Choose $S_3 \in (0, S_1 \wedge S_2]$ such that $C|\lambda|S_3 < 1$. From (3.2) and (5.13), we have $\|\xi_1 - \xi_2\|_{S_3} \le |\lambda| C S_3 \|\xi_1 - \xi_2\|_{S_3}$, so $\|\xi_1 - \xi_2\|_{S_3} = 0$, which means that $\xi_1(t) = \xi_2(t)$ for $0 \le t \le S_3$.

Let $T_0 = T_1 \wedge T_2$. We claim that $\xi_1(t) = \xi_2(t)$ for $t \in [0, T_0]$. Let $T \in [0, T_0]$ be the maximum such that $\xi_1(t) = \xi_2(t)$ for $t \in [0, T]$. Suppose $T < T_0$. Let $\breve{\xi}_1(t) = \xi_1(T+t)$, $\breve{\xi}_2(t) = \xi_2(T+t)$ for $t \in [0, T_0 - T]$. Then $\breve{\xi}_1$ and $\breve{\xi}_2$ both solve (5.15) for $t \in [0, T_0 - T]$. From the last paragraph, there is $S_3 \in (0, T_0 - T]$ such that $\breve{\xi}_1(t) = \breve{\xi}_2(t)$ for $0 \le t \le S_3$, which implies that $\xi_1(t) = \xi_2(t)$ for $0 \le t \le T + S_3$. This contradicts the maximum property of $T$. So $T = T_0$, and $\xi_1(t) = \xi_2(t)$ for $t \in [0, T_0]$.

Let $T_A = \sup\{T : (\xi, T) \in \mathcal{S}\}$. Define $\xi_A$ on $[0, T_A)$ as follows. For any $t \in [0, T_A)$, choose $(\xi, T) \in \mathcal{S}$ such that $t \le T$, and let $\xi_A(t) = \xi(t)$. From the last paragraph, $\xi_A$ is well defined, and solves (3.2) for $t \in [0, T_A)$. The uniqueness of $\xi_A$ also follows from the last paragraph. There is no solution to (3.2) defined on $[0, T_A]$. Otherwise, there exists some solution on $[0, T_A + \varepsilon]$ for some $\varepsilon > 0$, which contradicts the definition of $T_A$.

(i) Suppose $A_0 \in C([0,\infty))$, $a \in (0,\infty)$, and $T_{A_0} > a$. Then $K_a^{\xi_{A_0}} \subset \Omega \setminus \{p\}$. Choose a crosscut $\alpha$ in $\mathbb{H}$ such that $K_a^{\xi_{A_0}} \subset H(\alpha) \subset \Omega \setminus \{p\}$. Let $\delta_0, C_0 > 0$ be the $\delta, C$ given by Lemma 5.8 with $\zeta = \xi_{A_0}$. Let $C = \exp(C_0|\lambda|a)$ and $\delta = \delta_0/C$. Suppose $A \in C([0,\infty))$ and $\|A - A_0\|_a < \delta$. Then $|\xi_A(0) - \xi_{A_0}(0)| = |A(0) - A_0(0)| < \delta \le \delta_0$. Let $b \in [0, a \wedge T_A)$ be the maximal such that $|\xi_A(t) - \xi_{A_0}(t)| < \delta_0$ for $t \in [0, b)$. From the properties of $\delta_0$ and $C_0$, for $0 \le t < b$,

$$|X_t^{\xi_A} - X_t^{\xi_{A_0}}| \le C_0 \|\xi - \xi_{A_0}\|_t. \tag{5.16}$$

So $X_t^{\xi_A}$ is bounded on $[0, b)$. From (3.2), $\lim_{t \to b^-} \xi_A(t)$ exists. If $T_A = b$, we define $\xi_A(T) = \lim_{t \to b^-} \xi_A(t)$, then $\xi_A$ solves (3.2) for $t \in [0, T]$, which is a contradiction. Thus $T_A > b$. From (3.2) and (5.16), we have that for any $0 \le t < b$,

$$|\xi_A(t) - \xi_{A_0}(t)| \le \|A - A_0\|_a + C_0|\lambda| \int_0^t |\xi_A(s) - \xi_{A_0}(s)|\, ds.$$

Solving this inequality, we have that for any $0 \le t < b$,

$$|\xi_A(t) - \xi_{A_0}(t)| \le \exp(C_0|\lambda|t)\|A - A_0\|_a \le C\|A - A_0\|_a.$$

So $|\xi_A(b) - \xi_{A_0}(b)| \le C\|A - A_0\|_a < C\delta = \delta_0$. From the definition of $b$, we have $b = a$. Thus $T_A > a$ and $\|\xi_A - \xi_{A_0}\|_a \le C\|A - A_0\|_a$ if $\|A - A_0\|_a < \delta$. So $\{T_A > a\}$ is open w.r.t. $\|\cdot\|_a$, and $A \mapsto \xi_A$ is $(\|\cdot\|_a, \|\cdot\|_a)$ continuous on $\{T_A > a\}$.

(ii) Suppose $\alpha$ is a crosscut in $\mathbb{H}$ such that $\bigcup_{0\le t<T} K_t^\xi \subset H(\alpha) \subset \Omega \setminus \{p\}$. Then $T \le \mathrm{hcap}(H(\alpha))/2 < +\infty$. From the compactness of $\mathcal{H}(\alpha)$, $X_t^\xi$ is bounded on $[0,T)$. So from (3.2), $\xi(t) \to x$ for some $x \in \mathbb{R}$ as $t \to T$. Define $\xi(T) = x$. Then $\xi \in C([0,T])$, $K_T^\xi \subset H(\alpha) \subset \Omega \setminus \{p\}$, and so $J_t^\xi$ is defined for $t \in [0,T]$. Then $\xi(t)$ solves (3.2) for $0 \le t \le T$, which is a contradiction. □

## 6. Convergence of the driving functions.

From now on, we begin proving Theorem 4.2. We first study the case that the target is an interior point. In this section, we will show that the driving functions for the discrete LERW converge to those for the continuous LERW.

6.1. *Some estimates.* Suppose $\Omega$ is an almost $\mathbb{H}$ domain and $p \in \Omega$. We now use the notation in Sections 3 and 4 in the case that the target is an interior point. Let $\alpha$ be a crosscut in $\mathbb{H}$ such that $H(\alpha) \subset \Omega \setminus \{p\}$; and let $F$ be a compact subset of $\Omega \setminus H(\alpha)$. In the lemmas in this subsection, a uniform constant is a number that depends only on $\Omega, p, \alpha, F$. From the compactness of $\mathcal{H}(\alpha)$ (Lemma 5.4), there is a uniform constant $\mathbf{h} > 0$ such that if $K_a^\xi \subset H(\alpha)$, then for any $t \in [0,a]$, $\mathrm{dist}(\varphi_t^\xi(\partial\Omega \setminus \mathbb{R}) \cup \varphi_t^\xi(F), \mathbb{R}) \wedge \mathrm{dist}(\varphi_t^\xi(F), \varphi_t^\xi(\partial\Omega \setminus \mathbb{R})) \ge \mathbf{h}$.

LEMMA 6.1. *There are uniform constants $C_1, C_2 > 0$ such that if $K_a^\xi \subset H(\alpha)$, then for any $t_1 \le t_2 \in [0,a]$ and $z \in F$,*

$$|\varphi_{t_2}^\xi(z) - \varphi_{t_1}^\xi(z)| \le C_1|t_2 - t_1|;$$

$$\left|\varphi_{t_2}^\xi(z) - \varphi_{t_1}^\xi(z) - \frac{2(t_2 - t_1)}{\varphi_{t_1}^\xi(z) - \xi(t_1)}\right| \le C_2|t_2 - t_1|\left(|t_2 - t_1| + \sup_{t\in[t_1,t_2]}\{|\xi(t) - \xi(t_1)|\}\right).$$

PROOF. Suppose $K_a^\xi \subset H(\alpha)$. Then $|\varphi_t^\xi(z) - \xi(t)| \ge \mathbf{h}$ for any $t \in [0,a]$ and $z \in F$. Since $\varphi_{t_2}^\xi(z) - \varphi_{t_1}^\xi(z) = \int_{t_1}^{t_2} \frac{2}{\varphi_t^\xi(z) - \xi(t)}\,dt$, so $|\varphi_{t_2}^\xi(z) - \varphi_{t_1}^\xi(z)| \le$ $C_1|t_2-t_1|$ for any $t_1\le t_2\in[0,a]$ and $z\in F$, where $C_1=2/\mathbf{h}>0$. Thus for $t_1\le t_2\in[0,a]$ and $z\in F$,

$$\Big|\frac{2}{\varphi^\xi_{t_2}(z)-\xi(t_2)}-\frac{2}{\varphi^\xi_{t_1}(z)-\xi(t_1)}\Big|\le \frac{2}{\mathbf{h}^2}(|\varphi^\xi_{t_2}(z)-\varphi^\xi_{t_1}(z)|+|\xi(t_2)-\xi(t_1)|)$$
$$\le 2C_1/\mathbf{h}^2|t_2-t_1|+2/\mathbf{h}^2|\xi(t_2)-\xi(t_1)|.$$

Let $C_2:=2(C_1\vee 1)/\mathbf{h}^2>0$. Then for $t_1\le t_2\in[0,a]$ and $z\in F$,

$$\Big|\varphi^\xi_{t_2}(z)-\varphi^\xi_{t_1}(z)-\frac{2(t_2-t_1)}{\varphi^\xi_{t_1}(z)-\xi(t_1)}\Big|$$
$$=\Big|\int_{t_1}^{t_2}\Big(\frac{2}{\varphi^\xi_{t}(z)-\xi(t)}-\frac{2}{\varphi^\xi_{t_1}(z)-\xi(t_1)}\Big)\,dt\Big|$$
$$\le \int_{t_1}^{t_2}\Big|\frac{2}{\varphi^\xi_{t}(z)-\xi(t)}-\frac{2}{\varphi^\xi_{t_1}(z)-\xi(t_1)}\Big|\,dt$$
$$\le C_2|t_2-t_1|\Big(|t_2-t_1|+\sup_{t\in[t_1,t_2]}\{|\xi(t)-\xi(t_1)|\}\Big).\qquad\square$$

LEMMA 6.2. *For each $n_1\in\{0,1\}$, $n_2,n_3\in\mathbb{Z}_{\ge 0}$, there is a uniform constant $C>0$ depending on $n_1,n_2,n_3$, such that if $K^\xi_a\subset H(\alpha)$, then for any $t\in[0,a]$, $x\in[c_\alpha,d_\alpha]$, and $z\in F$, we have*

$$|\partial_1^{n_1}\partial_2^{n_2}\partial_{3,z}^{n_3}P^\xi(t,x,\varphi^\xi_t(z))|\le C.$$

PROOF. For $K\in H(\alpha)$, $x\in\mathbb{R}$ and $z\in\Omega_K$, let $P(K,x,z)$ be as in Section 4. Since $\partial\Omega_K$ is analytic, so $P(K,x,\cdot)$ extends harmonically across $\partial\Omega_K$. For $K\in H(\alpha)$ and $x,y\in\mathbb{R}$, let $Q_y(K,x,\cdot)$ be defined on $\overline{\Omega_K}\setminus\{x\}$ such that $Q_y(K,x,\cdot)$ is harmonic in $\Omega_K$; vanishes on $\mathbb{R}\setminus\{x\}$; behaves like $\operatorname{Im}\frac{c}{z-x}+O(1)$ near $x$ for some $c\in\mathbb{R}$; $Q_y(K,x,z)=-2\operatorname{Re}(\partial_{3,z}P(K,x,z)\cdot\frac{2}{z-y})$ for $z\in\partial\Omega_K\setminus\mathbb{R}$ and $z=\varphi_K(p)$. From the compactness of $H(\alpha)$, for any $n_2,n_3\in\mathbb{Z}_{\ge 0}$, there is a uniform constant $C>0$ depending on $n_2,n_3$, such that for any $K\in H(\alpha)$, $x,y\in[c_\alpha,d_\alpha]$, and $z\in F$, we have

$$|\partial_2^{n_2}\partial_{3,z}^{n_3}P(K,x,\varphi_K(z))|,\ |\partial_2^{n_2}\partial_{3,z}^{n_3}Q_y(K,x,\varphi_K(z))|\le C.$$

Note $P^\xi(t,x,z)=P(K^\xi_t,x,z)$ and $\partial_1P^\xi(t,x,z)=Q_{\xi(t)}(K^\xi_t,x,z)$, so we are done. $\square$

LEMMA 6.3. *There is a uniform constant $C>0$ such that if $K^\xi_a\subset H(\alpha)$, then for any $t,t'\in[0,a]$, $|X^\xi_t|\le C$ and $|X^\xi_t-X^\xi_{t'}|\le C(|t-t'|+|\xi(t)-\xi(t')|)$.*

PROOF. Suppose $K_a^\xi \subset H(\alpha)$. Let $J^\xi(t,x) = J_t^\xi(x)$. Note that $X_t^\xi = (\partial_{2,z}^2/\partial_{2,z})J^\xi(t,x)$. Since $\xi(t) \in [c_\alpha, d_\alpha]$ for $t \in [0,a]$, so it suffices to prove that there is a uniform constant $C > 0$ such that for any $t \in [0,a]$ and $x \in [c_\alpha, d_\alpha]$, $|\partial_1^{n_1}\partial_{2,z}^{n_2}(\partial_{2,z}^2/\partial_{2,z})J^\xi(t,x)| \le C$ for $n_1, n_2 \in \{0,1\}$. We need to show that $|\partial_{2,z}J^\xi(t,x)|$ is bounded from below by a positive uniform constant, and $|\partial_1^{n_1}\partial_{2,z}^{n_2+1}J^\xi(t,x)|$ is bounded from above by a positive uniform constant. The proof is similar to that of the above lemma. □

LEMMA 6.4. *There is a uniform constant $C > 0$ such that if $K_a^\xi \subset H(\alpha)$, then for any $t_1 \le t_2 \in [0,a]$ and $z \in F$, we have*

$$|\partial_1 P^\xi(t_2, \xi(t_2), \varphi_{t_2}^\xi(z)) - \partial_1 P^\xi(t_1, \xi(t_1), \varphi_{t_1}^\xi(z))| \le C(|t_2 - t_1| + |\xi(t_2) - \xi(t_1)|).$$

PROOF. This follows from Lemma 4.1, and the above three lemmas. □

LEMMA 6.5. *There is a uniform constant $d_1 > 0$ such that if $K_a^\xi \subset H(\alpha)$, then for any $t_1 < t_2 \in [0,a]$ that satisfy $|t_2 - t_1| \le d_1$, and for any $z \in F$, we have*

$$\begin{aligned} &P^\xi(t_2, \xi(t_2), \varphi_{t_2}^\xi(z)) - P^\xi(t_1, \xi(t_1), \varphi_{t_1}^\xi(z)) \\ &\quad = \partial_2 P^\xi(t_1, \xi(t_1), \varphi_{t_1}^\xi(z)) \cdot [(\xi(t_2) - \xi(t_1)) - (t_2 - t_1)X_{t_1}^\xi] \\ &\qquad + 1/2\partial_2^2 P^\xi(t_1, \xi(t_1), \varphi_{t_1}^\xi(z)) \cdot [(\xi(t_2) - \xi(t_1))^2 - 2(t_2 - t_1)] \\ &\qquad + O(A^2) + O(AB) + O(AB^2) + O(B^3), \end{aligned}$$

*where $A := |t_2 - t_1|$, $B := \sup_{s,t \in [t_1,t_2]}\{|\xi(s) - \xi(t)|\}$, and $O(X)$ is some number whose absolute value is bounded by $C|X|$ for some uniform constant $C > 0$.*

PROOF. We may choose a compact subset $F'$ of $\Omega \setminus H(\rho)$ such that $F$ is contained in the interior of $F'$. So from the compactness of $\mathcal{H}(\alpha)$, there is a uniform constant $d_0 > 0$ such that for any $K \in H(\alpha)$, $\mathrm{dist}(\varphi_K(F), \partial\varphi_K(F')) \ge d_0$. Suppose $K_a^\xi \subset H(\alpha)$. From Lemma 6.1, there is a uniform constant $d_1 > 0$ such that if $s, t \in [0,a]$ satisfy $|s - t| \le d_1$, then for any $z \in F$, $[\varphi_s^\xi(z), \varphi_t^\xi(z)] \subset \varphi_s^\xi(F')$.

Fix $z \in F$ and $t_1 < t_2 \in [0,a]$ with $|t_2 - t_1| \le d_1$. Let $P_1 = P^\xi(t_2, \xi(t_2), \varphi_{t_2}^\xi(z))$, $P_2 = P^\xi(t_1, \xi(t_2), \varphi_{t_2}^\xi(z))$, $P_3 = P^\xi(t_1, \xi(t_1), \varphi_{t_2}^\xi(z))$, $P_4 = P^\xi(t_1, \xi(t_1), \varphi_{t_1}^\xi(z))$. Then

$$P^\xi(t_2, \xi(t_2), \varphi_{t_2}^\xi(z)) - P^\xi(t_1, \xi(t_1), \varphi_{t_1}^\xi(z)) = (P_1 - P_2) + (P_2 - P_3) + (P_3 - P_4).$$

Now $P_1 - P_2 = \int_{t_1}^{t_2} \partial_1 P^\xi(t, \xi(t_2), \varphi_{t_2}^\xi(z))\,dt$. Fix any $t \in [t_1, t_2]$. Applying Lemmas 6.1 and 6.2 to $F'$, since $\xi(t), \xi(t_2) \in [c_\alpha, d_\alpha]$ and $[\varphi_t^\xi(z), \varphi_{t_2}^\xi(z)] \subset$

$\varphi_t^\xi(F')$, so we have

$$\partial_1 P^\xi(t,\xi(t_2),\varphi_{t_2}^\xi(z)) - \partial_1 P^\xi(t,\xi(t),\varphi_t^\xi(z)) = O(A)+O(B).$$

Applying Lemma 6.4 to $F$, we have

$$\partial_1 P^\xi(t,\xi(t),\varphi_t^\xi(z)) - \partial_1 P^\xi(t_1,\xi(t_1),\varphi_{t_1}^\xi(z)) = O(A)+O(B).$$

So we get

$$P_1 - P_2 = \partial_1 P^\xi(t_1,\xi(t_1),\varphi_{t_1}^\xi(z))(t_2-t_1) + O(A^2) + O(AB).$$

Applying Lemma 6.2 to $F'$, since $\varphi_{t_2}^\xi(z) \in \varphi_{t_1}^\xi(F')$, so we have

$$\begin{aligned} P_2 - P_3 = {} & \partial_2 P^\xi(t_1,\xi(t_1),\varphi_{t_2}^\xi(z))(\xi(t_2)-\xi(t_1)) \\ & + 1/2\partial_2^2 P^\xi(t_1,\xi(t_1),\varphi_{t_2}^\xi(z))(\xi(t_2)-\xi(t_1))^2 + O(B^3). \end{aligned}$$

Applying Lemmas 6.1 and 6.2 to $F'$, since $[\varphi_{t_1}^\xi(z),\varphi_{t_2}^\xi(z)] \subset \varphi_{t_1}^\xi(F')$, so we have

$$\partial_2^j P^\xi(t_1,\xi(t_1),\varphi_{t_2}^\xi(z)) - \partial_2^j P^\xi(t_1,\xi(t_1),\varphi_{t_1}^\xi(z)) = O(A),$$

for $j=1,2$. Thus

$$\begin{aligned} P_2 - P_3 = {} & \partial_2 P^\xi(t_1,\xi(t_1),\varphi_{t_1}^\xi(z))(\xi(t_2)-\xi(t_1)) \\ & + 1/2\partial_2^2 P^\xi(t_1,\xi(t_1),\varphi_{t_1}^\xi(z))(\xi(t_2)-\xi(t_1))^2 \\ & + O(AB) + O(AB^2) + O(B^3). \end{aligned}$$

Applying Lemmas 6.1 and 6.2 to $F'$, since $[\varphi_{t_1}^\xi(z),\varphi_{t_2}^\xi(z)] \subset \varphi_{t_1}^\xi(F')$, so we have

$$\begin{aligned} P_3 - P_4 &= 2\,\mathrm{Re}(\partial_{3,z} P^\xi(t_1,\xi(t_1),\varphi_{t_1}^\xi(z))(\varphi_{t_2}^\xi(z)-\varphi_{t_1}^\xi(z))) + O(A^2) \\ &= 2\,\mathrm{Re}(\partial_{3,z} P^\xi(t_1,\xi(t_1),\varphi_{t_1}^\xi(z))\frac{2(t_2-t_1)}{\varphi_{t_1}^\xi(z)-\xi(t_1)}) + O(AB) + O(A^2). \end{aligned}$$

The conclusion follows from Lemma 4.1 and the equalities for $P_j - P_{j+1}$, $j=1,2,3$. □

6.2. *Convergence.* We use the notation in Section 4.2. We may choose crosscuts $\rho_j$, $j=0,1,2$, in $D$ such that $H(\rho_0)$ is a neighborhood of $0_+$ in $D$, $H(\rho_0) \subset H(\rho_1) \subset H(\rho_2) \subset D\setminus\{z_e,\infty\}$, and

$$d_0 := \min\{\mathrm{dist}(0,\rho_0), \mathrm{dist}(\rho_0,\rho_1), \mathrm{dist}(\rho_1,\rho_2), \mathrm{dist}(\rho_2,z_e)\} > 0.$$

Now suppose $\delta < d_0$. Then $w_e^\delta \notin H(\rho_2)$ as $|w_e^\delta - z_e| < \delta$, any edge of $D^\delta$ can intersect at most one of $\rho_j$'s, and $(0,\delta] \subset H(\rho_0)$. Thus the LERW curve $q_\delta$

must cross all of these $\rho_j$'s. Let $F_D$ be a compact subset of $D \setminus \{\infty\} \setminus H(\rho_2)$ with nonempty interior. Suppose $f$ maps $D$ conformally onto an almost $\mathbb{H}$ domain $\Omega$ such that $f(0_+) = 0$. Let $p = f(z_e)$, $F_\Omega = f(F_D)$ and $\alpha_j = f(\rho_j)$, $j = 0, 1, 2$. Then $F_\Omega$ is a compact subset of $\Omega$ with nonempty interior; $\alpha_j$'s are crosscuts in $\mathbb{H}$; $\alpha_0$ strictly encloses 0; $\alpha_{j+1}$ strictly encloses $\alpha_j$; and $\{p\}, F_\Omega \subset \Omega \setminus H(\alpha_2)$.

In this subsection, a uniform constant is a number that depends only on $D$, $z_e$, $\rho_0$, $\rho_1$, $\rho_2$, $F_D$, $f$, and some other variables we will specify. We use $O(X)$ to denote a number whose absolute value is bounded by $C|X|$ for some positive uniform constant $C$. We use $o_\delta(X)$ to denote a number whose absolute value is bounded by $C(\delta)|X|$ for some positive uniform constant $C(\delta)$ depending on $\delta$, such that $C(\delta) \to 0$ as $\delta \to 0$.

Let $L^\delta$ denote the set of finite simple lattice paths $X = (X(-1), X(0), \ldots, X(s))$, $s \in \mathbb{N}$, on $D^\delta$, such that $X(-1) = 0$, $X(0) = \delta$, $X(k) \in D$ for $0 \le k \le s$, and $\bigcup_{k=0}^{s} (X(k-1), X(k)] \subset H(\rho_1)$. Let $\mathrm{Set}(X) = \{X(0), \ldots, X(s)\}$, $\mathrm{Tip}(X) = X(s)$, $D_X = D \setminus \bigcup_{k=0}^{s} (X(k-1), X(k)]$; $P_X$ be the generalized Poisson kernel in $D_X$ with the pole at $\mathrm{Tip}(X)$, normalized by $P_X(z_e) = 1$; and $g_X$ be defined on $V(D^\delta)$ such that $g_X \equiv 0$ on $V_\partial(D^\delta) \cup \mathrm{Set}(X) \setminus \{\mathrm{Tip}(X)\}$, $\Delta_{D^\delta} g_X \equiv 0$ on $V_I(D^\delta) \setminus \mathrm{Set}(X)$, and $g_X(w_e^\delta) = 1$.

LEMMA 6.6. *Suppose $G = (V, E)$ is a connected locally finite graph. Suppose $A, B \subset V$ are such that $B$ is finite and $A \cup B$ is reachable. Suppose $h$ is a nonnegative bounded function on $V$ such that $h$ vanishes on $A$, and is discrete harmonic on $V \setminus (A \cup B)$. Then we have $\sum_{w \in A} \Delta_G h(w) = -\sum_{w \in B} \Delta_G h(w)$.*

PROOF. For $w_0 \in B$, let $H_{w_0}$ be the bounded function on $V$, which is discrete harmonic in $V \setminus (A \cup B)$, vanishes on $A \cup B \setminus \{w_0\}$, and equals 1 at $w_0$. Then the lemma holds if $h = H_{w_0}$. Since $h(w) = \sum_{w_0 \in B} h(w_0) H_{w_0}(w)$, so we are done. □

PROPOSITION 6.1. *For any $\varepsilon > 0$, there is $\delta_0 > 0$ such that if $0 < \delta < \delta_0$, then for any $X \in L^\delta$, and any $w \in V(D^\delta) \cap (D \setminus H(\rho_2))$, we have $|g_X(w) - P_X(w)| < \varepsilon$.*

SKETCH OF THE PROOF. Suppose the proposition is not true. Then we can find $\varepsilon_0 > 0$, a sequence of lattice paths $X_n \in L^{\delta_n}$ with $\delta_n \to 0$, and a sequence of points $w_n \in V^{\delta_n} \cap (D \setminus H(\rho_2))$, such that $|g_{X_n}(w_n) - P_{X_n}(w_n)| > \varepsilon_0$ for all $n \in \mathbb{N}$. For simplicity of notation, we write $g_n$ for $g_{X_n}$, $P_n$ for $P_{X_n}$ and $D_n$ for $D_{X_n}$. Let $K_n = f(\bigcup_{j=0}^{l(X_n)} (X_n(j-1), X_n(j)])$. Then $K_n \in \mathcal{H}(\alpha_1)$. Write $\varphi_n$ for $\varphi_{K_n}$ and $\Omega_n$ for $\Omega_{K_n}$. Let $x_n = \varphi_n \circ f(\mathrm{Tip}(X_n))$. Then $x_n \in [c_{\alpha_1}, d_{\alpha_1}]$. Let $Q_n = P_n \circ f^{-1} \circ \varphi_n^{-1}$. Then $Q_n$ is the generalized Poisson kernel in $\Omega_n$ with the pole at $x_n$, normalized by $Q_n(\varphi_n(p)) = 1$. From the compactness of $\mathcal{H}(\alpha_1)$, by passing to a subsequence, we may assume that $K_n \xrightarrow{\mathcal{H}} K_0 \in \mathcal{H}(\alpha_1)$ and $x_n \to x_0 \in [c_{\alpha_1}, d_{\alpha_1}]$. Write $\Omega_0$ for $\Omega_{K_0}$ and $\varphi_0$ for $\varphi_{K_0}$. Let $Q_0$ be the generalized Poisson kernel in $\Omega_0$ with the pole at $x_0$, normalized by $Q_0(\varphi_0(p)) = 1$. Let $D_0 = f^{-1}(\Omega \setminus K_0)$ and $P_0 = Q_0 \circ \varphi_0 \circ f$. Then $P_0$ is the generalized Poisson kernel in $D_0$ with the pole at $f^{-1} \circ \varphi_0^{-1}(x_0)$, normalized by $P_0(z_e) = 1$. Moreover, $D_n \xrightarrow{\text{Cara}} D_0$, and $P_n \xrightarrow{\text{l.u.}} P_0$ in $D_0$.

We extend $g_n$ to $\text{CE}^n g_n$ that is defined on the union of lattice squares of $\delta\mathbb{Z}^2$ at whose four vertices $g_n$ is defined. Applying Harnack's inequality to the positive discrete harmonic function $g_n$, we find that $(\text{CE}^n g_n)$ is locally uniformly continuous in $D_0$. By the Arzela–Ascoli theorem, there is a subsequence of $(\text{CE}^n g_n)$, which converges locally uniformly to some $g_0$ in $D_0 \setminus \{\infty\}$. We may assume that the subsequence is $(\text{CE}^n g_n)$ itself. By applying Harnack's inequality to the discrete partial derivatives of $g_n$, we may assume that the continuation of all discrete partial derivatives of $g_n$ also converge to the corresponding partial derivatives of $g_0$. Then we conclude that $g_0$ is a positive harmonic function in $D_0 \setminus \{\infty\}$.

We may find a sequence of crosscuts $(\gamma^k)$ in $D_0$ such that $(H(\gamma^k))$ is a nesting neighborhood basis of the prime end $f^{-1} \circ \varphi_0^{-1}(x_0)$ in $D_0$, which is the pole of $P_0$. Fix $k \in \mathbb{N}$, for each $n \in \mathbb{N}$, we find a crosscut $\gamma_n^k$ in $D_n$ that bounds a neighborhood $H(\gamma_n^k)$ of $\text{Tip}(X_n)$, such that $\gamma_n^k$ converges to $\gamma^k$ in some sense as $n \to \infty$. For each $k \geq 2$, we may construct some "hook" in the area of $D_0$ between $\gamma^{k-1}$ and $\gamma^{k+1}$ that holds the boundary of $D_0$ and disconnects $\gamma^{k+1}$ from $\gamma^{k-1}$. We use these hooks to prove that the values of $g_n$ outside $H(\gamma^{k+1})$ are uniformly bounded, and $g_n(w) \to 0$ as $n \to \infty$ and $w \to \partial D_n$ in $V(D^{\delta_n}) \cap (D_n \setminus H(\gamma_n^{k+1}))$ in the spherical metric. Thus $g_0(z) \to 0$ as $z \to \widehat{\mathbb{C}} \setminus D_0$ in $D_0 \setminus H(\gamma^{k+1})$ in the spherical metric. Since $(H(\gamma^k))$ is a neighborhood basis of $f^{-1} \circ \varphi_0^{-1}(x_0)$ in $D_0$, so if $\infty \notin D$, then $g_0$ must be a generalized Poisson kernel in $D_0$ with the pole at $f^{-1} \circ \varphi_0^{-1}(x_0)$. Since $g_0(z_e) = \lim_n \text{CE}^n g_n(w_e^\delta) = \lim g_n(w_e^\delta) = 1 = P_0(z_e)$, so $g_0 \equiv P_0$ in $D_0$. The sequence $(w_n)$ has a subsequence $(w_{n_k})$ that converges to some $w_0 \in D$ or tends to $\widehat{\mathbb{C}} \setminus D$ in the spherical metric. In both cases, we can get a contradiction.

If $\infty \in D$, we only need to prove that $g_0$ is also harmonic at $\infty$. From Lemma 6.6, we have

$$\sum_{w \in \text{Set}(X_n) \cup (V(D^{\delta_n}) \cap \partial D)} \Delta_{D^{\delta_n}} g_n(w) = 0.$$

Choose a Jordan curve $\sigma$ in $D$ composed of line segments parallel to the $x$ or $y$ axis, such that $\partial D$ is enclosed by $\sigma$. Let $U(\sigma)$ denote the intersection of $D$ with the domain bounded by $J$. Let $G_n$ be a subgraph of $D^{\delta_n}$ spanned by edges in $D^{\delta_n}$ that is incident to at least one vertex in $U(\sigma)$. Let $A = \mathrm{Set}(X_n) \cup (V(D^{\delta_n}) \cap \partial D)$, and let $B$ be the set of vertices of $G$ in $D \setminus U(\sigma)$. Then from Lemma 6.6, we have

$$\sum_{(w,w')\in\mathcal{P}_\sigma^n} (g_n(w) - g_n(w')) = -\sum_{w\in\mathrm{Set}(X_n)\cup(V(D^{\delta_n})\cap\partial D)} \Delta_{D^{\delta_n}} g_n(w) = 0,$$

where $\mathcal{P}_\sigma^n = \{(w,w') : w \in V(D^{\delta_n}) \cap U(\sigma), w' \in V_I(D^{\delta_n}) \setminus U(\sigma), w \sim w'\}$. Since the discrete partial derivative of $g_n$ converges to the continuous partial derivative of $g_0$, so as $n \to \infty$,

$$\sum_{(w,w')\in\mathcal{P}_\sigma^n} (g_n(w) - g_n(w')) \to \int_\sigma \partial_{\mathbf{n}} g_0(z)\, ds(z).$$

Thus $\int_\sigma \partial_{\mathbf{n}} g_0(z)\, ds(z) = 0$, so $g_0$ is harmonic at $\infty$.

The reader may see Section 5 in [20] for the detailed proof of a similar proposition. □

Let the LERW curve $q_\delta$ on $[-1, \chi_\delta]$ be defined as in Section 4.2. For $-1 \le t \le \chi_\delta$, let $v_\delta(t) = \mathrm{hcap}(f \circ q_\delta((0,t]))/2$. Let $T_\delta = v_\delta(\chi_\delta)$ and $u_\delta = v_\delta^{-1}$. Let $\beta_\delta(t) = f(q_\delta(u_\delta(t)))$, $0 < t \le T_\delta$. Since $f(0_+) = 0$, so $\beta_\delta$ extends continuously to $[0, T_\delta]$ such that $\beta_\delta(0) = 0$. From Proposition 3.2, there is some $\xi_\delta \in C([0, T_\delta])$ such that $\beta_\delta((0,t]) = K_t^{\xi_\delta}$ for $0 \le t \le T_\delta$. For $n \in \mathbb{Z}_{\ge 0}$, let $\mathcal{F}_n$ be the $\sigma$-algebra generated by $\{n \le \chi_\delta\}$ and $q_\delta(j)$, $0 \le j \le n$. Let $n_\infty$ be the first $n$ such that $(q_\delta(n-1), q_\delta(n)]$ intersects $\rho_0$. Then $n_\infty$ is an $\mathcal{F}_n$-stopping time and $\bigcup_{k=0}^{n_\infty}(q_\delta(k-1), q_\delta(k_j)]$ is contained in $H(\rho_1)$ because $\delta < \mathrm{dist}(\rho_0, \rho_1)$. Let $T_{\alpha_0}^\delta = v_\delta(n_\infty)$. So $K_{T_{\alpha_0}^\delta}^{\xi_\delta} \subset H(\alpha_1)$. Then $T_{\alpha_0}^\delta \le \mathrm{hcap}(H(\alpha_1))/2$, so $T_{\alpha_0}^\delta = O(1)$.

Fix any $n \in \mathbb{Z}_{[-1, n_\infty - 1]}$. Then $(q_\delta(n), q_\delta(n+1)]$ can be disconnected from $\rho_1$ by an annulus $A = \{\delta < |z - q_\delta(n)| < d_0\}$. Let $\Gamma$ be the set of all crosscuts $\gamma$ in $D \setminus \bigcup_{k=0}^n [q_\delta(k-1), q_\delta(k)]$ that is contained in $A$, and disconnects $(q_\delta(n), q_\delta(n+1)]$ from $\rho_1$ in $D \setminus \bigcup_{k=0}^n [q_\delta(k-1), q_\delta(k)]$. Then the extremal length of $\Gamma$ is at most $2\pi/\ln(d_0/\delta)$. If $\gamma \in \Gamma$, then $\varphi_{v_\delta(n)}^{\xi_\delta} \circ f(\gamma)$ is a crosscut in $\mathbb{H}$, which disconnects $\varphi_{v_\delta(n)}^{\xi_\delta} \circ f((q_\delta(n), q_\delta(n+1)]) = \varphi_{v_\delta(n)}^{\xi_\delta}(K_{v_\delta(n+1)}^{\xi_\delta} \setminus K_{v_\delta(n)}^{\xi_\delta})$ from $\varphi_{v_\delta(n)}^{\xi_\delta}(\alpha_1)$ in $\mathbb{H}$. Since $K_{v_\delta(n)}^{\xi_\delta} \subset H(\alpha_0)$, and $\alpha_0$ is strictly enclosed by $\alpha_1$, so from the compactness of $\mathcal{H}(\alpha_0)$, the area of $H(\varphi_{v_\delta(n)}^{\xi_\delta}(\alpha_1))$ is bounded from above by a uniform constant $C_0 > 0$. By the conformal invariance, the extremal length of $f(\Gamma)$ is at most $2\pi/\ln(d_0/\delta)$. So there is $\gamma \in f(\Gamma)$ whose length is smaller than $l(\delta) := 2(C_0\pi/\ln(d_0/\delta))^{1/2}$. Then $l(\delta) = o_\delta(1)$. Since $\varphi_{v_\delta(n)}^{\xi_\delta}(K_{v_\delta(n+1)}^{\xi_\delta} \setminus K_{v_\delta(n)}^{\xi_\delta})$ is enclosed by $\gamma$, so its diameter is not bigger than $l(\delta)$. Thus there is $x_0 \in \mathbb{R}$ such that $\varphi_{v_\delta(n)}^{\xi_\delta}(K_{v_\delta(n+1)}^{\xi_\delta} \setminus K_{v_\delta(n)}^{\xi_\delta}) \subset \{z \in \mathbb{H} : |z - x_0| \le l(\delta)\}$. Thus $v_\delta(n+1) - v_\delta(n) \le \mathrm{hcap}(\{z \in \mathbb{H} : |z - x_0| \le l(\delta)\})/2 =$

$l(\delta)^2/2$ and $\xi_\delta(t) \in [x_0 - 2l(\delta), x_0 + 2l(\delta)]$ for any $t \in [v_\delta(n), v_\delta(n+1)]$, which implies that $|\xi_\delta(s) - \xi_\delta(t)| \le 4l(\delta)$ for any $s, t \in [v_\delta(n), v_\delta(n+1)]$.

Now fix a small $d > 0$. Define a nondecreasing sequence $(n_j)_{j\ge 0}$ inductively. Let $n_0 = 0$. Let $n_{j+1}$ be the first $n \ge n_j$ such that $n = n_\infty$, or $v_\delta(n) - v_\delta(n_j) \ge d^2$, or $|\xi_\delta(n) - \xi_\delta(n_j)| \ge d$, whichever comes first. Then $n_j$'s are stopping times w.r.t. $\{\mathcal{F}_n\}$, and are all bounded by $n_\infty$. From the result of the last paragraph, we may let $\delta > 0$ be smaller than some positive uniform constant depending on $d$, such that $v_\delta(n_{j+1}) - v_\delta(n_j) \le 2d^2$ and $|\xi_\delta(v_\delta(s)) - \xi_\delta(v_\delta(n_j))| \le 2d$ for any $s \in [n_j, n_{j+1}]$, $0 \le j < \infty$. Let $\mathcal{F}'_j = \mathcal{F}_{n_j}$, $0 \le j < \infty$. For $0 \le n \le n_\infty$, let $q_\delta^n$ be the subpath of $q_\delta$ up to time $n$; then $q_\delta^n \in L^\delta$. Let $(g_n)$ be the $(g_n)$ in Proposition 2.1 for the LERW $q_\delta$. Then $g_n = g_{q_\delta^n}$, where $g_{q_\delta^n}$ is as in Proposition 6.1. For simplicity, we write $P_n$ for $P_{q_\delta^n}$.

From Proposition 2.1, for any $w \in V(D^\delta) \cap F_D$, $(g_{n_j}(w))_{j\ge 0}$ is a martingale w.r.t. $\{\mathcal{F}'_j\}$, so $\mathbf{E}[g_{n_{j+1}}(w)|\mathcal{F}'_j] = g_{n_j}(w)$ for any $j \in \mathbb{Z}_{\ge 0}$. From Proposition 6.1, we have $\mathbf{E}[P_{n_{j+1}}(w)|\mathcal{F}'_j] = P_{n_j}(w) + o_\delta(1)$. From Harnack's inequality, the absolute values of the gradients of $P_{n_j}$ on $F_D$ are bounded by a positive uniform constant. Since for any $z \in F_D$, there is $w \in V(D^\delta) \cap F_D$ with $|z - w| \le o_\delta(1)$, so for any $z \in F_D$, $\mathbf{E}[P_{n_{j+1}}(z)|\mathcal{F}'_j] = P_{n_j}(z) + o_\delta(1)$. Note that $P_n \circ f^{-1} = P^{\xi_\delta}(v_\delta(n), \xi_\delta(v_\delta(n)), \varphi^{\xi_\delta}_{v_\delta(n)}(\cdot))$. So for any $z \in F_\Omega = f(F_D)$,

$$\begin{aligned} &\mathbf{E}[P^{\xi_\delta}(v_\delta(n_{j+1}), \xi_\delta(v_\delta(n_{j+1})), \varphi^{\xi_\delta}_{v_\delta(n_{j+1})}(z))|\mathcal{F}'_j] \\ &\quad = P^{\xi_\delta}(v_\delta(n_j), \xi_\delta(v_\delta(n_j)), \varphi^{\xi_\delta}_{v_\delta(n_j)}(z)) + o_\delta(1). \end{aligned} \tag{6.1}$$

PROPOSITION 6.2. *There are a uniform constant $d_2 > 0$, and a uniform constant $\delta(d) > 0$ depending on $d$, such that if $d < d_2$ and $\delta < \delta(d)$, then for all $j \in \mathbb{Z}_{\ge 0}$,*

$$\mathbf{E}\left[(\xi_\delta(v_\delta(n_{j+1})) - \xi_\delta(v_\delta(n_j))) - \int_{v_\delta(n_j)}^{v_\delta(n_{j+1})} X_t^{\xi_\delta}\, dt \Big| \mathcal{F}'_j\right] = O(d^3);$$

$$\mathbf{E}[(\xi_\delta(v_\delta(n_{j+1})) - \xi_\delta(v_\delta(n_j)))^2 - 2(v_\delta(n_{j+1}) - v_\delta(n_j))|\mathcal{F}'_j] = O(d^3).$$

PROOF. Note that $K^\xi_{T^\delta_{\alpha_0}} \subset H(\alpha_1)$. Let $d_1 > 0$ be the uniform constant given by Lemma 6.5 with $\alpha = \alpha_1$. Let $d_2 = (d_1/2)^{1/2}$. Suppose $d < d_2$. Fix $j \in \mathbb{Z}_{\ge 0}$. Let $a = v_\delta(n_j)$, $b = v_\delta(n_{j+1})$. Then $0 \le b - a \le 2d^2 \le 2d_2^2 = d_1$, and $|\xi_\delta(s) - \xi_\delta(t)| \le 4d$ for any $s, t \in [a, b]$. Fix $z \in F_\Omega$. From Lemma 6.5, we have

$$\begin{aligned} &P^{\xi_\delta}(b, \xi_\delta(b), \varphi_b^{\xi_\delta}(z)) - P^{\xi_\delta}(a, \xi_\delta(a), \varphi_a^{\xi_\delta}(z)) \\ &\quad = \partial_2 P^{\xi_\delta}(a, \xi_\delta(a), \varphi_a^{\xi_\delta}(z))((\xi_\delta(b) - \xi_\delta(a)) - (b-a)X_a^{\xi_\delta}) \\ &\qquad + \tfrac{1}{2}\partial_2^2 P^{\xi_\delta}(a, \xi_\delta(a), \varphi_a^{\xi_\delta}(z))((\xi_\delta(b) - \xi_\delta(a))^2 - 2(b-a)) + O(d^3). \end{aligned}$$

Take the conditional expectation of this equality with respect to $\mathcal{F}'_j$. From (6.1), we have

$$\partial_2 P^{\xi_\delta}(a,\xi_\delta(a),\varphi_a^{\xi_\delta}(z))\mathbf{E}[(\xi_\delta(b)-\xi_\delta(a))-(b-a)X_a^{\xi_\delta}|\mathcal{F}'_j]$$
$$+\tfrac{1}{2}\partial_2^2 P^{\xi_\delta}(a,\xi_\delta(a),\varphi_a^{\xi_\delta}(z))\mathbf{E}[(\xi_\delta(b)-\xi_\delta(a))^2-2(b-a)|\mathcal{F}'_j]$$
$$=O(d^3)+o_\delta(1).$$

Since $o_\delta(1)\to 0$ uniformly as $\delta\to 0$, so there is a positive uniform function $\delta(d)$ depending only on $d$ such that if $\delta<\delta(d)$, then $|o_\delta(1)|\le d^3$. From Lemma 6.3, we have $X_t^{\xi_\delta}-X_a^{\xi_\delta}=O(d)$ for any $t\in[a,b]$. Thus for $\delta<\delta(d)$,

$$\partial_2 P^{\xi_\delta}(a,\xi_\delta(a),\varphi_a^{\xi_\delta}(z))\mathbf{E}\left[(\xi_\delta(b)-\xi_\delta(a))-\int_a^b X_t^{\xi_\delta}\,dt\Big|\mathcal{F}'_j\right]$$
$$+\tfrac{1}{2}\partial_2^2 P^{\xi_\delta}(a,\xi_\delta(a),\varphi_a^{\xi_\delta}(z))\mathbf{E}[(\xi_\delta(b)-\xi_\delta(a))^2-2(b-a)|\mathcal{F}'_j]=O(d^3).$$

Note that this is true for any $z\in F_\Omega$. We may choose $z_1\ne z_2\in F_\Omega$ and solve the linear equations to get the estimates of the two conditional expectations. We already know that $\partial_2^j P^{\xi_\delta}(a,\xi_\delta(a),\varphi_a^{\xi_\delta}(z))=O(1)$ for $j=1,2$. So the proof will be completed if we prove that there is a uniform positive constant $C_0$ such that there are $z_1,z_2\in F_\Omega$ that satisfy

$$|\partial_2 P^{\xi_\delta}(a,\xi_\delta(a),\varphi_a^{\xi_\delta}(z_1))\cdot\partial_2^2 P^{\xi_\delta}(a,\xi_\delta(a),\varphi_a^{\xi_\delta}(z_2))$$
$$-\partial_2 P^{\xi_\delta}(a,\xi_\delta(a),\varphi_a^{\xi_\delta}(z_2))\cdot\partial_2^2 P^{\xi_\delta}(a,\xi_\delta(a),\varphi_a^{\xi_\delta}(z_1))|\ge C_0.$$

This follows from the compactness of $\mathcal{H}(\alpha_1)$, and the fact that for every $K\in\mathcal{H}(\alpha_1)$ and $x\in[c_{\alpha_1},d_{\alpha_1}]$, there are $z_1,z_2\in F_\Omega$ such that

$$\begin{aligned}&\partial_2 P(K,x,\varphi_K(z_1))\partial_2^2 P(K,x,\varphi_K(z_2))\\&\quad-\partial_2 P(K,x,\varphi_K(z_2))\partial_2^2 P(K,x,\varphi_K(z_1))\ne 0.\end{aligned}\tag{6.2}$$

Here, if (6.2) does not hold for some $K\in\mathcal{H}(\alpha_1)$ and $x\in[c_{\alpha_1},d_{\alpha_1}]$, then there is $C=C(K,x,F_\Omega)$ such that $\partial_2^2 P(K,x,z)=C\partial_2 P(K,x,z)$ for $z\in\varphi_K(F_\Omega)$. Since $\varphi_K(F_\Omega)$ contains an interior point, and $\partial_2^j P(K,x,\cdot)$, $j=1,2$, are harmonic in $\Omega_K$, so $\partial_2^2 P(K,x,z)=C\partial_2 P(K,x,z)$ for $z\in\Omega_K$, which cannot be true because $x$ is a pole of $\partial_2^j P(K,x,\cdot)$ of order $j+1$ for $j=1,2$. □

Let $\eta_\delta(t)=\xi_\delta(t)-2\int_0^t X_s^{\xi_\delta}\,ds$, $0\le t\le T_{\alpha_0}^\delta=v_\delta(n_\infty)$. From Lemma 6.3, we have $\int_{v_\delta(n_j)}^{v_\delta(n_{j+1})}X_s^{\xi_\delta}\,ds=O(d^2)$ for $0\le t\le T_{\alpha_0}^\delta$. Thus

$$\mathbf{E}[(\eta_\delta(v_\delta(n_{j+1}))-\eta_\delta(v_\delta(n_j)))|\mathcal{F}'_j]=O(d^3);$$
$$\mathbf{E}[(\eta_\delta(v_\delta(n_{j+1}))-\eta_\delta(v_\delta(n_j)))^2-2(v_\delta(n_{j+1})-v_\delta(n_j))|\mathcal{F}'_j]=O(d^3).$$

The following theorem can be deduced by using the Skorokhod embedding theorem. It is very similar to Theorem 3.7 in [10], so we omit the proof.

THEOREM 6.1. *For every $\varepsilon > 0$, there is a uniform constant $\delta_0 > 0$ depending on $\varepsilon$ such that if $\delta < \delta_0$, then there is a coupling of the processes $\eta_\delta(t)$ and a Brownian motion $B(t)$ such that*

$$\mathbf{P}[\sup\{|\eta_\delta(t) - \sqrt{2}B(t)| : t \in [0, T^\delta_{\alpha_0}]\} < \varepsilon] > 1 - \varepsilon.$$

Note that for $t \in [0, T^\delta_{\alpha_0}]$, $\xi_\delta(t)$ solves the equation

$$\xi_\delta(t) = \eta_\delta(t) + 2\int_0^t X_s^{\xi_\delta} ds. \tag{6.3}$$

Suppose $B(t)$ is a Brownian motion, and $\xi_0(t)$, $0 \le t < T_0$, is the maximal solution to

$$\xi_0(t) = \sqrt{2}B(t) + 2\int_0^t X_s^{\xi_0} ds. \tag{6.4}$$

Then there is a.s. a simple curve $\beta_0$ such that $\beta_0(0) = 0$, $\beta_0(t) \in \mathbb{H}$ for $0 < t < T_0$, and $K_t^{\xi_0} = \beta_0((0, t])$ for $0 \le t < T$, and there is a continuous increasing function $u_0$ such that $\gamma_0(t) := f^{-1}(\beta_0(u_0^{-1}(t)))$, $0 \le t < S_0 = u_0(T_0)$, is an LERW$(D; 0_+ \to z_e)$ trace.

If $\alpha$ is a crosscut in $\mathbb{H}$, and $\beta$ defined on $[0, T)$ is a curve in $\overline{\mathbb{H}}$, let $T_\alpha(\beta)$ be the first $t$ such that $\beta(t) \in \alpha$, if such $t$ exists; otherwise let $T_\alpha(\beta) = T$. Since $\beta_\delta([0, T^\delta_{\alpha_0}])$ intersects $\alpha_0$, so $T_{\alpha_0}(\beta_\delta) \le T^\delta_{\alpha_0}$.

THEOREM 6.2. *Suppose $\alpha$ is a crosscut in $\mathbb{H}$ that strictly encloses $0$, and $H(\alpha) \subset \Omega \setminus \{p\}$. If $\infty \in D$, we also assume that $f(\infty) \notin H(\alpha)$. For every $\varepsilon > 0$, there is $\delta_0 > 0$ such that if $\delta < \delta_0$, then there is a coupling of the processes $\xi_\delta(t)$ and $\xi_0(t)$ such that*

$$\mathbf{P}[\sup\{|\xi_\delta(t) - \xi_0(t)| : t \in [0, T_\alpha(\beta_\delta) \vee T_\alpha(\beta_0)]\} < \varepsilon] > 1 - \varepsilon. \tag{6.5}$$

*If $\xi_\delta$ or $\xi_0$ is not defined on $[0, T_\alpha(\beta_\delta) \vee T_\alpha(\beta_0)]$, we set the value of* sup *to be $+\infty$.*

PROOF. Let $\rho_j$ and $\alpha_j = f(\rho_j)$ be as in the beginning of this subsection such that $\alpha$ is strictly enclosed by $\alpha_0$. From Lemma 5.5, there is $\delta_1 > 0$ such that if $K_a^\zeta \subset H(\alpha)$ and $\|\zeta - \eta\|_a < \delta_1$, then $K_a^\eta$ is strictly enclosed by $\alpha_0$. Since $K^{\xi_\delta}_{T^\delta_{\alpha_0}}$ intersects $\alpha_0$, so if $\xi_\delta$ and $\xi_0$ are coupled, then on the event that $|\xi_\delta(t) - \xi_0(t)| < \delta_1$ for $0 \le t \le T^\delta_{\alpha_0}$, we have $\beta_0((0, T^\delta_{\alpha_0}]) = K^{\xi_0}_{T^\delta_{\alpha_0}} \not\subset H(\alpha)$, which implies that $T_\alpha(\beta_\delta) \vee T_\alpha(\beta_0) \le T^\delta_{\alpha_0}$. We may assume $\varepsilon < \delta_1$. Then we suffice to prove this theorem with (6.5) replaced by

$$\mathbf{P}[\sup\{|\xi_\delta(t) - \xi_0(t)| : t \in [0, T^\delta_{\alpha_0}]\} < \varepsilon] > 1 - \varepsilon. \tag{6.6}$$

Since $K^{\xi_\delta}_{T^\delta_{\alpha_0}} \subset H(\alpha_1)$, so from Lemmas 5.5 and 5.8, there are $\delta_2, C_1 > 0$ such that for any $t \in [0, T^\delta_{\alpha_0}]$, if $\|\xi_0 - \xi_\delta\|_t < \delta_2$, then $K^{\xi_0}_t \subset H(\alpha_2)$, and

$$|X^{\xi_0}_t - X^{\xi_\delta}_t| \le C_1 \|\xi_0 - \xi_\delta\|_t. \tag{6.7}$$

Let $C_2 = e^{C_1 h_1}/(2C_1)$, where $h_1 = \mathrm{hcap}(H(\alpha_1))$. From Theorem 6.1, there is $\delta_0 > 0$ such that if $\delta < \delta_0$, then there is a coupling of $\eta_\delta$ with $\sqrt{2}B$ such that the probability that $|\eta_\delta(t) - \sqrt{2}B(t)| < (\varepsilon \wedge \delta_2)/C_2$ for $t \in [0, T^\delta_{\alpha_0}]$ is greater than $1 - \varepsilon$. Let $\mathcal{E}^\delta$ denote this event. Assume $\mathcal{E}^\delta$ occurs.

Now $\xi_0(0) = 0 = \xi_\delta(0)$. Let $[0, b)$ be maximal subinterval of $[0, T^\delta_{\alpha_0}) \cap [0, T_0)$, on which $|\xi_0(t) - \xi_\delta(t)| < \varepsilon \wedge \delta_2$. Then from (6.3), (6.4) and (6.7), we have

$$\|\xi_0 - \xi_\delta\|_t \le \|\eta_\delta - \sqrt{2}B\|_{T^\delta_{\alpha_0}} + 2C_1 \int_0^t \|\xi_0 - \xi_\delta\|_s \, ds,$$

for any $t \in [0, b]$. Solving this inequality, since $b \le T^\delta_{\alpha_0} \le h_1/2$ and $\mathcal{E}^\delta$ occurs, so

$$\|\xi_0 - \xi_\delta\|_b \le (e^{2C_1 b} - 1)/(2C_1) \|\eta_\delta - \sqrt{2}B\|_{T^\delta_{\alpha_0}} \le C_2 \|\eta_\delta - \sqrt{2}B\|_{T^\delta_{\alpha_0}} < \varepsilon \wedge \delta_2.$$

Thus $K^{\xi_0}_t \subset H(\alpha_2)$ for $0 \le t < b$. From Theorem 3.1(ii), we have $b < T_0$. Since $\|\xi_0 - \xi_\delta\|_b < \varepsilon \wedge \delta_2$, so $b = T^\delta_{\alpha_0}$. Thus $\xi_0(t)$ is defined on $[0, T^\delta_{\alpha_0}]$, and $|\xi_\delta(t) - \xi_0(t)| < \varepsilon$ for $t \in [0, T^\delta_{\alpha_0}]$ if $\mathcal{E}^\delta$ occurs. So we have (6.6). □

## 7. Convergence of the curves.

7.1. *Local convergence.* We use the notation in Section 4.2. First we introduce a well-known lemma about random walks on $\delta\mathbb{Z}^2$.

LEMMA 7.1. *Suppose $w \in \delta\mathbb{Z}^2$ and $K \subset \mathbb{C}$ is a connected set that satisfies* $\mathrm{diam}(K) \ge R$ [*resp.* $\mathrm{diam}^\#(K) \ge R$]. *Then the probability that a random walk on $\delta\mathbb{Z}^2$ started from $w$ will exit $\mathbf{B}(w; R)$ [resp. $\mathbf{B}^\#(w; R)$] before using an edge of $\delta\mathbb{Z}^2$ that intersects $K$ is at most $C_0((\delta + \mathrm{dist}(w, K))/R)^{C_1}$ [resp. $C_0((\delta + \mathrm{dist}^\#(w, K))/R)^{C_1}$] for some absolute constants $C_0, C_1 > 0$.*

For $w \in V(D^\delta)$, let $X_w$ be a random walk on $D^\delta$ started from $w$, stopped when it hits $V_\partial(D^\delta) \cup \{w^\delta_e\}$. Let $Y_w$ be $X_w$ conditioned to hit $w^\delta_e$. Then $q_\delta = \mathrm{LE}(Y_\delta)$. Lemma 7.1 will be applied because if $w \in D$, $X_w$ is not different from a random walk on $\delta\mathbb{Z}^2$ started from $w$ stopped when it uses an edge that intersects $\partial D$ or hits $w^\delta_e$.

DEFINITION 7.1. Let $z \in \mathbb{C}$, $r, \varepsilon > 0$. A $(z, r, \varepsilon)$-quasi-loop in a path $\omega$ is a pair $a, b \in \omega$ such that $a, b \in \mathbf{B}(z; r)$, $|a - b| \le \varepsilon$, and the subarc of $\omega$ with endpoints $a$ and $b$ is not contained in $\mathbf{B}(z; 2r)$. Let $\mathcal{L}_\delta(z, r, \varepsilon)$ denote the event that $q_\delta$ has a $(z, r, \varepsilon)$-quasi-loop.

LEMMA 7.2. *Suppose* $r>0$ *and* $\mathbf{B}(z_0;5r)\subset D$. *Then* $\mathbf{P}[\mathcal{L}_\delta(z_0,r,\varepsilon)]\to 0$, *as* $\varepsilon\to 0$, *uniformly in* $\delta$.

PROOF. We will use the idea in the proof of Lemma 3.4 in [16]. However, that proof does not apply here immediately, because we are dealing with the loop-erasure of a *conditional* random walk, and Wilson's algorithm does not apply to a *conditional* UST.

We will argue on the reversal path. Let $X^r_w$ be a random walk on $D^\delta$ started from $w$, stopped when it hits $\partial D$. Let $Y^r_w$ be $X^r_w$ conditioned to hit the boundary vertex $\langle\delta,0\rangle$. Let $q^r_\delta=\mathrm{LE}(Y^r_{w^\delta_e})$. Then $q^r_\delta$ has the same distribution as the reversal of $q_\delta$. Let $\mathcal{L}^r_\delta(z_0,r,\varepsilon)$ denote the event that $q^r_\delta$ has a $(z_0,r,\varepsilon)$-quasi-loop. Then $\mathbf{P}[\mathcal{L}^r_\delta(z_0,r,\varepsilon)]=\mathbf{P}[\mathcal{L}_\delta(z_0,r,\varepsilon)]$. It suffices to show that $\lim_{\varepsilon\to0}\mathbf{P}[\mathcal{L}^r_\delta(z_0,r,\varepsilon)]=0$, uniformly in $\delta\in(0,\delta_1]$ for some absolute constant $\delta_1>0$ because if $\delta>\delta_1$, then $\mathcal{L}^r_\delta(z_0,r,\varepsilon)$ does not happen when $\varepsilon<\delta_1$.

Let $\mathbf{B}_k=\mathbf{B}(z_0;kr)$, $k=1,2,3,4,5$. Let $t_0=0$ and $j=0$. If $t_j$ is defined, then define $s_{j+1}$ to be the first time $s>t_j$ such that $Y^r_{w^\delta_e}(s)\in\mathbf{B}_1$, if such $s$ exists; otherwise, let $M=j$ and stop here. If $s_{j+1}$ is defined, then define $t_{j+1}$ to be the first time $t>s_{j+1}$ such that $Y^r_{w^\delta_e}(t)\notin\mathbf{B}_2$. Let $j=j+1$ and iterate the definition. Then we get a sequence $s_1<t_1<\cdots<s_M<t_M$. Such $M$ is a random number. Finally, for each $s\ge0$, let $(Y^r_{w^\delta_e})^s$ be the subpath of $Y^r_{w^\delta_e}$ up to time $s$.

For $j\in\mathbb{N}$, let $\mathcal{Y}_j$ be the event that $j\le M$ and $\mathrm{LE}((Y^r_{w^\delta_e})^{t_j})$ has a $(z_0,r,\varepsilon)$-quasi-loop. Then $\mathcal{Y}_1$ is empty, and it is clear that for any $m\in\mathbb{N}$,

$$\mathcal{L}^r_\delta(z_0,r,\varepsilon)\subset\bigcup_{j=1}^{\infty}\mathcal{Y}_j\subset\{M\ge m+1\}\cup\bigcup_{j=1}^{m}\mathcal{Y}_j. \tag{7.1}$$

We first estimate $\mathbf{P}[M\ge j+1|(Y^r_{w^\delta_e})^{t_j}]$. For $w\in V(D^\delta)$, let $Q(w)$ or $Q^\delta(w)$ be the probability that $X^r_w$ leaves $D$ through $[\delta,0]$; let $Q_1(w)$ or $Q^\delta_1(w)$ be the probability that $X^r_w$ avoids $\mathbf{B}_1$ and leaves $D$ through $[\delta,0]$. Then the probability that $Y^r_w$ does not hit $\mathbf{B}_1$ is equal to $Q_1(w)/Q(w)$. From the Markov property of $Y$, we have

$$\mathbf{P}[M\ge j+1|(Y^r_{w^\delta_e})^{t_j}]=1-Q_1(Y^r_{w^\delta_e}(t_j))/Q(Y^r_{w^\delta_e}(t_j)).$$

Let $F=\{2r\le|z-z_0|\le3r\}$. Then $F$ is a compact subset of $D\setminus\overline{\mathbf{B}_1}$, and if $\delta<r$, then $Y^r_{w^\delta_e}(t_j)\in F$. We claim that there are absolute constants $\delta_0\in(0,r)$ and $C_2>0$ such that $Q_1(w)/Q(w)\ge C_2$ for any $w\in V(D^\delta)\cap F$, if $\delta<\delta_0$. If the claim is not true, then we can find $\delta_n\to0$, $w_n\in V(D^{\delta_n})\cap F$, and $w_n\to w_0\in F$, such that $Q^{\delta_n}_1(w_n)/Q^{\delta_n}(w_n)\to0$. Let $I^{\delta_n}=Q^{\delta_n}(\cdot)/Q^{\delta_n}(w_n)$ and $J^{\delta_n}=(Q^{\delta_n}(\cdot)-Q^{\delta_n}_1(\cdot))/Q^{\delta_n}(w_n)$. Let $P$ be the generalized Poisson kernels in $D$ with the pole at $0_+$, normalized by $P(w_0)=1$. Then $I^{\delta_n}$ converges to $P$ locally uniformly in $D$. Since $J^{\delta_n}$ vanishes on the boundary vertices of $D^\delta$ including 0, agrees with $I^{\delta_n}$ on the vertices in $\mathbf{B}_1$, and is discrete harmonic in $D\setminus\overline{\mathbf{B}_1}$, so $J_n^\delta$ converges to a continuous function $H$ locally uniformly in $\overline{D}\setminus\mathbf{B}_1$, where $H$ vanishes on $\partial D$, agrees with $P$ on $\partial\mathbf{B}_1$, and is harmonic in $D\setminus\overline{\mathbf{B}_1}$. Then $H\le P$ in $D\setminus\overline{\mathbf{B}_1}$. From $J^{\delta_n}(w_n)\to 1$, we have $H(w_0)=1=P(w_0)$. From the maximum principle of harmonic functions, we have $P(w)-H(w)=0$ for any $w\in D\setminus\overline{\mathbf{B}_1}$, which is impossible. So the claim is justified. Suppose $\delta<\delta_0$. Then $\mathbf{P}[M\ge j+1|(Y^r_{w_e^\delta})^{t_j}]\le 1-C_2$. By induction, we find that

$$\mathbf{P}[M\ge m+1]\le (1-C_2)^m. \tag{7.2}$$

We now estimate $\mathbf{P}[\mathcal{Y}_{j+1}|\neg\mathcal{Y}_j,(Y^r_{w_e^\delta})^{t_j}]$. Let $\mathcal{Q}_j$ be the set of components of intersection of $\mathbf{B}_2$ with $\mathrm{LE}((Y^r_{w_e^\delta})^{s_{j+1}})$ that do not contain $Y^r_{w_e^\delta}(s_{j+1})$. Observe that if $\mathcal{Y}_j$ does not occur, then for $\mathcal{Y}_{j+1}$ to occur, there must be a $K\in\mathcal{Q}_j$ such that $Y^r_{w_e^\delta}$ comes at some time $t\in[s_{j+1},t_{j+1}]$ within distance $\varepsilon$ of $K\cap\mathbf{B}_1$ but $Y^r_{w_e^\delta}(t)\notin K$ for all $t\in[s_{j+1},t_{j+1}]$. But if $Y^r_{w_e^\delta}(t)$ is close to $K$ for some $t\in[s_{j+1},t_{j+1}]$, then Lemma 7.1 can be applied, to estimate the probability that $Y^r_{w_e^\delta}(t)$ will not hit $K$ before time $t_{j+1}$.

Suppose $\delta<\delta_1:=\delta_0\wedge\mathrm{dist}(0,\mathbf{B}_5)$; then $\delta\notin\mathbf{B}_5$, so $Q$ is discrete harmonic inside $\mathbf{B}_5$, and $Q(w)>0$ for any $w\in V(D^\delta)\cap\mathbf{B}_5$. Applying Harnack's inequality to $Q$, we get an absolute constant $C_1\ge 1$ such that $Q(w_1)\le C_1Q(w_2)$ for any $w_1,w_2\in V(D^\delta)\cap\mathbf{B}_4$. Let $T_3$ be the first time that a path leaves $\mathbf{B}_3$ or hits $K$. Then for any $w\in V(D^\delta)\cap\mathbf{B}_3$, $X^r_w(t)$ and $Y^r_w(t)$, $t=0,1,\dots,T_3$, are contained in $\mathbf{B}_4$ because $\delta<\delta_0<r$. Note that for any path $(w_0,w_1,\dots,w_n)$ on $D^\delta$ that is contained in $\mathbf{B}_4$,

$$\mathbf{P}[Y^r_{w_0}(j)=w_j,1\le j\le n]\,/\,\mathbf{P}[X^r_{w_0}(j)=w_j,1\le j\le n]=Q(w_n)/Q(w_0)\le C_1.$$

Therefore, conditioned on $Y^r_{w_e^\delta}(s_{j+1})$, for each given $K\in\mathcal{Q}_j$, the probability that $Y^r_{w_e^\delta}([s_{j+1},t_{j+1}])$ gets to within distance $\varepsilon$ of $K$ but does not hit $K$ is at most $C_3((\delta+\varepsilon)/r)^{C_4}$ for some absolute constant $C_3,C_4>0$. Note that if $\delta>\varepsilon$, then the above event cannot happen, so the probability is at most $C_3(2\varepsilon/r)^{C_4}$. Observe that $|\mathcal{Q}_j|$, the cardinality of $\mathcal{Q}_j$, is at most $j$. Let $C_5=C_3(2/r)^{C_4}$. Then

$$\mathbf{P}[\mathcal{Y}_{j+1}|\neg\mathcal{Y}_j]\le jC_5\varepsilon^{C_4}.$$

This gives

$$\begin{aligned}\mathbf{P}\left[\bigcup_{j=1}^m\mathcal{Y}_j\right]&=\sum_{j=1}^{m-1}\mathbf{P}[\mathcal{Y}_{j+1}\cap\neg\mathcal{Y}_j]\le\sum_{j=1}^{m-1}\mathbf{P}[\mathcal{Y}_{j+1}|\neg\mathcal{Y}_j]\\&\le\sum_{j=1}^{m-1}jC_5\varepsilon^{C_4}\le m^2C_5\varepsilon^{C_4}.\end{aligned}$$

Combining this with (7.1) and (7.2), we find that

$$\mathbf{P}[\mathcal{L}^r_\delta(z_0,r,\varepsilon)] \le (1-C_2)^m + m^2 C_5 \varepsilon^{C_4}.$$

Since $C_2>0$, the lemma follows by taking $m=\lfloor \varepsilon^{-C_4/3}\rfloor$, say. □

Definition 7.2. Let $F\subset\mathbb{C}$, and $r,\varepsilon>0$. An $(F,r,\varepsilon)$-quasi-loop in a path $\omega$ is a pair $a,b\in\omega$ such that $a\in F$, $|a-b|<\varepsilon$, and the subarc of $\omega$ with endpoints $a$ and $b$ is not contained in $\mathbf{B}(a;r)$.

Corollary 7.1. *Suppose $F$ is a compact subset of $D\setminus\{\infty\}$, and $r>0$. Then the probability that $q_\delta$ contains an $(F,r,\varepsilon)$-quasi-loop tends to 0 as $\varepsilon\to 0$, uniformly in $\delta$.*

Proof. Let $\mathcal{L}_\delta(F,r,\varepsilon)$ denote this event. We may find $r_0\in(0,r/3)$ and finitely many points $z_1,\dots,z_n\in F$, such that $\mathbf{B}(z_j;5r_0)\subset D$ for each $j\in\mathbb{Z}_{[1,n]}$, and $F\subset\bigcup_{j=1}^n \mathbf{B}(z_j;r_0/2)$. It is easy to check that if $\varepsilon<r_0/2$, then $\mathcal{L}_\delta(F,r,\varepsilon)\subset\bigcup_{j=1}^n\mathcal{L}_\delta(z_j,r_0,\varepsilon)$. The conclusion follows from Lemma 7.2. □

Corollary 7.2. *Suppose $F$ is a compact subset of $\Omega\setminus\{f(\infty)\}$, and $r>0$. Then the probability that $\beta_\delta$ contains an $(F,r,\varepsilon)$-quasi-loop tends to 0 as $\varepsilon\to 0$, uniformly in $\delta$.*

Proof. This follows from the last corollary, and the facts that $f$ maps $D$ conformally onto $\Omega$, $f$ (resp. $f^{-1}$) is uniformly continuous on each compact subset of $D\setminus\{\infty\}$ (resp. $\Omega\setminus\{f(\infty)\}$), and that $\beta_\delta$ is a time-change of $f\circ q_\delta$. □

For a domain $E$ and $\varepsilon>0$, let $\partial^{\#}_\varepsilon E := \{z\in E : \mathrm{dist}^{\#}(z,\widehat{\mathbb{C}}\setminus E)<\varepsilon\}$. For any $\varepsilon>0$ there are $\varepsilon_1,\varepsilon_2>0$ such that $f(\partial^{\#}_{\varepsilon_1}D)\subset\partial^{\#}_\varepsilon\Omega$ and $f^{-1}(\partial^{\#}_{\varepsilon_2}\Omega)\subset\partial^{\#}_\varepsilon D$. In the following lemmas, let $F_D$ (resp. $F_\Omega$) be a compact subset of $D\setminus\{z_e,\infty\}$ [resp. $\Omega\setminus\{p,f(\infty)\}$].

Lemma 7.3. *The probability that $Y_\delta$ or $q_\delta$ visits $\partial^{\#}_\varepsilon D$ after visiting $F_D$ tends to 0 as $\varepsilon,\delta\to 0$.*

Proof. Since $q_\delta$ is the loop-erasure of $Y_\delta$, so we only need to consider $Y_\delta$. By the Markov property of $Y$, we need to prove that the probability that $Y_w$ visits $\partial^{\#}_\varepsilon D$ tends to 0 as $\varepsilon,\delta\to 0$, uniformly in $w\in F_D$. For $w\in V(D^\delta)$, let $Q(w)$ be the probability that $X_w$ visits $w_e^\delta$. Let $P_\varepsilon(w)$ be the probability that $Y_w$ hits $\partial^{\#}_\varepsilon D$. Then $Q(w)P_\varepsilon(w)$ equals the probability that $X_w$ first hits $\partial^{\#}_\varepsilon D$ and then $w_e^\delta$, which is not bigger than $\sup\{Q(w): w\in\partial^{\#}_\varepsilon D\}$.

Choose $z_0\in F_D$. Let $w_0^\delta$ be the vertex of $D^\delta$ closest to $z_0$. As $\delta\to 0$, $Q(\cdot)/Q(w_0^\delta)$ converges to $G(D,z_e;\cdot)/G(D,z_e;z_0)$ uniformly on any subset of $\overline{D}$ bounded away from $z_e$. Thus $\sup\{Q(w): w\in\partial^{\#}_\varepsilon(D)\}/\inf\{Q(w): w\in F_D\}\to 0$ as $\delta,\varepsilon\to 0$. So $P_\varepsilon(w)\to 0$ as $\varepsilon,\delta\to 0$, uniform on $w\in F_D$. □

COROLLARY 7.3. *The probability that $\beta_\delta$ visits $\partial^\#_\varepsilon\Omega$ after $F_\Omega$ tends to* 0 *as* $\varepsilon,\delta\to 0$.

LEMMA 7.4. *For any $\varepsilon>0$, there are $M,\delta_0>0$ such that if $\delta<\delta_0$, then with probability greater than $1-\varepsilon$, $q_\delta$ stays in $\mathbf{B}(0;M)$ after visiting $F_D$.*

PROOF. This follows from Lemma 7.1 and the idea in the proof of Lemma 7.3. □

LEMMA 7.5. *Let $T^\delta_{F_\Omega}$ be the first time that $\beta_\delta$ visits $F_\Omega$. For any $\varepsilon>0$, there are $\varepsilon_0,\delta_0>0$ such that for $\delta<\delta_0$, with probability greater than $1-\varepsilon$, $\beta_\delta$ satisfies that if $|\beta_\delta(t_1)-\beta_\delta(t_2)|<\varepsilon_0$ for some $t_1,t_2\ge T^\delta_{F_\Omega}$, then* $\operatorname{diam}(\beta_\delta([t_1,t_2]))<\varepsilon$.

PROOF. From Lemma 7.4, there are $M,\delta_1>0$ such that if $\delta<\delta_1$, then with probability greater than $1-\varepsilon/3$, $q_\delta$ stays in $\mathbf{B}(0;M)$ after visiting $f^{-1}(F_\Omega)$, so $\beta_\delta$ stays in $f(D\cap\mathbf{B}(0;M))$ after $T^\delta_{F_\Omega}$. Let $\mathcal{E}^\delta_1$ denote this event. From Corollary 7.3, there are $\delta_2,\varepsilon_1>0$ such that if $\delta<\delta_2$, then with probability greater than $1-\varepsilon/3$, $\beta_\delta(t)\in F:=\Omega\setminus\partial^\#_{\varepsilon_1}\Omega$ for $t\ge a$. Let $\mathcal{E}^\delta_2$ denote this event. Let $F_0=F\setminus f(D\cap\{|z|>M\})$. Then $F_0$ is a compact subset of $\Omega\setminus\{f(\infty)\}$, so from Corollary 7.2, there is $\varepsilon_0>0$ such that with probability greater than $1-\varepsilon/3$, $\beta_\delta$ does not contain an $(F_0,\varepsilon/3,\varepsilon_0)$-quasi-loop. Let $\mathcal{E}^\delta_3$ denote this event. Let $\delta_0=\delta_1\wedge\delta_2$ and $\mathcal{E}^\delta=\mathcal{E}^\delta_1\cap\mathcal{E}^\delta_2\cap\mathcal{E}^\delta_3$. Suppose $\delta<\delta_0$. Then $\mathbf{P}[\mathcal{E}^\delta]>1-\varepsilon$. Assume $\mathcal{E}^\delta$ occurs. Suppose $t_1,t_2\ge T^\delta_{F_\Omega}$ and $|\beta_\delta(t_1)-\beta_\delta(t_2)|<\varepsilon_0$. Since $\mathcal{E}^\delta_1$ and $\mathcal{E}^\delta_2$ occur, so $\beta_\delta(t_1)\in F_0$. Since $\mathcal{E}^\delta_3$ occurs, so $\beta_\delta$ does not contain an $(F_0,\varepsilon/3,\varepsilon_0)$-quasi-loop. Thus $\beta_\delta([t_1,t_2])\subset\mathbf{B}(\beta_\delta(t_1);\varepsilon/3)$, whose diameter is less than $\varepsilon$. □

THEOREM 7.1. *Let $\alpha$ be a crosscut in $\mathbb{H}$ that strictly encloses* 0, *such that $H(\alpha)\subset\Omega\setminus\{p,f(\infty)\}$. For every $\varepsilon>0$, there is $\delta_0>0$ depending on $\alpha$ and $\varepsilon$, such that if $\delta<\delta_0$, then there is a coupling of the processes $\beta_\delta(t)$ and $\beta_0(t)$ such that*

$$\mathbf{P}[\sup\{|\beta_\delta(t)-\beta_0(t)|:t\in[0,T_\alpha(\beta_\delta)\vee T_\alpha(\beta_0)]\}<\varepsilon]>1-\varepsilon. \tag{7.3}$$

PROOF. Let $\alpha_0$ be a crosscut in $\mathbb{H}$ that strictly encloses $\alpha$ such that $H(\alpha_0)\subset\Omega\setminus\{p,f(\infty)\}$. Let $d_0=\operatorname{dist}(\alpha,\alpha_0)>0$. Since $\beta_0((0,T_{\alpha_0}(\beta_0)])$ intersects $\alpha_0$, so if $\beta_\delta$ and $\beta_0$ are coupled, then on the event that $|\beta_\delta(t)-\beta_0(t)|<d_0$ for $0\le t\le T_{\alpha_0}(\beta_0)$, we have $\beta_\delta((0,T_{\alpha_0}(\beta_0)])\notin H(\alpha)$, which implies that $T_\alpha(\beta_\delta)\vee T_\alpha(\beta_0)\le T_{\alpha_0}(\beta_0)$. We may assume $\varepsilon<d_0$. Then we suffice to prove this theorem with (7.3) replaced by

$$\mathbf{P}[\sup\{|\beta_\delta(t)-\beta_0(t)|:t\in[0,T_{\alpha_0}(\beta_0)]\}<\varepsilon]>1-\varepsilon. \tag{7.4}$$

Choose a crosscut $\alpha_1$ in $\mathbb{H}$ that strictly encloses $\alpha_0$, such that $H(\alpha_1) \subset \Omega \setminus \{p, f(\infty)\}$. Suppose the theorem is not true; then there exist $\varepsilon_0 > 0$ and a sequence $\delta_n \to 0$, such that for each $\delta_n$, there is no coupling of $\beta_{\delta_n}$ with $\beta_0$ such that (7.4) holds with $\delta = \delta_n$. From Theorem 6.2, and by passing to a subsequence, we may assume that for each $n$, there is a coupling of $\xi_{\delta_n}$ and $\xi_0$ such that

$$\mathbf{P}[\sup\{|\xi_{\delta_n}(t) - \xi_0(t)| : t \in [0, T_{\alpha_1}(\beta_0)]\} \geq 1/2^n] < 1/2^n. \tag{7.5}$$

We may assume that all $\xi_{\delta_n}$ and $\xi_0$ are defined in the same probability space, and (7.5) is satisfied. By discarding a null event, we have

$$\|\xi_{\delta_n} - \xi_0(t)\|_{T_{\alpha_1}(\beta_0)} \to 0. \tag{7.6}$$

Fix any $t \in [0, T_{\alpha_1}(\beta_0)]$. Suppose $F$ is any compact subset of $\mathbb{H} \setminus \beta_0((0,t])$. From $\|\xi_{\delta_n} - \xi_0(t)\|_t \to 0$, we see that $\varphi_t^{\xi_{\delta_n}} \to \varphi_t^{\xi_0}$ uniformly on $F$, and $F \subset \mathbb{H} \setminus \beta_{\delta_n}((0,t])$ for all but finitely many $n$. Thus $(\mathbb{H} \setminus \beta_{\delta_n}((0,t])) \cap (\mathbb{H} \setminus \beta_0((0,t])) \xrightarrow{\text{Cara}} \mathbb{H} \setminus \beta_0((0,t])$. From Lemma 5.1, $(\varphi_t^{\xi_{\delta_n}})^{-1} \xrightarrow{\text{l.u.}} (\varphi_t^{\xi_0})^{-1}$ in $\mathbb{H} = \varphi_t^{\xi_0}(\mathbb{H} \setminus \beta_0((0,t]))$. Thus we have $\mathbb{H} \setminus \beta_{\delta_n}((0,t]) \xrightarrow{\text{Cara}} \mathbb{H} \setminus \beta_0((0,t])$ for any $t \in [0, T_{\alpha_1}(\beta_0)]$.

We may assume that $\mathbf{B}(0; \varepsilon_0) \cap \mathbb{H} \subset H(\alpha_0)$. Since $\beta_0$ is a continuous curve started from 0, so there is $b > 0$ such that with probability greater than $1 - \varepsilon_0/5$, $\beta_0$ is defined on $[0, b]$, and $\beta_0([0, b]) \subset \mathbf{B}(0; \varepsilon_0/4)$. Let $\mathcal{E}_1^0$ denote this event. If $\mathcal{E}_1^0$ occurs, then $b < T_{\alpha_1}(\beta_0)$. For each $n \in \mathbb{N}$, let $\mathcal{E}_1^n$ denote the event that $\beta_n$ is defined on $[0, b]$ and $\beta_n([0, b]) \subset \mathbf{B}(0; \varepsilon_0/3)$. From (7.6) and Lemma 5.5, we have $\mathcal{E}_1^0 \subset \liminf \mathcal{E}_1^n$. So there is $N_1 \in \mathbb{N}$ such that $\mathbf{P}[\mathcal{E}_1^n] > 1 - \varepsilon_0/5$ if $n > N_1$.

Let $a = b/2$. If $\mathcal{E}_1^0$ occurs, then $\beta_0((0, a]) \subset H(\alpha_1) \subset \Omega \setminus \{p, f(\infty)\}$. So there is a nonempty compact subset $F_1$ of $\Omega \setminus \{p, f(\infty)\}$ such that $\mathbf{P}[\mathcal{E}_2^0] > 1 - \varepsilon_0/5$, where $\mathcal{E}_2^0$ is the subevent of $\mathcal{E}_1^0$ on which $\beta_0((0, a]) \cap F_1 \neq \varnothing$. Choose another compact subset $F_2$ of $\Omega \setminus \{p, f(\infty)\}$ such that $F_1$ is contained in the interior of $F_2$. Let $\mathcal{E}_2^n$ denote the event that $\beta_{\delta_n}$ is defined on $[0, a]$, and $\beta_{\delta_n}((0, a]) \cap F_2 \neq \varnothing$. If $\mathcal{E}_2^0$ occurs, then $a \leq T_{\alpha_1}(\beta_0)$, so $\mathbb{H} \setminus \beta_{\delta_n}((0, a]) \xrightarrow{\text{Cara}} \mathbb{H} \setminus \beta_0((0, a])$, and so $\operatorname{dist}(z_0, \beta_{\delta_n}((0, a])) \to 0$ for any $z_0 \in \beta_0((0, a])$. Thus $\mathcal{E}_2^0 \subset \liminf \mathcal{E}_2^n$. So there is $N_2 \in \mathbb{N}$ such that $\mathbf{P}[\mathcal{E}_2^n] > 1 - \varepsilon_0/5$ if $n > N_2$. Note that if $\mathcal{E}_2^n$ occurs, then $a \geq T_{F_2}^{\delta_n}$, where $T_{F_2}^{\delta_n}$ is the first time that $\beta_{\delta_n}$ visits $F_2$.

From Theorem 6.2 and Lemma 7.5, there are $\varepsilon_1 \in (0, \varepsilon_0)$ and $N_3 \in \mathbb{N}$ such that if $n \geq N_3$, then with probability at least $1 - \varepsilon_0/5$, $\xi_{\delta_n}$ is defined on $[0, T_{\alpha_1}(\beta_0)]$, and if $|\beta_{\delta_n}(t_2) - \beta_{\delta_n}(t_1)| < \varepsilon_1$ for some $t_1, t_2 \geq T_{F_2}^{\delta_n}$, then $\operatorname{diam}(\beta_{\delta_n}([t_1, t_2])) < \varepsilon_0/3$. Let $\mathcal{E}_3^n$ denote this event.

Since $\beta_0$ is continuous on $[a, T_{\alpha_1}(\beta_0)]$, $\operatorname{dist}(\beta_0([a, T_{\alpha_1}(\beta_0)]), \mathbb{R}) > 0$ and $T_{\alpha_0}(\beta_0) < T_{\alpha_1}(\beta_0)$, so there is $\Delta, h > 0$ such that with probability at least $1 - \varepsilon_0/5$, the followings hold: $T_{\alpha_1}(\beta_0) - T_{\alpha_0}(\beta_0) > \Delta$, $\operatorname{Im} \beta_0(t) \geq h$ for any

$t \in [a, T_{\alpha_1}(\beta_0)]$, and if $t_1, t_2 \in [a, T_{\alpha_1}(\beta_0)]$ and $|t_1 - t_2| \le \Delta$, then $|\beta_0(t_1) - \beta_0(t_2)| < \varepsilon_1/3$. Let $\mathcal{E}_4$ denote this event.

Let $A = \text{hcap}(H(\alpha_1))/2$. Then $T_{\alpha_1}(\beta_0) \le A$. Choose $N \in \mathbb{N}$ such that $A/N < (\Delta \wedge b)/2$, and define $t_k = a + (T_{\alpha_1}(\beta_0) - a)k/N$, $k = 0, 1, \ldots, N$. Then $t_0 = a$, $t_N = T_{\alpha_1}(\beta_0)$ and $t_1 \le b$, $t_{N-1} \ge T_{\alpha_0}(\beta_0)$. Fix $k \in \mathbb{Z}_{[1,N]}$. Since $\beta_0(t_k) \in \mathbb{H} \setminus \beta_0((0, t_{k-1}])$ and $\mathbb{H} \setminus \beta_{\delta_n}((0, t_{k-1}]) \xrightarrow{\text{Cara}} \mathbb{H} \setminus \beta_0((0, t_{k-1}])$, so there is $M_k^1 \in \mathbb{N}$ such that $\beta_0(t_k) \notin \beta_{\delta_n}((0, t_{k-1}])$ when $n > M_k^1$. Since $\beta_0(t_k)$ is a boundary point of $\mathbb{H} \setminus \beta_0((0, t_k])$ and $\mathbb{H} \setminus \beta_{\delta_n}((0, t_k]) \xrightarrow{\text{Cara}} \mathbb{H} \setminus \beta_0((0, t_k])$, so there is $M_k^2 \in \mathbb{N}$ such that when $n > M_k^2$, there is $z_n \in \partial(\mathbb{H} \setminus \beta_{\delta_n}((0, t_k]))$ with $|z_n - \beta_0(t_k)| < (\varepsilon_1/3) \wedge h$. If event $\mathcal{E}_4$ occurs, and $n > M_k^1 \vee M_k^2$, then $z_n \notin \mathbb{R}$ and $z_n \notin \beta((0, t_{k-1}])$, which implies that $z_n = \beta_{\delta_n}(s_k)$ for some $s_k \in (t_{k-1}, t_k]$. Thus if $\mathcal{E}_4$ occurs and $n > M := \bigvee_{k=1}^N (M_k^1 \vee M_k^2)$, then we have $s_k \in (t_{k-1}, t_k]$, $k = 1, 2, \ldots, N$, such that $|\beta_{\delta_n}(s_k) - \beta_0(t_k)| < \varepsilon_1/3$.

Let $L = \bigvee_{j=1}^3 N_j \vee M$ and $\mathcal{E}^n = \bigcap_{j=1}^3 \mathcal{E}_j^n \cap \mathcal{E}_1^0 \cap \mathcal{E}_4$. Suppose $n > L$. Then $\mathbf{P}[\mathcal{E}^n] > 1 - \varepsilon_0$. Assume $\mathcal{E}^n$ occurs. Fix $t \in [0, T_{\alpha_0}(\beta_0)]$. If $t \le b$, then $\beta_{\delta_n}(t)$, $\beta_0(t) \in \mathbf{B}(0; \varepsilon_0/3)$ because $\mathcal{E}_1^0$ and $\mathcal{E}_1^n$ both occur and $n > N_1$, so $|\beta_{\delta_n}(t) - \beta_0(t)| < \varepsilon_0$. Now suppose $t \ge b$. Then $t \in [b, T_{\alpha_0}(\beta_0)] \subset [t_1, t_{N-1}] \subset [s_1, s_N]$. Thus $t \in [s_k, s_{k+1}]$ for some $k \in \mathbb{Z}_{[1,N-1]}$. Since $n > M$, $t_k, t_{k+1} \in [a, T_{\alpha_1}(\beta_0)]$, $|t_k - t_{k+1}| < \Delta$, and $\mathcal{E}_4$ occurs, so

$$\begin{aligned} |\beta_{\delta_n}(s_k) - \beta_{\delta_n}(s_{k+1})| &\le |\beta_{\delta_n}(s_k) - \beta_0(t_k)| + |\beta_0(t_k) - \beta_0(t_{k+1})| \\ &\quad + |\beta_0(t_{k+1}) - \beta_{\delta_n}(s_{k+1})| \\ &< \varepsilon_1/3 + \varepsilon_1/3 + \varepsilon_1/3 = \varepsilon_1. \end{aligned}$$

Since $n > N_2$ and $\mathcal{E}_2^n$ occurs, so $s_k, s_{k+1} \ge a \ge T_{F_2}^{\delta_n}$. Since $n > N_3$, $\mathcal{E}_3^n$ occurs, and $t \in [s_k, s_{k+1}]$, so $|\beta_{\delta_n}(t) - \beta_{\delta_n}(s_k)| < \varepsilon_0/3$. Since $t \in [s_k, s_{k+1}] \subset [t_{k-1}, t_{k+1}]$, so $|t - t_k| < \Delta$. Since $\mathcal{E}_4$ occurs, so $|\beta_0(t) - \beta_0(t_k)| < \varepsilon_1/3$. Thus

$$\begin{aligned} |\beta_{\delta_n}(t) - \beta_0(t)| &\le |\beta_{\delta_n}(t) - \beta_{\delta_n}(s_k)| + |\beta_{\delta_n}(s_k) - \beta_0(t_k)| + |\beta_0(t_k) - \beta_0(t)| \\ &\le \varepsilon_0/3 + \varepsilon_1/3 + \varepsilon_1/3 < \varepsilon_0/3 + \varepsilon_0/3 + \varepsilon_0/3 = \varepsilon_0. \end{aligned}$$

Thus with probability greater than $1 - \varepsilon_0$, $|\beta_{\delta_n}(t) - \beta_0(t)| < \varepsilon_0$ for $0 \le t \le T_{\alpha_0}(\beta_0)$, which contradicts the choice of $(\delta_n)$. □

7.2. *Global convergence.* We restrict $\beta_\delta$ to $[0, T_\delta)$. Then $\lim_{t \to T_\delta} \beta_\delta(t) = f(w_e^\delta)$. Recall that $\beta_0$ is defined on $[0, T_0)$, where $[0, T_0)$ is the maximal interval on which the solution to (6.4) exists. Let $\mathcal{B}$ denote the set of continuous curves $\beta : [0, T(\beta)) \to \Omega \cup \mathbb{R}$, for some $T(\beta) \in (0, \infty]$, with $\beta(0) = 0$ and $\beta(t) \in \Omega$ for $t \in (0, T(\beta))$. So $T$ is a function taking values in $(0, \infty]$ on $\mathcal{B}$ that describes the length of lifetime. Then $\beta_0$ and $\beta_\delta$ are $\mathcal{B}$-valued random variables, and $T(\beta_\delta) = T_\delta$, $T(\beta_0) = T_0$.

Let $\mathcal{A}$ denote the set of crosscuts $\alpha$ in $\mathbb{H}$ that strictly enclose 0, and such that $H(\alpha) \subset \Omega \setminus \{p, f(\infty)\}$. For $\alpha_1, \alpha_2 \in A$, we write $\alpha_1 \prec \alpha_2$ or $\alpha_2 \succ \alpha_1$ if $\alpha_1$ is strictly enclosed by $\alpha_2$. For any $\beta \in \mathcal{B}$ and $\alpha \in \mathcal{A}$, let $T_\alpha(\beta)$ be the biggest $T \in (0, T(\beta)]$ such that $\beta(t) \notin \alpha$ for $0 \le t < T$. It is clear that $T_{\alpha_1} \le T_{\alpha_2}$ if $\alpha_1 \prec \alpha_2$. Define $T_\alpha^+ = \bigwedge_{\alpha' \succ \alpha} T_{\alpha'}$. If $\beta$ does not leave $H(\alpha)$ immediately after hitting $\alpha$, then $T_\alpha(\beta) < T_\alpha^+(\beta)$.

Suppose $\alpha \in \mathcal{A}$. For $\beta_1, \beta_2 \in \mathcal{B}$, let $\Delta(\beta_1, \beta)$ be 0 if $\beta_1 = \beta_2$ and 1 otherwise, where $\beta_1 = \beta_2$ means that $T(\beta_1) = T(\beta_2)$ and $\beta_1(t) = \beta_2(t)$ for $0 \le t < T(\beta_1)$, and define

$$d_\alpha^\vee(\beta_1, \beta_2) = \Delta(\beta_1, \beta_2) \wedge \sup\{|\beta_1(t) - \beta_2(t)| : t \in [0, T_\alpha(\beta_1) \vee T_\alpha(\beta_2)]\},$$

where the value of the sup is set to be $\infty$ if either $\beta_1(t)$ or $\beta_2(t)$ is not defined at some $t$ in the interval of the formula. Then $0 \le d_\alpha^\vee \le 1$. Now define

$$d_\alpha(\beta_1, \beta_2) = \inf\left\{\sum_{k=1}^n d_\alpha^\vee(\gamma_{k-1}, \gamma_k) : \gamma_0 = \beta_1, \gamma_n = \beta_2, \gamma_k \in \mathcal{B}, k \in \mathbb{Z}_{[1,n-1]}, n \in \mathbb{N}\right\}.$$

Then $d_\alpha$ is a pseudo-metric on $\mathcal{B}$, and $d_\alpha \le d_\alpha^\vee$. For $\alpha \in \mathcal{A}$, $\beta_1 \in \mathcal{B}$ and $r > 0$, let $\mathbf{B}_\alpha(\beta_1; r) = \{\beta \in \mathcal{B} : d_\alpha(\beta, \beta_1) < r\}$. Let $\mathcal{T}_\alpha$ denote the topology generated by $d_\alpha$. It is clear that if $\alpha_1 \prec \alpha_2$ , then $d_{\alpha_1}^\vee \le d_{\alpha_2}^\vee$, so $d_{\alpha_1} \le d_{\alpha_2}$, from which follows that $\mathcal{T}_{\alpha_1} \subset \mathcal{T}_{\alpha_2}$. Let $\mathcal{T}_\alpha^+ = \bigcap_{\alpha' \succ \alpha} \mathcal{T}_{\alpha'}$.

LEMMA 7.6. *Suppose* $\alpha_1 \prec \alpha_2 \in \mathcal{A}$ *and* $d_0 = 1 \wedge \operatorname{dist}(\alpha_1, \alpha_2) > 0$. *Suppose* $\beta_1, \beta_2 \in \mathcal{B}$, *and* $d_{\alpha_2}(\beta_1, \beta_2) < d_0$. *Then* $d_{\alpha_1}^\vee(\beta_1, \beta_2) \le d_{\alpha_2}(\beta_1, \beta_2)$.

PROOF. Choose $d_1 \in (d_{\alpha_2}(\beta_1, \beta_2), d_0)$. Then there are $\gamma_0, \gamma_1, \dots, \gamma_n \in \mathcal{B}$ such that $\gamma_0 = \beta_1$, $\gamma_n = \beta_2$ and $\sum_{j=1}^n d_{\alpha_2}^\vee(\gamma_{j-1}, \gamma_j) < d_1$. For each $j \in \mathbb{Z}_{[1,n]}$, since $d^\vee(\gamma_{j-1}, \gamma_j) < d_1 < 1$, so

$$d_{\alpha_2}^\vee(\gamma_{j-1}, \gamma_j) = \sup\{|\gamma_{j-1}(t) - \gamma_j(t)| : 0 \le t \le T_{\alpha_2}(\gamma_{j-1}) \vee T_{\alpha_2}(\gamma_j)\}.$$

Let $t_0 = T_{\alpha_1}(\beta_1) \vee T_{\alpha_1}(\beta_2)$. Assume, for example, that $t_0 = T_{\alpha_1}(\beta_1) = T_{\alpha_1}(\gamma_0)$.

We claim that $t_0 \le T_{\alpha_2}(\gamma_j)$ for any $0 \le j \le n$. Since $t_0 = T_{\alpha_1}(\gamma_0) < T_{\alpha_2}(\gamma_0)$, if the claim is not true, then there is $k \in \mathbb{Z}_{[1,n]}$ such that $t_0 > T_{\alpha_2}(\gamma_k)$ and $t_0 \le T_{\alpha_2}(\gamma_j)$ for $0 \le j \le k-1$. Let $t_1 = T_{\alpha_2}(\gamma_k)$. So $t_1 \in [0, T_{\alpha_2}(\gamma_j)]$, $0 \le j \le k$. Then we have

$$d_0 > d_1 > \sum_{j=1}^k d_{\alpha_2}^\vee(\gamma_{j-1}, \gamma_j) \ge \sum_{j=1}^k 1 \wedge |\gamma_{j-1}(t_1) - \gamma_j(t_1)| \ge 1 \wedge |\gamma_0(t_1) - \gamma_k(t_1)|.$$

Since $t_0$ is the first $t$ such that $\beta_1(t) \in \alpha_1$, and $t_1 < t_0$, so $\gamma_0(t_1) = \beta_1(t_1)$ is enclosed by $\alpha_1$. Since $\gamma_k(t_1) \in \alpha_2$ and $\alpha_1 \prec \alpha_2$, so $|\gamma_0(t_1) - \gamma_k(t_1)| \ge \operatorname{dist}(\alpha_1, \alpha_2)$. This implies that $1 \wedge |\gamma_0(t_1) - \gamma_k(t_1)| \ge d_0$, which is a contradiction. So the claim is justified.

Thus for any $t \in [0, t_0]$, we have $t \in [0, T_{\alpha_2}(\gamma_j)]$ for any $0 \le j \le n$, so

$$|\beta_1(t) - \beta_2(t)| \le \sum_{j=1}^{n} |\gamma_{j-1}(t) - \gamma_j(t)| \le \sum_{j=1}^{n} d^{\vee}_{\alpha_2}(\gamma_{j-1}, \gamma_j) < d_1.$$

Since this is true for any $t \in [0, t_0] = [0, T_{\alpha_1}(\beta_1) \vee T_{\alpha_1}(\beta_2)]$ and $d_1 \in (d_{\alpha_2}(\beta_1, \beta_2), d_0)$, so $d^{\vee}_{\alpha_1}(\beta_1, \beta_2) \le d_{\alpha_2}(\beta_1, \beta_2)$. □

LEMMA 7.7. *$\{T^+_{\alpha_1} < T_{\alpha_2}\} \in \mathcal{T}^+_{\alpha_1}$ for any $\alpha_1, \alpha_2 \in \mathcal{A}$.*

PROOF. Fix any $\alpha''_1 \in \mathcal{A}$ such that $\alpha''_1 \succ \alpha_1$. There is $\alpha'_1 \in \mathcal{A}$ with $\alpha''_1 \succ \alpha'_1 \succ \alpha_1$. Suppose $\beta_1 \in \{T^+_{\alpha_1} < T_{\alpha_2}\}$. Then there is $a > 0$ such that $a < T_{\alpha'_1}(\beta_1) \wedge T_{\alpha_2}(\beta_1)$ and $\beta_1(a) \notin H(\alpha_1)$. Let $d_0 = 1 \wedge \operatorname{dist}(\beta_1(a), H(\alpha_1)) \wedge \operatorname{dist}(\alpha'_1, \alpha''_1) \wedge \operatorname{dist}(\beta_1([0,a]), \alpha_2) > 0$. Suppose $\beta_2 \in \mathbf{B}_{\alpha''_1}(\beta_1; d_0)$. From Lemma 7.6, $d^{\vee}_{\alpha'_1}(\beta_2, \beta_1) < d_0$. Since $d_0 \le 1$, so $|\beta_2(t) - \beta_1(t)| < d_0$ for $0 \le t \le T_{\alpha'_1}(\beta_1)$. Since $a < T_{\alpha'_1}(\beta_1)$, so $|\beta_2(t) - \beta_1(t)| < d_0$ for any $t \in [0, a]$. Since $\beta_1([0,a])$ is strictly enclosed by $\alpha_2$, and $d_0 \le \operatorname{dist}(\beta_1([0,a]), \alpha_2)$, so $\beta_2([0,a])$ is also strictly enclosed by $\alpha_2$, which implies that $a < T_{\alpha_2}(\beta_2)$. Since $|\beta_2(a) - \beta_1(a)| < d_0$ and $d_0 \le \operatorname{dist}(\beta_1(a), H(\alpha_1))$, so $\beta_2(a) \notin H(\alpha_1)$, which implies that $T^+_{\alpha_1}(\beta_2) < a$. Thus $T^+_{\alpha_1}(\beta_2) < T_{\alpha_2}(\beta_2)$, that is, $\beta_2 \in \{T^+_{\alpha_1} < T_{\alpha_2}\}$. So $\mathbf{B}_{\alpha''_1}(\beta_1; d_0) \subset \{T^+_{\alpha_1} < T_{\alpha_2}\}$. Thus $\{T^+_{\alpha_1} < T_{\alpha_2}\} \in \mathcal{T}_{\alpha''_1}$. Since $\alpha''_1 \succ \alpha_1$ is chosen arbitrarily, so $\{T^+_{\alpha_1} < T_{\alpha_2}\} \in \mathcal{T}^+_{\alpha_1}$. □

LEMMA 7.8. *Suppose $\alpha_1, \alpha_2 \in \mathcal{A}$ and $B \in \mathcal{T}^+_{\alpha_1}$. Then $B \cap \{T^+_{\alpha_1} < T_{\alpha_2}\} \in \mathcal{T}_{\alpha_2}$.*

PROOF. Fix $\beta_1 \in B \cap \{T^+_{\alpha_1} < T_{\alpha_2}\}$. Then there is $a > 0$ such that $a < T_{\alpha_2}(\beta_1)$ and $\beta_1(a) \notin H(\alpha_1)$. We may choose $\alpha'_1 \succ \alpha$ and $\alpha'_2 \prec \alpha_2$ such that $\beta_1(a) \notin H(\alpha'_1)$ and $\beta_1([0,a])$ is strictly enclosed by $\alpha'_2$. Since $B \in \mathcal{T}^+_{\alpha_1} \subset \mathcal{T}_{\alpha'_1}$, so there is $d_0 > 0$ such that $\mathbf{B}_{\alpha'_1}(\beta_1; d_0) \subset B$. Let $d_1 = 1 \wedge d_0 \wedge \operatorname{dist}(\beta_1(a), H(\alpha'_1)) \wedge \operatorname{dist}(\alpha'_2, \alpha_2)$. Suppose $\beta_2 \in \mathbf{B}_{\alpha_2}(\beta_1; d_1)$. From Lemma 7.6, $d^{\vee}_{\alpha'_2}(\beta_2, \beta_1) < d_1$. Since $d_1 \le 1$, so $|\beta_2(t) - \beta_1(t)| < d_1$ for $0 \le t \le T_{\alpha'_2}(\beta_1)$. Since $a < T_{\alpha'_2}(\beta_1)$, so $|\beta_2(t) - \beta_1(t)| < d_1$ for $0 \le t \le a$. Since $d_1 \le \operatorname{dist}(\beta_1(a), H(\alpha'_1))$, so $\beta_2(a) \notin H(\alpha'_1)$. Thus $T_{\alpha'_1}(\beta_2) \vee T_{\alpha'_1}(\beta_1) < a$. So we have

$$d_{\alpha'_1}(\beta_2, \beta_1) \le d^{\vee}_{\alpha'_1}(\beta_2, \beta_1) \le \sup\{|\beta_2(t) - \beta_1(t)| : 0 \le t \le a\} < d_1 \le d_0.$$

Thus $\beta_2 \in \mathbf{B}_{\alpha'_1}(\beta_1; d_0) \subset B$. Since $\beta_1([0,a])$ is strictly enclosed by $\alpha'_2$, $\alpha'_2 \prec \alpha_2$, and $|\beta_2(t) - \beta_1(t)| < d_1 \le \operatorname{dist}(\alpha'_2, \alpha_2)$ for $0 \le t \le a$, so $\beta_2([0,a])$ is strictly enclosed by $\alpha_2$. Thus $T^+_{\alpha_1}(\beta_2) \le T_{\alpha'_1}(\beta_2) < a < T_{\alpha_2}(\beta_2)$, that is, $\beta_2 \in \{T^+_{\alpha_1} < T_{\alpha_2}\}$. So $\mathbf{B}_{\alpha_2}(\beta_1; d_1) \subset B \cap \{T^+_{\alpha_1} < T_{\alpha_2}\}$. Thus $B \cap \{T^+_{\alpha_1} < T_{\alpha_2}\} \in \mathcal{T}_{\alpha_2}$. □

COROLLARY 7.4. *$\{T^+_{\alpha_1} < T_{\alpha_2}\} \in \mathcal{T}_{\alpha_2}$ for any $\alpha_1, \alpha_2 \in \mathcal{A}$.*

Let $\mu_\delta$ and $\mu_0$ be the distribution of $\beta_\delta$ and $\beta_0$, respectively. From Theorem 7.1, for any $\alpha\in\mathcal{A}$, $\mu_\delta\to\mu_0$ weakly w.r.t. $d_\alpha$, as $\delta\to 0$. Suppose $A$ is a nonempty finite subset of $\mathcal{A}$. Let $d_A=\bigvee_{\alpha\in A}d_\alpha$ and $\mathcal{T}_A$ be the topology generated by $d_A$. So $\mathcal{T}_A=\bigvee_{\alpha\in A}\mathcal{T}_A$. For $\beta_1\in\mathcal{B}$ and $r>0$, let $\mathbf{B}_A(\beta_1;r):=\{\beta\in\mathcal{B}: d_A(\beta,\beta_1)<r\}=\bigcap_{\alpha\in A}\mathbf{B}(\beta_1;r)$. Let $\mathcal{B}_A^+:=\{\bigvee_{\alpha\in A}T_\alpha^+<T\}$, that is, the set of $\beta\in\mathcal{B}$ that are not contained in $\bigcup_{\alpha\in A}H(\alpha)$.

Lemma 7.9. *$(\mathcal{B}_A^+,d_A)$ is separable.*

Proof. For $r\in\mathbb{Q}_{>0}$, let $\mathcal{C}_r$ denote the set of continuous curves $\gamma:[0,r]\to\Omega\cup\mathbb{R}$ with $\gamma(0)=0$ and $\gamma(t)\in\Omega$ for $t\in(0,r]$. Then $\mathcal{C}_r$ is a subset of $C([0,r],\mathbb{C})$. Let $d_r$ be the restriction of $\|\cdot\|_r$ to $\mathcal{C}_r$, that is, $d_r(\gamma_1,\gamma_2)=\sup\{|\gamma_1(t)-\gamma_2(t)|:0\le t\le r\}$. Then $(\mathcal{C}_r,d_r)$ is a subspace of $(C([0,r],\mathbb{C}),\|\cdot\|_r)$, so is separable. Let $\{\gamma_{r,n}:n\in\mathbb{N}\}$ be dense in $(\mathcal{C}_r,d_r)$. For each $r\in\mathbb{Q}_{>0}$ and $n\in\mathbb{N}$, we choose $\beta_{r,n}\in\mathcal{B}$ such that $T(\beta_{r,n})>r$ and $\beta_{r,n}(t)=\gamma_{r,n}(t)$ for $0\le t\le r$. Then $\{\beta_{r,n}:r\in\mathbb{Q}_{>0},n\in\mathbb{N}\}$ is countable.

Suppose $\beta_1\in\mathcal{B}_A^+$, and $d_0>0$. There is $r_0\in\mathbb{Q}_{>0}$ such that $\bigvee_{\alpha\in A}T_\alpha^+(\beta_1)<r_0<T(\beta_1)$. For each $\alpha\in A$, there is $t_\alpha\in(0,r_0)$ such that $\beta_1(t_\alpha)\notin H(\alpha)$. Let $d_1=\bigwedge_{\alpha\in A}\operatorname{dist}(\beta_1(t_\alpha),H(\alpha))\wedge d_0>0$. From the denseness of $\{\gamma_{r_0,n}:n\in\mathbb{N}\}$ in $(\mathcal{C}_{r_0},d_{r_0})$, we have $n_0\in\mathbb{N}$ such that $|\beta_{r_0,n_0}(t)-\beta_1(t)|=|\gamma_{r_0,n_0}(t)-\beta_1(t)|<d_1$ for $0\le t\le r_0$. Fix $\alpha\in A$. Since $|\beta_{r_0,n_0}(t_\alpha)-\beta_1(t_\alpha)|<d_1\le\operatorname{dist}(\beta_1(t_\alpha),H(\alpha))$, so $\beta_{r_0,n_0}(t_\alpha)\notin H(\alpha)$. Thus $T_\alpha(\beta_{r_0,n_0})<r_0<T(\beta_{r_0,n_0})$. Since this is true for any $\alpha\in A$, so $\beta_{r_0,n_0}\in\mathcal{B}_A^+$. Since

$$d_\alpha(\beta_{r_0,n_0},\beta_1)\le d_\alpha^\vee(\beta_{r_0,n_0},\beta_1)\le\sup\{|\beta_{r_0,n_0}(t)-\beta_1(t)|:0\le t\le r_0\}<d_1\le d_0$$

for any $\alpha\in A$, so $d_A(\beta_{r_0,n_0},\beta_1)<d_0$. Thus $\{\beta_{r,n}\}\cap\mathcal{B}_A^+$ is dense in $(\mathcal{B}_A^+,d_A)$. □

Theorem 7.2. *$\mu_\delta\to\mu_0$ weakly w.r.t. $d_A$, as $\delta\to 0$.*

Proof. Suppose $A=\{\alpha_1,\dots,\alpha_n\}$. The case $n=1$ follows from Theorem 7.1. Now suppose $n\ge 2$. We suffice to show that for any $G\in\mathcal{T}_A$, $\liminf_{\delta\to 0}\mu_\delta(G)\ge\mu_0(G)$.

We may find polygonal paths $\alpha_j^0\in\mathcal{A}$, $1\le j\le n$, such that $\alpha_j^0\succ\alpha_j$ for each $j$, and such that for $j\ne k$, any line segment on $\alpha_j^0$ is not parallel to any line segment on $\alpha_k^0$. Fix $j\in\mathbb{Z}_{[1,n]}$. List the vertices on $\overline{\alpha_j^0}$ in the counterclockwise order as $z_0^0,z_1^0,\dots,z_m^0$. We may find $z_0^1>0>z_m^1$, and $z_k^1\in\Omega$, $1\le k\le m-1$, and let $\alpha_j^1=\bigcup_{k=1}^{m-1}(z_{k-1}^1,z_k^1]\cup(z_{m-1}^1,z_m^1)$, such that $\mathcal{A}\ni\alpha_j^1\succ\alpha_j^0$, $[z_{k-1}^1,z_k^1]$ is parallel to $[z_{k-1}^0,z_k^0]$ for $1\le k\le m$, and $[z_l^0,z_l^1]\cap[z_k^0,z_k^1]=\varnothing$ for $1\le l<k\le m$. For $r\in[0,1]$, let $z_k(r)=z_k^0+r(z_k^1-z_k^0)$, $0\le k\le m$, and let $\alpha_j(r)=\bigcup_{k=1}^{m-1}(z_{k-1}(r),z_k(r)]\cup(z_{m-1}(r),z_m(r))$. Then $\alpha_j(r)\in\mathcal{A}$ for all $r\in[0,1]$, and $\alpha_j(s)\prec\alpha_j(r)$ if $0\le s<r\le 1$. And for any $s\in[0,1)$, if $\alpha_j(s)\prec\alpha\in\mathcal{A}$, then there is $r \in (s,1)$ such that $\alpha_j(r) \prec \alpha$. Thus for any $\beta \in \mathcal{B}$, we have that $r \mapsto T_{\alpha_j(r)}(\beta)$ is increasing on $[0,1]$, and for any $s \in [0,1)$, $T^+_{\alpha_j(s)} = \lim_{r\downarrow s} T_{\alpha_j(r)}$, so there are at most countably many $r \in [0,1]$ such that $T^+_{\alpha_j(r)}(\beta) > T_{\alpha_j(r)}(\beta)$. So there is $r_j \in (0,1)$ such that $\mu_0(\{T^+_{\alpha_j(r_j)} > T_{\alpha_j(r_j)}\}) = 0$. For $j = 1,\dots,k$, let $\alpha^2_j = \alpha_j(r_j)$, then $\alpha_j \prec \alpha^2_j$, and $\mu_0(\{T^+_{\alpha^2_j} > T_{\alpha^2_j}\}) = 0$.

Suppose $j \ne k \in \mathbb{Z}_{[1,n]}$. Since any line segment on $\alpha^2_j$ is not parallel to any line segment on $\alpha^2_k$, so $S_{j,k} := \alpha^2_j \cap \alpha^2_k$ is a finite set. If for some $j \ne k$ and $\beta \in \mathcal{B}$, $T_{\alpha^2_j}(\beta) = T_{\alpha^2_k}(\beta) < T(\beta)$, then $\beta$ must pass through $S_{j,k}$. From Theorem 3.1(ii), we have $T_{\alpha^2_j}(\beta_0), T_{\alpha^2_k}(\beta_0) < T(\beta_0)$. Thus $\{T_{\alpha^2_j}(\beta_0) = T_{\alpha^2_k}(\beta_0)\} \subset \{\beta_0 \text{ passes through } S_{j,k}\}$. From the property of chordal $\mathrm{SLE}_2$, for any $z_0 \in \Omega$, the probability that $\beta_0$ passes through $z_0$ is 0, which implies $\mathbf{P}[\beta_0 \text{ passes through } S_{j,k}] = 0$, so $\mu_0(\{T_{\alpha^2_j} = T_{\alpha^2_k}\}) = 0$.

For $j \in \mathbb{Z}_{[1,n]}$, let $I_j = \mathbb{Z}_{[1,n]} \setminus \{j\}$ and $B_j = \{\bigvee_{k\in I_j} T^+_{\alpha^2_k} < T_{\alpha^2_j}\} = \bigcap_{k\in I_j}\{T^+_{\alpha^2_k} < T_{\alpha^2_j}\}$, which belongs to $\mathcal{T}_{\alpha^2_j}$ from Corollary 7.4. Then $B_1,\dots,B_n$ are mutually disjoint. Let $N = \mathcal{B} \setminus \bigcup_{j=1}^m B_j$. Then

$$N \subset \bigcup_{1\le j\le n} \{T^+_{\alpha^2_j} > T_{\alpha^2_j}\} \cup \bigcup_{1\le j<k\le n} \{T_{\alpha^2_j} = T_{\alpha^2_k}\}.$$

Thus $\mu_0(N) = 0$. Fix $j \in \mathbb{Z}_{[1,n]}$. If $B \in \mathcal{T}_{\alpha_j}$, then $B \in \mathcal{T}_{\alpha^2_j}$, so $B \cap B_j \in \mathcal{T}_{\alpha^2_j}$. If $B \in \mathcal{T}_{\alpha_k}$ for some $k \in I_j$, then $B \in \mathcal{T}^+_{\alpha^2_k}$. From Lemma 7.8, we have $B \cap \{T^+_{\alpha^2_k} < T_{\alpha^2_j}\} \in \mathcal{T}_{\alpha^2_j}$. Thus $B \cap B_j = B \cap \{T^+_{\alpha^2_k} < T_{\alpha^2_j}\} \cap B_j \in \mathcal{T}_{\alpha^2_j}$. Let $\mathcal{T}_j$ denote the collection of sets $B \subset \mathcal{B}$ such that $B \cap B_j \in \mathcal{T}_{\alpha^2_j}$. Then $\mathcal{T}_j$ is a topology. We have proved that $\mathcal{T}_{\alpha_k} \subset \mathcal{T}_j$ for any $k \in \mathbb{Z}_{[1,n]}$. Thus $\mathcal{T}_A = \bigvee_{k=1}^n \mathcal{T}_{\alpha_j} \subset \mathcal{T}_j$.

Suppose $G \in \mathcal{T}_A$. Let $G_j = G \cap B_j$, $1 \le j \le n$. For each $j \in \mathbb{Z}_{[1,n]}$, since $G \in \mathcal{T}_A \subset \mathcal{T}_j$, so $G_j = G \cap B_j \in \mathcal{T}_{\alpha^2_j}$. Since $\mu_\delta \to \mu_0$ w.r.t. $d_{\alpha^2_j}$, so $\liminf_{\delta\downarrow 0} \mu_\delta(G_j) \ge \mu_0(G_j)$. Since $G$ is the disjoint union of $G \cap N$ and $G_j$, $1 \le j \le n$, and $\mu_0(G \cap N) = 0$, so

$$\liminf_{\delta\to 0} \mu_\delta(G) \ge \sum_{j=1}^n \liminf_{\delta\to 0} \mu_\delta(G_j) \ge \sum_{j=1}^n \mu_0(G_j) = \mu_0(G).$$

Since this is true for any $G \in \mathcal{T}_A$, so we have $\mu_\delta \to \mu_0$ weakly w.r.t. $d_A$, as $\delta \to 0$. □

We may find a sequence $\{\breve{\alpha}_n : n \in \mathbb{N}\}$ in $\mathcal{A}$ such that for any $\alpha \in \mathcal{A}$, there is $n \in \mathbb{N}$ such that $\breve{\alpha}_n \succ \alpha$. For $n \in \mathbb{N}$, let $T_n = \bigvee_{j=1}^n T_{\breve{\alpha}_j}$. Then for any $\beta \in \mathcal{B}$, $\bigvee_{n=1}^\infty T_n(\beta) = \bigvee_{\alpha\in\mathcal{A}} T_\alpha(\beta)$. If $\beta_0$ does not visit $f(\infty)$, then $\bigvee_{n=1}^\infty T_n(\beta_0) = T(\beta_0) = T_0$. From the property of chordal $\mathrm{SLE}_2$, $\beta_0$ does not visit $f(\infty)$ a.s., so $\bigvee_{n=1}^\infty T_n(\beta_0) = T_0$ a.s.

THEOREM 7.3. *For any $n \in \mathbb{N}$ and $\varepsilon > 0$, there is $\delta_0 > 0$ such that if $\delta < \delta_0$, then there is a coupling of $\beta_\delta$ and $\beta_0$ such that with probability greater than $1-\varepsilon$, $|\beta_\delta(t) - \beta_0(t)| < \varepsilon$ for $t \in [0, T_n(\beta_0)]$.*

PROOF. For each $1 \le j \le n$, choose $\alpha_j \succ \breve{\alpha}_j$. Let $A = \{\alpha_1, \dots, \alpha_n\}$. From Theorem 3.1(ii), we have $\beta_0 \in \mathcal{B}_A^+$. As $\delta \to 0$, $w_e^\delta \to z_e$, so $f(w_e^\delta) \to p \notin \bigcup_{\alpha \in A} H(\alpha)$. There is $\delta_1 > 0$, such that if $\delta < \delta_1$, then $f(w_e^\delta) \notin \bigcup_{\alpha\in A} H(\alpha)$, so $\beta_\delta \in \mathcal{B}_A^+$. Thus $\mu_0$ and $\mu_\delta$ are supported by $\mathcal{B}_A^+$ when $\delta < \delta_1$. From Theorem 7.2, $\mu_\delta \to \mu_0$ weakly as $\delta \to 0$, w.r.t. $d_A$. From Lemma 7.9, $(\mathcal{B}_A^+, d_A)$ is separable. So from the coupling theorem in [3], there is $\delta_0 \in (0, \delta_1)$ such that if $\delta < \delta_0$, there is a coupling of $\beta_\delta$ and $\beta_0$ such that

$$\mathbf{P}\left[d_A(\beta_\delta, \beta_0) < \bigwedge_{j=1}^n \operatorname{dist}(\alpha_j, \breve{\alpha}_j) \wedge 1 \wedge \varepsilon\right] > 1 - \varepsilon. \tag{7.7}$$

Assume $d_A(\beta_\delta, \beta_0) < \bigwedge_{j=1}^n \operatorname{dist}(\alpha_j, \breve{\alpha}_j) \wedge 1 \wedge \varepsilon$. Then for each $j \in \{1, \dots, n\}$, we have $d_{\alpha_j}(\beta_\delta, \beta_0) < \operatorname{dist}(\alpha_j, \breve{\alpha}_j) \wedge 1 \wedge \varepsilon$, which implies $d^\vee_{\breve{\alpha}_j}(\beta_\delta, \beta_0) < 1 \wedge \varepsilon$ from Lemma 7.6, so $|\beta_\delta(t) - \beta_0(t)| < 1 \wedge \varepsilon$ for $0 \le t \le T_{\breve{\alpha}_j}(\beta_\delta) \vee T_{\breve{\alpha}_j}(\beta_0)$. Since $T_n(\beta_0) = \bigvee_{j=1}^n T_{\breve{\alpha}_j}(\beta_0)$, so $|\beta_\delta(t) - \beta_0(t)| < \varepsilon$ for $t \in [0, T_n(\beta_0)]$. $\square$

THEOREM 7.4. (i) *For any $\alpha \in \mathcal{A}$, $n \in \mathbb{N}$ and $\varepsilon > 0$, there is $\delta_0 > 0$ such that if $\delta < \delta_0$, then there is a coupling of $\beta_\delta$ and $\beta_0$ such that with probability greater than $1 - \varepsilon$, $|f^{-1}(\beta_\delta(t)) - f^{-1}(\beta_0(t))| < \varepsilon$ for $t \in [T_\alpha(\beta_0), T_n(\beta_0)]$.*

(ii) *Suppose $0_+$ is degenerate. Then for any $n \in \mathbb{N}$ and $\varepsilon > 0$, there is $\delta_0 > 0$ such that if $\delta < \delta_0$, then there is a coupling of $\beta_\delta$ and $\beta_0$ such that with probability greater than $1-\varepsilon$, $|f^{-1}(\beta_\delta(t)) - f^{-1}(\beta_0(t))| < \varepsilon$ for $t \in (0, T_n(\beta_0)]$.*

PROOF. (i) Since $\beta_0([T_\alpha(\beta_0), T_n(\beta_0)])$ is a compact subset of $\Omega \setminus \{f(\infty)\}$, on which $f^{-1}$ is continuous in Euclidean metric, so there is $\varepsilon_0 > 0$ such that $\mathbf{P}[\mathcal{E}_1] > 1 - \varepsilon/2$, where $\mathcal{E}_1$ is the event that $|f^{-1}(z_2) - f^{-1}(z_1)| < \varepsilon$ for any $z_1 \in \beta_0([T_\alpha(\beta_0), T_n(\beta_0)])$ and $z_2 \in \Omega$ with $|z_2 - z_1| < \varepsilon_0$. From Theorem 7.3 there is $\delta_0 > 0$ such that if $\delta < \delta_0$, then $\beta_\delta$ and $\beta_0$ can be coupled such that with probability greater than $1 - \varepsilon/2$, $|\beta_\delta(t) - \beta_0(t)| < \varepsilon_0$ for $t \in [0, T_n(\beta_0)]$. Let $\mathcal{E}_2^\delta$ denote this event. Let $\mathcal{E}^\delta = \mathcal{E}_1 \cap \mathcal{E}_2^\delta$. Suppose $\delta < \delta_0$. Then $\mathbf{P}[\mathcal{E}^\delta] > 1 - \varepsilon$. Assume $\mathcal{E}^\delta$ occurs. Then for $t \in [T_\alpha(\beta_0), T_n(\beta_0)]$, $|\beta_\delta(t) - \beta_0(t)| < \varepsilon_0$, so $|f^{-1}(\beta_\delta(t)) - f^{-1}(\beta_0(t))| < \varepsilon$.

(ii) Suppose $0_+$ is degenerate. From [13], $f^{-1}$ extends continuously to $\Omega \cup \{0\}$. Since $\beta_0([0, T_n(\beta_0)])$ is a compact subset of $(\Omega \setminus \{f(\infty)\}) \cup \{0\}$, so the above argument still works here. $\square$

Let $\bar\gamma_0 = f^{-1} \circ \beta_0$ and $\bar\gamma_\delta = f^{-1} \circ \beta_\delta$. Then $\bar\gamma_0$ is a time-change of $\gamma_0$, and $\bar\gamma_\delta$ is a time-change of $q_\delta$.

THEOREM 7.5. *$\lim_{t\to S_0}\gamma_0(t)=\lim_{t\to T_0}\bar\gamma_0(t)=z_e$ almost surely.*

PROOF. Let $L$ be the set of spherical subsequential limits of $\bar\gamma_0(t)$ as $t\to T_0$. We first claim that $L\cap\partial^\# D=\varnothing$ a.s. If the claim is not true, then there is $\varepsilon_0>0$ such that $\mathbf{P}[L\cap\partial^\# D\ne\varnothing]>\varepsilon_0$. Since $\bar\gamma_0([T_1(\beta_0),T_2(\beta_0)])\subset D\setminus\{z_e,\infty\}$, so for every $\varepsilon>0$ there is a compact subset $F_1$ of $D\setminus\{z_e,\infty\}$ such that $\mathbf{P}[\mathcal{E}_0]>1-\varepsilon_0/3$, where $\mathcal{E}_0$ is the event that $\bar\gamma_0([T_1(\beta_0),T_2(\beta_0)])$ intersects $F_1$. Let $F_2$ be a compact subset of $D\setminus\{z_e,\infty\}$ such that $F_1$ is contained in the interior of $F_2$. Let $d_0=\operatorname{dist}(F_1,\partial F_2)>0$. From Lemma 7.3, there are $\varepsilon_1,\delta_1>0$ such that if $\delta<\delta_1$, then the probability that $\bar\gamma_\delta$ visits $\partial^\#_{\varepsilon_1}D$ after $F_2$ is smaller than $\varepsilon_0/3$. Since $\mathbf{P}[\bar\gamma_0([T_2(\beta_0),T_0))\cap\partial^\#_{\varepsilon_1/2}D\ne\varnothing]>\varepsilon_0$ and $T_0=\bigvee_{n=1}^\infty T_n(\beta_0)$ a.s., so there is $n_0\in\mathbb{N}$ such that $\mathbf{P}[\mathcal{E}_1]>\varepsilon_0$, where $\mathcal{E}_1$ is the event that $\bar\gamma_0([T_2(\beta_0),T_{n_0}(\beta_0)])\cap\partial^\#_{\varepsilon_1/2}D\ne\varnothing$. Note that $T_1=T_{\breve\alpha_1}$. From Theorem 7.4(i), there are $\delta_2<\delta_1$ and a coupling of $\bar\gamma_{\delta_2}$ with $\bar\gamma_0$ such that with probability greater than $1-\varepsilon_0/3$, $|\bar\gamma_{\delta_2}(t)-\bar\gamma_0(t)|<(\varepsilon_1/4)\wedge d_0$ for $t\in[T_1(\beta_0),T_{n_0}(\beta_0)]$. Let $\mathcal{E}_2$ denote this event. Since $\delta_2<\delta_1$, so the probability that $\bar\gamma_{\delta_2}$ does not visit $\partial^\#_{\varepsilon_1}D$ after $F_2$ is greater than $1-\varepsilon_0/3$. Let $\mathcal{E}_3$ denote this event. Let $\mathcal{E}=\bigcap_{j=0}^3\mathcal{E}_j$. Then $\mathbf{P}[\mathcal{E}]>0$. So $\mathcal{E}$ is nonempty. Assume $\mathcal{E}$ occurs. Since $\mathcal{E}_0$ occurs, so there is $t_0\in[T_1(\beta_0),T_2(\beta_0)]$ such that $\bar\gamma_0(t_0)\in F_1$. Since $\mathcal{E}_2$ occurs, so $|\bar\gamma_{\delta_2}(t_0)-\bar\gamma_0(t_0)|<d_0$, which implies that $\bar\gamma_{\delta_2}(t_0)\in F_2$. Since $\mathcal{E}_1$ occurs, there is $t_1\in[T_2(\beta_0),T_{n_0}(\beta_0)]$ such that $\bar\gamma_0(t_1)\in\partial^\#_{\varepsilon_1/2}D$. Since $\mathcal{E}_2$ occurs, so $\operatorname{dist}^\#(\bar\gamma_{\delta_2}(t_1),\bar\gamma_0(t_1))\le 2\operatorname{dist}(\bar\gamma_{\delta_2}(t_1)-\bar\gamma_0(t_1))<\varepsilon_1/2$, which implies that $\bar\gamma_{\delta_2}(t_1)\in\partial^\#_{\varepsilon_1}D$. Since $t_0\le T_2(\beta_0)\le t_1$, so $\bar\gamma_{\delta_2}$ visits $\partial^\#_{\varepsilon_1}D$ after $F_2$, which means that $\mathcal{E}_3$ cannot occur. So we get a contradiction. Thus $L\cap\partial^\# D=\varnothing$ a.s.

Second, we claim that $\operatorname{diam}^\#(L)=0$ a.s. If the claim is not true, then from the last paragraph we have $\mathbf{P}[\operatorname{diam}^\#(L)>0,L\subset D]>0$. Then there are $z_0\in D\setminus\{\infty\}$ and $r_0,\varepsilon_0>0$ such that $\overline{\mathbf{B}(z_0;4r_0)}\subset D$ and the probability that $L\cap\mathbf{B}(z_0;r_0/2)\ne\varnothing$ and $L\setminus\mathbf{B}(z_0;4r_0)\ne\varnothing$ is greater than $\varepsilon_0$. Let $\mathcal{E}_0$ denote this event. From Corollary 7.1, there is $\varepsilon_1>0$ such that with probability greater than $1-\varepsilon_0/2$, $\bar\gamma_\delta$ does not contain a $(\overline{\mathbf{B}(z_0;r_0)},r_0,\varepsilon_1)$-quasi-loop. For $n\in\mathbb{N}$, let $\mathcal{E}_0^n$ denote the event that there are $t_1<t_0<t_2<T_n(\beta_0)$ with $\bar\gamma_0(t_1),\bar\gamma_0(t_2)\in\mathbf{B}(z_0;r_0/2)$, $|\bar\gamma_0(t_1)-\bar\gamma_0(t_2)|<\varepsilon_1/3$, and $\bar\gamma_0(t_0)\notin\mathbf{B}(z_0;3r_0)$. If $\mathcal{E}_0$ occurs, then since $T_0=\bigvee_{n=1}^\infty T_n(\beta_0)$ a.s., and $\beta_0(t)$ has subsequential limits, as $t\to T_0$, inside $\mathbf{B}(z_0;r_0/2)$ and outside $\mathbf{B}(z_0;4r_0)$, so some $\mathcal{E}_0^n$, $n\in\mathbb{N}$, must occur. Thus $\mathcal{E}_0\subset\bigcup_{n=1}^\infty\mathcal{E}_0^n$. Since $\mathbf{P}[\mathcal{E}_0]>\varepsilon_0$, and $(\mathcal{E}_n^0)$ is increasing, so there is $n_0\in\mathbb{N}$ such that $\mathbf{P}[\mathcal{E}_0^{n_0}]>\varepsilon_0$. Choose $\alpha\in\mathcal{A}$ such that $f^{-1}(H(\alpha))\cap\overline{\mathbf{B}(z_0;4r_0)}=\varnothing$. From Theorem 7.4, there are $\delta_0>0$ and a coupling of $\bar\gamma_{\delta_0}$ and $\bar\gamma_0$ such that with probability greater than $1-\varepsilon_0/2$, $|\bar\gamma_{\delta_0}(t)-\bar\gamma_0(t)|<(r_0/2)\wedge(\varepsilon_1/3)$ for $t\in[T_\alpha(\beta_0),T_{n_0}(\beta_0)]$. Let $\mathcal{E}_1$ denote this event. Let $\mathcal{E}_2$ denote the event that $\bar\gamma_{\delta_0}$ does not contain a $(\overline{\mathbf{B}(z_0;r_0)},r_0,\varepsilon_1)$-quasi-loop. Then $\mathbf{P}[\mathcal{E}_2]>1-\varepsilon_0/2$ from the choice of $\varepsilon_1$. Let $\mathcal{E}=\mathcal{E}_0^{n_0}\cap$ $\mathcal{E}_1 \cap \mathcal{E}_2$. Then $\mathbf{P}[\mathcal{E}] > 0$. So $\mathcal{E}$ is nonempty. Assume $\mathcal{E}$ occurs. Since $\mathcal{E}_0^{n_0}$ occurs, so there are $t_1 < t_0 < t_2 < T_{n_0}(\beta_0)$ with $\bar{\gamma}_0(t_1), \bar{\gamma}_0(t_2) \in \mathbf{B}(z_0; r_0/2)$, $|\bar{\gamma}_0(t_1) - \bar{\gamma}_0(t_2)| < \varepsilon_1/3$, and $\bar{\gamma}_0(t_0) \notin \mathbf{B}(z_0; 3r_0)$. For $j = 1, 2$, since $\bar{\gamma}_0(t_j) \in \mathbf{B}(z_0; r_0/2)$, so $\beta_0(t_j) \not\subset H(\alpha)$, which implies that $t_j \geq T_\alpha(\beta_0)$. Since $\mathcal{E}_1$ occurs, so $|\bar{\gamma}_{\delta_0}(t_j) - \bar{\gamma}_0(t_j)| < (r_0/2) \wedge (\varepsilon_1/3)$, $j = 1, 2$, and $|\bar{\gamma}_{\delta_0}(t_0) - \bar{\gamma}_0(t_0)| < r_0$, which implies that $\bar{\gamma}_{\delta_0}(t_1) \in \mathbf{B}(z_0; r_0)$, $|\bar{\gamma}_{\delta_0}(t_1) - \bar{\gamma}_{\delta_0}(t_2)| < \varepsilon_1$ and $\bar{\gamma}_{\delta_0}(t_0) \notin \mathbf{B}(z_0; 2r_0)$, so $|\bar{\gamma}_{\delta_0}(t_0) - \bar{\gamma}_{\delta_0}(t_1)| \geq r_0$. So we find a $(\overline{\mathbf{B}(z_0; r_0)}, r_0, \varepsilon_1)$-quasi-loop on $\bar{\gamma}_{\delta_0}$, which contradicts $\mathcal{E}_2$. So $\mathbf{P}[\operatorname{diam}^{\#}(L) > 0] = 0$.

Thus almost surely $L$ is a single point in $D$, which means that $\lim_{t \to T_0} \bar{\gamma}_0(t)$ exists in the spherical metric and lies in $D$. Now we claim that $\lim_{t \to T_0} \bar{\gamma}_0(t) \notin \bar{\gamma}_0([0, T_0))$ a.s. If the claim is not true, then there exist $z_0 \in D$ and $r_0 > 0$ such that with a positive probability, we have $\lim_{t \to T_0} \bar{\gamma}_0(t) \in \bar{\gamma}_0([0, T_0)) \cap \mathbf{B}(z_0; r_0/2)$ and $\bar{\gamma}_0([0, T_0))) \not\subset \mathbf{B}(z_0; 4r_0)$, so we can use an argument that is similar to the last paragraph to find a contradiction. Note that almost surely $\bar{\gamma}_0$ does not visit $\infty$. Thus almost surely we may extend $\bar{\gamma}_0$ to be a simple continuous curve defined on $[0, T_0]$ such that $\bar{\gamma}_0(T_0) \in D \setminus \{\infty\}$. If $\mathbf{P}[\bar{\gamma}_0(T_0) \neq z_e] > 0$, then there is $n_0 \in \mathbb{N}$ such that the probability that $\bar{\gamma}_0([0, T_0])$ is enclosed by $f^{-1}(\breve{\alpha}_{n_0})$ is positive, which contradicts Theorem 3.1(ii). Thus $\mathbf{P}[\bar{\gamma}_0(T_0) = z_e] = 1$. Since $\gamma_0$ is a time-change of $\bar{\gamma}_0$, so $\lim_{t \to S_0} \gamma_0(t) = \lim_{t \to T_0} \bar{\gamma}_0(t) = z_e$ a.s. □

Proof of Theorem 4.2. (i) Choose $r > 0$ such that $\mathbf{B} := \overline{\mathbf{B}(z_e; r)} \subset D$. From Corollary 7.1, there is $\varepsilon_0 \in (0, \varepsilon)$ such that the probability that $\bar{\gamma}_\delta$ does not contain a $(\mathbf{B}, \varepsilon/6, \varepsilon_0)$-quasi-loop is greater than $1 - \varepsilon/3$. Let $\mathcal{E}_0^\delta$ denote this event. There is $\delta_1$ such that if $\delta < \delta_1$, then $|w_e^\delta - z_e| < r \wedge (\varepsilon_0/3)$. From Theorem 7.5, we have $\lim_{t \to T_0} \bar{\gamma}_0(t) = z_e$ a.s. Since $T_0 = \bigvee_{n=1}^\infty T_n(\beta_0)$ a.s., so there is $n_0 \in \mathbb{N}$ such that with probability greater than $1 - \varepsilon/3$, the diameter of $\bar{\gamma}_0([T_{n_0}(\beta_0), T_0))$ is less than $\varepsilon_0/3$. Let $\mathcal{E}_1$ denote this event. Choose $\alpha \in \mathcal{A}$ such that $f^{-1}(H(\alpha)) \subset U$. Then $T_\alpha(\beta_0) \leq T_U(\bar{\gamma}_0)$. From Theorem 7.4(i), there is $\delta_0 < \delta_1$ such that if $\delta < \delta_0$, then there is a coupling of $\bar{\gamma}_\delta$ and $\bar{\gamma}_0$ such that with probability greater than $1 - \varepsilon/3$, $|\bar{\gamma}_{\delta_2}(t) - \bar{\gamma}_0(t)| < \varepsilon_0/3$ for $t \in [T_U(\bar{\gamma}_0), T_{n_0}(\beta_0)]$. Let $\mathcal{E}_2$ denote this event. Let $\mathcal{E}^\delta = \mathcal{E}_0^\delta \cap \mathcal{E}_1 \cap \mathcal{E}_2$. Suppose $\delta < \delta_0$. Then $\mathbf{P}[\mathcal{E}^\delta] > 1 - \varepsilon$. Assume $\mathcal{E}^\delta$ occurs. Let $T_e = T_{n_0}(\beta_0)$. Then $|\bar{\gamma}_\delta(t) - \bar{\gamma}_0(t)| < \varepsilon_0/3 < \varepsilon/3$ for $T_U(\bar{\gamma}_0) \leq t \leq T_e$. And $|\bar{\gamma}_\delta(T_e) - w_e^\delta| \leq |\bar{\gamma}_\delta(T_e) - \bar{\gamma}_0(T_e)| + |\bar{\gamma}_0(T_e) - z_e| + |z_e - w_e^\delta| < \varepsilon_0$. Since $\bar{\gamma}_\delta(T_\delta) = w_e^\delta \in \mathbf{B}$ and $\bar{\gamma}_\delta$ does not contain a $(\mathbf{B}, \varepsilon/6, \varepsilon_0)$-quasi-loop, so the diameter of $\bar{\gamma}_\delta([T_e, T_\delta))$ is less than $\varepsilon/3$. Choose $\mathring{u}$ that maps $[T_U(\bar{\gamma}_0), T_\delta)$ onto $[T_U(\bar{\gamma}_0), T_0)$ such that $\mathring{u}(t) = t$ for $T_U(\bar{\gamma}_0) \leq t \leq T_e$; then $|\bar{\gamma}_\delta(\mathring{u}^{-1}(t)) - \bar{\gamma}_0(t)| < \varepsilon$ for $T_U(\bar{\gamma}_0) \leq t < T_0$. Since $\bar{\gamma}_\delta$ and $\bar{\gamma}_0$ are time-changes of $q_\delta$ and $\gamma_0$, respectively, so the proof of (i) is finished.

(ii) If $0_+$ is degenerate, then we use Theorem 7.4(ii) in the above proof. □

## 8. Other kinds of targets.

8.1. *When the target is a prime end.* Now we consider the case that the target is a prime end. We use the notation and boundary conditions given in Section 4.2 for the discrete LERW aimed at a prime end $w_e$. Suppose $f$ maps $D$ conformally onto an almost $\mathbb{H}$ domain $\Omega$ such that $f(0_+)=0$.

We will go through the propositions in Sections 6 and 7, and explain how they can be modified to prove Theorem 4.2 in this case. We only consider $D^\delta$ for $\delta \in \mathcal{M}$, so the words "$\delta < *$" should be replaced by "$\delta \in \mathcal{M}$ and $\delta < *$," and the words "$\delta \to 0$" should be replaced by "$\delta \to 0$ along $\mathcal{M}$."

Let $X_t^\xi$ and $P^\xi(t,x,\cdot)$ be notation in the case that the target is a prime end defined in Sections 3.4 and 4.1. Then all lemmas in Section 6.1 still hold. For Proposition 6.1, redefine $P_X$ to be the generalized Poisson kernel in $D_X$ with the pole at $\mathrm{Tip}(X)$, normalized by $\partial_{\mathbf{n}} P_X(w_e)=1$; let $h_X$ be defined on $V(D^\delta)$ that satisfies $h_X \equiv 0$ on $V_\partial(D^\delta) \cup \mathrm{Set}(X) \setminus \{\mathrm{Tip}(X)\}$, $\Delta_{D^\delta} h_X \equiv 0$ on $V_I(D^\delta) \setminus \mathrm{Set}(X)$, and $\Delta_{D^\delta} h_X(w_e)=1$. Proposition 6.1 should be restated as Proposition 8.1 below, which together with Proposition 2.1 implies Proposition 6.2, and then all theorems in Section 6.2.

PROPOSITION 8.1. *For any $\varepsilon > 0$, there is $\delta_0 > 0$ such that if $\delta \in \mathcal{M}$ and $\delta < \delta_0$, then for any $X \in L^\delta$, and any $w \in V(D^\delta) \cap (D \setminus H(\rho_2))$, we have $|\delta \cdot h_X(w) - P_X(w)| < \varepsilon$.*

PROOF. Fix $z_0 \in D \setminus H(\rho_2)$ and let $w_0^\delta$ be a vertex on $D^\delta$ that is closest to $z_0$. For $\delta \in \mathcal{M}$ and $X \in L^\delta$, let $g_X^0(w) = h_X(w)/h_X(w_0^\delta)$. Then from Proposition 6.1, $g_X^0$ converges to the generalized Poisson kernel $P_X^0$ in $D_X$ with the pole at $\mathrm{Tip}(X)$, normalized by $P_X^0(z_0)=1$, uniformly on $D \setminus H(\rho_3)$ for any crosscut $\rho_3$ in $D$ such that $H(\rho_1) \subset H(\rho_3)$ and $\overline{\rho_1} \cap \overline{\rho_3} = \varnothing$. Since $\partial D$ is flat near $w_e$, and $g_X^0$ vanishes on $\partial D$ near $w_e$, so $g_X^0$ can be naturally extended to be a discrete harmonic function on $\delta\mathbb{Z}^2 \cap \mathbf{B}(w_e; r_0)$ for some $r_0 > 0$. We may also extend $P_X^0$ to be a harmonic function defined in $\mathbf{B}(w_e; r_0)$ by the Schwarz reflection principle. Then we can prove that the discrete partial derivatives of $g_X^0$ approximate the corresponding partial derivatives of $P_X^0$ locally uniformly in $\mathbf{B}(w_e; r_0)$. Especially, we have $(g_X^0(w_e^\delta) - g_X^0(w_e))/\delta \to \partial_{\mathbf{n}} P_X^0(w_e)$ as $\delta \to 0$, because $w_e^\delta$ is the unique adjacent vertex of $w_e$ in $D^\delta$. Note that $\Delta_{D^\delta} g_X^0(w_e) = g_X^0(w_e^\delta) - g_X^0(w_e)$. From the definition of $g_X^0$, we have $\Delta_{D^\delta} h_X(w_e)/(\delta \cdot h_X(w_0^\delta)) \to \partial_{\mathbf{n}} P_X^0(w_e)$ as $\delta \to 0$. Since $\Delta_{D^\delta} h_X(w_e)=1$, so $1/(\delta \cdot h_X(w_0^\delta)) \to \partial_{\mathbf{n}} P_X^0(w_e)$ as $\delta \to 0$. Thus $\delta \cdot h_X(w) = g_X(w) \cdot \delta \cdot h_X(w_0^\delta)$ converges to $P_X^0(w)/\partial_{\mathbf{n}} P_X^0(w_e) = P_X(w)$ uniformly on $D \setminus H(\rho_2)$. □

In Section 7, redefine $X_w$ to be a random walk on $D^\delta$ started from $w$, stopped when it hits $V_\partial(D^\delta)$, and $Y_w$ to be that $X_w$ conditioned to hit $V_\partial(D^\delta)$ at $w_e$. Then $q_\delta$ is the loop-erasure of $Y_\delta$. Lemma 7.2 still holds.

For the proof, we argue on $Y_w$ instead of the reversal path. Then Corollary 7.1 and Corollary 7.2 immediately follow. Let $F_D$ (resp. $F_\Omega$) be a compact subset of $D\setminus\{\infty\}$ [resp. $\Omega\setminus\{f(\infty)\}$]. Lemma 7.4 still holds. Lemma 7.3, Corollary 7.3 and Lemma 7.5 should be restated as Lemma 8.1, Corollary 8.1 and Lemma 8.2, respectively, whose proofs are similar. Then we have Theorem 7.1.

LEMMA 8.1. *Suppose $U_e$ is a neighborhood of $w_e$ in $D$. Then the probability that $Y_\delta$ or $q_\delta$ visits $(D\setminus U_e)\cap\partial^\#_\varepsilon D$ after visiting $F_D$ tends to 0 as $\varepsilon\to 0$ and $\delta\to 0$ along $\mathcal{M}$.*

COROLLARY 8.1. *Suppose $U_e$ is a neighborhood of $f(w_e)$ in $\Omega$. Then the probability that $\beta_\delta$ visits $(\Omega\setminus U_e)\cap\partial^\#_\varepsilon\Omega$ after visiting $F_\Omega$ tends to 0 as $\varepsilon\to 0$ and $\delta\to 0$ along $\mathcal{M}$.*

LEMMA 8.2. *Suppose $U_e$ is a neighborhood of $f(w_e)$ in $\Omega$. Let $T^\delta_{F_\Omega}$ (resp. $T^\delta_e$) be the first time $\beta_\delta$ hits $F_\Omega$ (resp. $\overline{U_e}$). For any $\varepsilon>0$, there are $\varepsilon_0,\delta_0>0$ such that for $\delta<\delta_0$, with probability greater than $1-\varepsilon$, $\beta_\delta$ satisfies that if $|\beta_\delta(t_1)-\beta_\delta(t_2)|<\varepsilon_0$ for some $t_1,t_2\in[T^\delta_{F_\Omega},T^\delta_e]$, then $\operatorname{diam}(\beta_\delta([t_1,t_2]))<\varepsilon$.*

In Section 7.2, keep $\mathcal{B}$ unchanged, but redefine $\mathcal{A}$ to be the set of crosscuts $\alpha$ in $\mathbb{H}$ such that $\alpha$ strictly encloses 0, $H(\alpha)\subset\Omega\setminus\{f(\infty)\}$, and $H(\alpha)$ is bounded away from $f(w_e)$. Then Theorems 7.2, 7.3 and 7.4 still hold. Using Lemma 8.3 below, we can prove Theorem 7.5 with $z_e$ replaced by $w_e$, and finally Theorem 4.2.

LEMMA 8.3. *For $r>0$, the probability that $q_\delta$ visits $D\setminus\mathbf{B}(w_e;r)$ after $D\cap\mathbf{B}(w_e;\varepsilon)$ tends to 0 as $\varepsilon\to 0$ and $\delta\to 0$ along $\mathcal{M}$.*

PROOF. Let $Y^r_w$ be that $X_w$ conditioned to leave $D$ through $[\delta,0]$. Let $q^r_\delta=\mathrm{LE}(Y^r_{w^\delta_e})$. Then $q^r_\delta$ has the same distribution as the reversal of $q_\delta$. Let $P_Y$ be the probability that $Y^r_{w^\delta_e}$ visits $D\cap\mathbf{B}(w_e;\varepsilon)$ after $D\setminus\mathbf{B}(w_e;r)$. We suffice to prove that $P_Y$ tends to 0 as $\varepsilon\to 0$ and $\delta\to 0$ along $\mathcal{M}$.

We may assume that $\varepsilon<r<r_e/2$, where $r_e>0$ satisfies $\mathbf{B}(w_e;r_e)\cap D=(w_e+a\mathbb{H})\cap\mathbf{B}(w_e;r_e)$ for some $a\in\{\pm 1,\pm i\}$. Let $Q(w)$ be the probability that $X_w$ leaves $D$ through $[\delta,0]$. Let $P_X$ be the probability that $X_{w^\delta_e}$ visits $D\cap\mathbf{B}(w_e;\varepsilon)$ after $D\setminus\mathbf{B}(w_e;r)$, and leaves $D$ through $[\delta,0]$. Then $P_Y=P_X/Q(w^\delta_e)$. Let $Q_r(w)$ be the probability that $X_w$ reaches $D\setminus\mathbf{B}(w_e;r)$. Then $P_X\le Q_r(w^\delta_e)\sup\{Q(w):w\in\mathbf{B}(w_e;\varepsilon)\cap D\}$. Choose $z_0\in D$ and $r_0>0$ such that $B:=\mathbf{B}(z_0;r_0)\subset D$. Let $Q_B(w)$ be the probability that $X_w$ visits $B$ before $\partial D$. Then $Q(w^\delta_e)\ge Q_B(w^\delta_e)\inf\{Q(w):w\in B\}$. Thus

$$P_Y\le\frac{Q_r(w^\delta_e)}{Q_B(w^\delta_e)}\cdot\frac{\sup\{Q(w):w\in\mathbf{B}(w_e;\varepsilon)\cap D\}}{\inf\{Q(w):w\in B\}}. \tag{8.1}$$

Let $w_0^\delta$ be a vertex of $D^\delta$ closest to $z_0$. As $\delta \to 0$ along $\mathcal{M}$, $Q(\cdot)/Q(w_0^\delta)$ converges to the generalized Poisson kernel $P$ in $D$ with the pole at $0_+$, normalized by $P(z_0) = 1$, uniformly on any subset of $\overline{D}$ that is bounded away from $0_+$. Thus

$$\sup\{Q(w) : w \in \mathbf{B}(w_e; \varepsilon) \cap D\} / \inf\{Q(w) : w \in B\} \to 0 \tag{8.2}$$

as $\varepsilon \to 0$ and $\delta \to 0$ along $\mathcal{M}$.

As $\delta \to 0$ along $\mathcal{M}$, $Q_B$ converges to $H(D \setminus B, \partial B; \cdot)$ in $D \setminus B$. Let $U = D \cap \mathbf{B}(w_e; r)$ and $\rho = \{|z - w_e| = r\} \cap D$. As $\delta \to 0$ along $\mathcal{M}$, $Q_r$ converges to $H(U, \rho; \cdot)$ in $U$. Since $\partial D$ is flat near $w_e$, so $Q_r$ and $Q_B$ extend to be a discrete harmonic function on $\delta\mathbb{Z}^2 \cap (D \cup \mathbf{B}(w_e; r))$. So the discrete partial derivatives of $Q_B$ and $Q_r$ converge to the continuous partial derivatives of $H(D \setminus B, \partial B; \cdot)$ and $H(U, \rho; \cdot)$, respectively, in $D \cup \mathbf{B}(w_e; r)$. Thus $Q_B(w_e^\delta)/\delta \to \partial_{\mathbf{n}} H(D \setminus B, \partial B; w_e)$ and $Q_r(w_e^\delta)/\delta \to \partial_{\mathbf{n}} H(U, \rho; w_e)$ as $\delta \to 0$ along $\mathcal{M}$. So we have

$$Q_r(w_e^\delta)/Q_B(w_e^\delta) \to \partial_{\mathbf{n}} H(U, \rho; w_e) / \partial_{\mathbf{n}} H(D \setminus B, \partial B; w_e) \tag{8.3}$$

as $\delta \to 0$ along $\mathcal{M}$. The conclusion follows from (8.1), (8.2) and (8.3). $\square$

In the proof of Lemma 8.3, we consider the LERW curve $q_\delta^r$, which has the same distribution as the reversal of $q_\delta$. If $\partial D$ is flat near 0, then we have the convergence of $q_\delta^r$ to a continuous LERW$(D; w_e \to 0_+)$ trace. From the conformal invariance of continuous LERW, we have the reversibility of continuous LERW.

COROLLARY 8.2. *Suppose $w_1 \neq w_2$ are two prime ends of $D$. For $j = 1, 2$, suppose $\gamma_j(t)$, $0 < t < S_j$, is an* LERW$(D; w_j \to w_{3-j})$ *trace. Then there is a random continuous decreasing function $u_r$ that maps $(0, S_1)$ onto $(0, S_2)$ such that $(\gamma_1 \circ u_r^{-1}(t), 0 < t < S_2)$ has the same distribution as $(\gamma_2(t), 0 < t < S_2)$.*

8.2. *When the target is a side arc.* Now we consider the case that the target is a side arc. We use the notation and boundary conditions given in Section 4.2 for the discrete LERW aimed at a side arc $I_e$. Let $f$ map $D$ conformally onto an almost $\mathbb{H}$ domain $\Omega$ such that $f(0_+) = 0$.

We will modify the propositions in Sections 6 and 7 to prove Theorem 4.2 in this case. Recall that if $I_e$ is not a whole side, then we only consider $D^\delta$ for $\delta \in \mathcal{M}$, so the words "$\delta < *$" should be replaced by "$\delta \in \mathcal{M}$ and $\delta < *$," and the words "$\delta \to 0$" should be replaced by "$\delta \to 0$ along $\mathcal{M}$." If $I_e$ is a whole side, we may consider $D^\delta$ for any small $\delta$. For consistency, let $\mathcal{M} = (0, \infty)$ in this case.

Let $X_t^\xi$ and $P^\xi(t, x, \cdot)$ be notation in the case that the target is a side arc defined in Sections 3.4 and 4.1. Then all lemmas in Section 6.1 still hold. For Proposition 6.1, redefine $P_X$ to be the generalized Poisson kernel in $D_X$ with the pole at $\mathrm{Tip}(X)$, normalized by $\int_{I_e} \partial_{\mathbf{n}} P_X(z)\, ds(z) = 1$; let $h_X$ be defined on $V(D^\delta)$ that satisfies $h_X \equiv 0$ on $V_\partial(D^\delta) \cup \mathrm{Set}(X) \setminus \{\mathrm{Tip}(X)\}$, $\Delta_{D^\delta} h_X \equiv 0$ on $V_I(D^\delta) \setminus \mathrm{Set}(X)$, and $\sum_{w\in I_e^\delta} \Delta_{D^\delta} h_X(w) = 1$. Then Proposition 6.1 should be restated as Proposition 8.2, which together with Proposition 2.1 implies Proposition 6.2, and then all theorems in Section 6.2.

PROPOSITION 8.2. *For any $\varepsilon > 0$, there is $\delta_0 > 0$ such that if $\delta \in \mathcal{M}$ and $\delta < \delta_0$, then for any $X \in L^\delta$, and any $w \in V(D^\delta) \cap (D \setminus H(\rho_2))$, we have $|h_X(w) - P_X(w)| < \varepsilon$.*

PROOF. Let $z_0$, $w_0^\delta$, $h_X^0$ and $P_X^0$ be as in the proof of Proposition 8.1. Then we have the convergence of $h_X^0$ to $P_X^0$. Now we suffice to prove that $\sum_{w\in I_e^\delta} \Delta_{D^\delta} h_X^0(w) \to \int_{I_e} \partial_{\mathbf{n}} P_X^0(z)\, ds(z)$ as $\delta \to 0$ along $\mathcal{M}$.

We first consider the case that $I_e$ is a whole side. Then we may choose a polygonal Jordan curve $\sigma$ in $D$ that disconnects $I_e$ from other sides of $D$, such that $\sigma$ is disjoint from $\rho_2$, and every line segment on $\sigma$ is parallel to either the $x$ or $y$ axis. Let $U(\sigma)$ denote the doubly connected domain bounded by $I_e$ and $\sigma$. Since $P_X^0$ is bounded and harmonic in $U(\sigma)$, so we have

$$\int_{I_e} \partial_{\mathbf{n}} P_X^0(z)\, ds(z) = -\int_\sigma \partial_{\mathbf{n}} P_X^0(z)\, ds(z), \tag{8.4}$$

where $\mathbf{n}$ is the inward unit normal vector on the boundary of $U(\sigma)$.

Suppose $\delta$ is smaller than the Euclidean distance from $\sigma$ to $\rho_2$ and any side of $D$. Let $G$ be the subgraph of $D^\delta$ spanned by the set of edges in $D^\delta$ that is incident to at least one vertex in $U(\sigma)$. Let $A$ be the set of vertices of $G$ on $I_e$, and let $B$ be the set of vertices of $G$ in $D \setminus U(\sigma)$. From Lemma 6.6, we have

$$\sum_{w\in I_e^\delta} \Delta_{D^\delta} h_X^0(w) = -\sum_{(w,w')\in\mathcal{P}_\sigma} (h_X^0(w) - h_X^0(w')), \tag{8.5}$$

where $\mathcal{P}_\sigma = \{(w,w') : w \in V(D^\delta) \cap U(\sigma), w' \in V_I(D^\delta) \setminus U(\sigma), w \sim w'\}$.

Since the discrete partial derivatives of $h_X^0$ converge to the corresponding partial derivatives of $P_X^0$ uniformly on $\sigma$, so as $\delta \to 0$, we have

$$\sum_{(w,w')\in\mathcal{P}_\sigma} (h_X^0(w) - h_X^0(w')) \to \int_\sigma \partial_{\mathbf{n}} P_X^0(z)\, ds(z).$$

This together with (8.4) and (8.5) finishes the proof of the first case.

The second case is that $I_e$ is not a whole side. We assume that $\partial D$ is flat near the two ends $z_e^1$ and $z_e^2$ of $I$. We may choose a polygonal cross-cut $\sigma$ in $D$ composed of line segments parallel to $x$ or $y$ axis, such that its two ends approach to $z_e^1$ and $z_e^2$, respectively, and $\sigma$ disconnects $I_e$ from $\partial D\setminus\overline{I_e}$. Since $P_X^0$ is bounded and harmonic in $H(\sigma)$, so $\int_{I_e}\partial_{\mathbf{n}}P_X^0(z)\,ds(z)=-\int_\sigma\partial_{\mathbf{n}}P_X^0(z)\,ds(z)$, where $\mathbf{n}$ is the inward unit normal vector on the boundary of $H(\sigma)$. An argument similar to the last paragraph gives

$$\sum_{w\in I_e^\delta}\Delta_{D^\delta}h_X^0(w)=-\sum_{(w,w')\in\mathcal{P}_\sigma}(h_X^0(w)-h_X^0(w')),$$

where $\mathcal{P}_\sigma=\{(w,w'):w\in V(D^\delta)\cap H(\sigma),w'\in V_I(D^\delta)\setminus H(\sigma),w\sim w'\}$. So we suffice to show that

$$\sum_{(w,w')\in\mathcal{P}_\sigma}(h_X^0(w)-h_X^0(w'))\to\int_\sigma\partial_{\mathbf{n}}P_X^0(z)\,ds(z) \tag{8.6}$$

as $\delta\to0$ along $\mathcal{M}$. To prove this, we use the flat boundary conditions at $w_e^1$ and $w_e^2$ to extend $h_X^0$ and $P_X^0$ harmonically across $\partial D$ near $w_e^1$ and $w_e^2$. Since $\overline{\sigma}$ is compact in the extended domain: $D$ unions two balls centered at $w_e^1$ and $w_e^2$, respectively, so we get the uniform convergence of the discrete partial derivatives of $h_X^0$ to the corresponding partial derivatives of $P_X^0$ on $\sigma$. Then we are done. □

In Section 7, redefine $X_w$ to be a random walk on $D^\delta$ started from $w$, stopped when it hits $V_\partial(D^\delta)$, and $Y_w$ to be that $X_w$ conditioned to hit $V_\partial(D^\delta)$ at $I_e^\delta$. Then $q_\delta$ is the loop-erasure of $Y_\delta$. Lemma 7.2 still holds, and Corollaries 7.1 and 7.2 immediately follow from this lemma. Let $F_D$ (resp. $F_\Omega$) be a compact subset of $D\setminus\{\infty\}$ [resp. $\Omega\setminus\{f(\infty)\}$]. Lemma 7.4 still holds, and Lemma 8.1, Corollary 8.1, and Lemma 8.2 hold with $w_e$ replaced by $I_e$. Using this, we can obtain Theorem 7.1.

In Section 7.2, keep $\mathcal{B}$ unchanged, but redefine $\mathcal{A}$ to be the set of crosscuts $\alpha$ in $\mathbb{H}$ that strictly encloses 0, such that $H(\alpha)\subset\Omega\setminus\{f(\infty)\}$ and $H(\alpha)$ is bounded away from $f(I_e)$. Then we have Theorems 7.2, 7.3 and 7.4. Let $\bar\gamma_\delta=f^{-1}\circ\beta_\delta$ and $\bar\gamma_0=f^{-1}\circ\beta_0$. Using Lemma 8.4 and Theorem 8.1 below, we can prove Theorem 4.2 in this case.

LEMMA 8.4. *Let $T_{F_D}^\delta$ be the first time that $\bar\gamma_\delta$ visits $F_D$. For $a>0$, let $\partial_aD=\{z\in D:\mathrm{dist}(z,\partial D)<a\}$. For any $\varepsilon\in(0,1)$, there are $\varepsilon_0,\delta_0>0$ such that if $\delta\in\mathcal{M}$ and $\delta<\delta_0$, then with probability greater than $1-\varepsilon$, if $\bar\gamma_\delta(t_0)\in\partial_{\varepsilon_0}D$ for some $t_0\ge T_{F_D}^\delta$, then $\bar\gamma_\delta(t)\in\mathbf{B}(\bar\gamma_\delta(t_0),\varepsilon)$ for $t\ge t_0$.*

PROOF. Since $\bar\gamma_\delta$ is a time-change of $q_\delta$, which is the loop-erasure of $Y_\delta$, so we suffice to prove this lemma with $Y_\delta$ replacing $\bar\gamma_\delta$. We first consider the case that $I_e$ is not a whole side. Choose $r>0$ such that $D\cap\mathbf{B}(w_e^j,3r)=(w_e^j+a^j\mathbb{H})\cap\mathbf{B}(w_e^j,3r)$ for $j=1,2$, where $a^1,a^2\in\{\pm1,\pm i\}$, and $\mathbf{B}(w_e^1;3r)$ is disjoint from $\mathbf{B}(w_e^2;3r)$. Let $B^j=\overline{\mathbf{B}(w_e^j;r)}$ and $\sigma^j=D\cap\partial B^j$, $j=1,2$.

Let $Q(w)$ be the probability that $X_w$ hits $\partial D$ at $I_e^\delta$. Let $Q_r(w)$ be the probability that $X_w$ visits $B_1 \cup B_2$ before leaving $D$. Then $Q$ and $Q_r$ converge to $H(D, I_e; \cdot)$ and $H(D \setminus (B^1 \cup B^2), \sigma^1 \cup \sigma^2; \cdot)$, respectively, uniformly on $F_D$. Thus $Q_r(w)/Q(w) \to 0$ as $r \to 0$ and $\delta \to 0$ along $\mathcal{M}$, uniformly in $w \in F_D$. Note that $Q_r(w)/Q(w)$ is the probability that $Y_w$ visits $B^1 \cup B^2$. From the Markov property of $Y$, the probability that $Y_\delta$ visits $B^1 \cup B^2$ after $F_D$ tends to 0 as $r \to 0$ and $\delta \to 0$ along $\mathcal{M}$. So we may choose $r, \delta_e > 0$ such that $\mathbf{P}[\mathcal{E}_e^\delta] < \varepsilon/3$ if $\delta \in \mathcal{M}$ and $\delta < \delta_e$, where $\mathcal{E}_e^\delta$ is the event that $Y_\delta$ visits $B^1 \cup B^2$ after $F_D$.

For $j = 1, 2$, every point on $[w_e^j - 2a^j r, w_e^j + 2a^j r]$ corresponds to a prime end of $D$. Since $w_e^1$ and $w_e^2$ are end points of $I_e$, so $I_e \cap [w_e^j - 2a^j r, w_e^j + 2a^j r] = [w_e^j, w_e^j + 2c^j a^j r]$ for some $c^j \in \{\pm 1\}$, $j = 1, 2$. For $j = 1, 2$, let $z^j = w_e^j - c^j a^j r$; then $z^j$ is the end point of $\sigma^j$ that does not lie on $I_e$. For $j = 1, 2$, choose $\theta_1^j \ne \theta_2^j \in \sigma^j$ such that $\theta_1^j$ is closer to $z^j$ than $\theta_2^j$, and let $\rho_k^j$ denote the open arc on $\sigma^j$ bounded by $z^j$ and $\theta_k^j$, $k = 1, 2$. We may find two closed simple curves $\rho_1^0$ and $\rho_2^0$ in $D$ such that for $k = 1, 2$, $\theta_k^1$ and $\theta_k^2$ are end points of $\rho_k^0$, $\rho_k^0 \cap \sigma^j = \{\theta_k^j\}$, $j = 1, 2$; $\rho_1^0 \cap \rho_2^0 = \varnothing$; and $\rho_1 := \rho_1^0 \cup \rho_1^1 \cup \rho_1^2$ disconnects $I_e$ from any side of $D$ that does not contain $I_e$, and so $\rho_1$ is a crosscut in $D$, and $H(\rho_1)$ is a neighborhood of $I_e$. Let $\rho_2 = \rho_2^0 \cup \rho_2^1 \cup \rho_2^2$. Then $\rho_2$ is also a crosscut in $D$, and $H(\rho_2) \subset H(\rho_1)$.

For $j = 1, 2$, let $\rho_3^j = \sigma^j \setminus \rho_2^j$. Let $\rho_3 = \rho_2^0 \cup \rho_3^1 \cup \rho_3^2$ and $\rho_{1.5} = \rho_1^0 \cup (w_e^1, \theta_1^1] \cup (w_e^2, \theta_1^2]$. Then $\rho_3$ and $\rho_{1.5}$ are also crosscuts in $D$, $H(\rho_3) \subset H(\rho_{1.5})$, and $d_1 := \operatorname{dist}(\rho_3, \rho_{1.5}) > 0$. From Lemma 7.1, there are $\delta_1, \varepsilon_1 > 0$ such that if $\delta \in \mathcal{M}$, $\delta < \delta_1$, and $w \in \partial_{\varepsilon_1} D$, then the probability that $X_w$ leaves $\mathbf{B}(w; (d_1/2) \wedge (\varepsilon/3))$ is less than $\varepsilon/6$. For $w \in H(\rho_3)$, if $X_w$ hits $V_\partial(D^\delta) \setminus I_e^\delta$, then $X_w$ must intersect both $\rho_3$ and $\rho_{1.5}$, so $X_w$ must leave $\mathbf{B}(w; d_1/2)$ before it hits $\partial D$. Thus if $\delta \in \mathcal{M}$, $\delta < \delta_1$ and $w \in H(\rho_3) \cap \partial_{\varepsilon_1} D$, then $Q(w) \ge 1 - \varepsilon/6 \ge 1/2$. Since $Y_w$ is $X_w$ conditioned to hit $I_e^\delta$, so the probability that $Y_w$ leaves $\mathbf{B}(w; \varepsilon/3)$ before it hits $\partial D$ is at most 2 times the probability that $X_w$ leaves $\mathbf{B}(w; \varepsilon/3)$ before it hits $\partial D$, and so is less than $\varepsilon/3$ when $\delta < \delta_1$. From the Markov property of $Y_w$, if $\delta \in \mathcal{M}$ and $\delta < \delta_1$, then with probability greater than $1 - \varepsilon/3$, $Y_\delta$ satisfies that if $Y_\delta(t_1) \in H(\rho_3) \cap \partial_{\varepsilon_1} D$, then $Y_\delta(t) \in \mathbf{B}(Y_\delta(t_1); \varepsilon)$ for $t \ge t_1$. Let $\mathcal{E}_1^\delta$ denote this event.

Let $U_e = H(\rho_2) \setminus \rho_2$. Then $U_e$ is a neighborhood of $I_e$ in $D$. From Lemma 8.1, there are $\delta_2, \varepsilon_2 > 0$ such that if $\delta \in \mathcal{M}$ and $\delta < \delta_2$, then with probability greater than $1 - \varepsilon/3$, $Y_\delta$ does not visit $\partial_{\varepsilon_2} D \setminus U_e$ after $T_{F_D}^\delta$. Let $\mathcal{E}_2^\delta$ denote this event.

Let $\delta_0 = \delta_e \wedge \delta_1 \wedge \delta_2$, $\varepsilon_0 = \varepsilon_1 \wedge \varepsilon_2$, and $\mathcal{E}^\delta = \mathcal{E}_1^\delta \cap \mathcal{E}_2^\delta \setminus \mathcal{E}_e^\delta$. Suppose $\delta \in \mathcal{M}$ and $\delta < \delta_0$. Then $\mathbf{P}[\mathcal{E}^\delta] > 1 - \varepsilon$. Assume $\mathcal{E}^\delta$ occurs. Suppose $Y_\delta(t_0) \in \partial_{\varepsilon_0} D$ for some $t_0 \ge T_{F_D}^\delta$. Since $\delta < \delta_2$ and $\mathcal{E}_2^\delta$ occurs, so $Y_\delta(t_0) \in U_e$. Since $\delta < \delta_e$ and $\mathcal{E}_e^\delta$ does not occur, so $Y_\delta(t_0) \in H(\rho_2) \setminus (B^1 \cup B^2) \subset H(\rho_3)$. Since $\delta < \delta_1$ and $\mathcal{E}_1^\delta$ occurs, and $Y_\delta(t_0) \in H(\rho_3) \cap \partial_{\varepsilon_3} D$, so $Y_\delta(t) \in \mathbf{B}(Y_\delta^\delta(t_0); \varepsilon)$ for $t \ge t_0$.

The case that $I_e$ is a whole side is easier. We may choose a Jordan curve $\rho$ in $D$ that disconnects $I_e$ from other sides of $D$. Let $U_e$ denote the domain bounded by $I_e$ and $\rho$. From the argument used in the first part of the proof, we have $\delta_1, \varepsilon_1 > 0$ such that if $\delta \in \mathcal{M}$ and $\delta < \delta_1$, then with probability greater than $1 - \varepsilon/3$, $Y_\delta$ satisfies that if $Y_\delta(t_1) \in U_e$, then $Y_\delta(t) \in \mathbf{B}(Y_\delta(t_1); \varepsilon)$ for $t \geq t_1$. Let $\mathcal{E}_1^\delta$ denote this event. From Lemma 8.1, there are $\delta_2, \varepsilon_2 > 0$ such that if $\delta \in \mathcal{M}$ and $\delta < \delta_2$, then with probability greater than $1 - \varepsilon/3$, $Y_\delta$ does not visit $\partial_{\varepsilon_2} D \setminus U_e$ after $T_{F_D}^\delta$. Let $\mathcal{E}_2^\delta$ denote this event. Let $\delta_0 = \delta_1 \wedge \delta_2$, $\varepsilon_0 = \varepsilon_1 \wedge \varepsilon_2$ and $\mathcal{E}^\delta = \mathcal{E}_1^\delta \cap \mathcal{E}_2^\delta$. Assume $\delta \in \mathcal{M}$ and $\delta < \delta_0$, then $\mathbf{P}[\mathcal{E}^\delta] > 1 - \varepsilon$. If $\mathcal{E}^\delta$ occurs and $Y_\delta(t_0) \in \partial_{\varepsilon_0} D$ for some $t_0 \geq T_{F_D}^\delta$, then $Y_\delta(t) \in \mathbf{B}(Y_\delta(t_0), \varepsilon)$ for $t \geq t_0$. $\square$

THEOREM 8.1. *Almost surely* $\lim_{t \to S_0} \gamma_0(t) = \lim_{t \to T_0} \bar{\gamma}_0(t)$ *exists and lies on* $\partial D$.

PROOF. Let $L$ be the set of subsequential limits of $\bar{\gamma}_0(t)$ as $t \to T_0$, in the spherical metric. From Lemma 7.4, Theorem 7.4, and the idea in the first paragraph of the proof of Theorem 7.5, we have $\infty \notin L$ a.s. So $L$ is the set of subsequential limits of $\bar{\gamma}_0(t)$ as $t \to T_0$, in the Euclidean metric. From Theorem 3.1(ii), we have $L \cap \partial D \neq \varnothing$ a.s. From Theorem 7.4, Lemma 8.4, and the idea in the second paragraph of the proof of Theorem 7.5, we have $\operatorname{diam}(L) = 0$ a.s. So we are done. $\square$

From the property of discrete LERW and the conformal invariance of continuous LERW, we then have the following corollary.

COROLLARY 8.3. *Suppose* $\gamma(t)$, $0 \leq t < S$, *is an* $\mathrm{LERW}(D; w_0 \to I_e)$ *trace; then almost surely* $\widehat{\lim}_{t \to S} \gamma(t)$, *the limit of* $\gamma(t)$ *in* $\widehat{D}$, *as* $t \to S$, *exists and lies on* $I_e$, *and the distribution of* $\widehat{\lim}_{t \to S} \gamma(t)$ *is the same as the distribution of the limit point in* $\widehat{D}$ *of the Brownian excursion in* $D$ *started from* $w_0$ *conditioned to hit* $I_e$. *And if* $J_e$ *is a subarc of* $I_e$, *then after a time-change,* $\gamma(t)$ *conditioned on the event that* $\widehat{\lim}_{t \to S} \gamma(t) \in J_e$ *has the same distribution as an* $\mathrm{LERW}(D; w_0 \to J_e)$ *trace.*

QUESTION. Can we prove Theorem 7.5, Corollary 8.2 and Corollary 8.3 directly from the definition of continuous LERW?

**Acknowledgments.** I would like to thank Professor Nikolai Makarov, who introduced me to this subject, and gave me many valuable suggestions. I also thank Professor Yuvel Peres for his interesting courses and talks about random walks.

Department of Mathematics
Yale University
P.O. Box 208283
New Haven, Connecticut 06520-8283
USA
E-mail: dapeng.zhan@yale.edu